\documentclass[12pt]{amsart}

\usepackage{amsfonts,amsmath,amssymb,amsthm}
\usepackage{hyperref}
\hypersetup{breaklinks=true}

\numberwithin{equation}{section}
\theoremstyle{plain}
\newtheorem{theorem}{Theorem}[section]    % theorem with number
\newtheorem{lemma}[theorem]{Lemma}       % lemma with number
\newtheorem{proposition}[theorem]{Proposition}      % lemma with number
\newtheorem{corollary}[theorem]{Corollary}

\newtheorem{definition}[theorem]{Definition}

\theoremstyle{remark}
\newtheorem*{remark*}{Remark}  % theorem without number
\newcommand{\acts}{\operatorname{\curvearrowright}}

\newcommand{\Aut}{\mathop{\mathrm{Aut}}\nolimits}

\newcommand{\Confdim}{\mathop{\mathrm{Confdim}}\nolimits}

\newcommand{\D}{\partial}
\newcommand{\De}{\Delta}

\newcommand{\diam}{\mathop{\mathrm{diam}}\nolimits}
\newcommand{\dist}{\mathop{\mathrm{dist}}\nolimits}
\newcommand{\eps}{\epsilon}
\newcommand{\F}{\mathcal{F}}
\newcommand{\Ga}{\Gamma}
\newcommand{\ga}{\gamma}

\newcommand{\hr}{\mathbb{R}}
\newcommand{\hh}{\mathbb{H}}
\newcommand{\hn}{\mathbb{N}}
\newcommand{\hp}{\mathbb{P}}
\newcommand{\hs}{\mathbb{S}}
\newcommand{\hm}{\mathbb{M}}

\newcommand{\Isom}{\mathop{\mathrm{Isom}}\nolimits}
\newcommand{\Mod}{\mathop{\mathrm{Mod}}\nolimits}

\newcommand{\ra}{\rightarrow}
\newcommand{\R}{\mathbb{R}}

\newcommand{\si}{\operatorname{\sigma}}

\renewcommand{\S}{\mathbb{S}}

\newcommand{\T}{\mathcal{T}}

\begin{document}

\title[Combinatorial modulus and Coxeter groups]{Combinatorial modulus, 
The Combinatorial Loewner Property, and Coxeter groups} 
\author{ Marc Bourdon and Bruce Kleiner} 
%\date{22 October 2007}
\maketitle

\begin{abstract}
We study combinatorial modulus on  self-similar metric spaces.
We give new examples
of hyperbolic groups whose boundaries satisfy a combinatorial version
of the Loewner property, and prove Cannon's conjecture for Coxeter
groups.  We also establish some connections
with  $\ell_p$-cohomology.
\end{abstract}

%\begin{divers}
%  {Tools} Quasi-conformal geometry, quasi-isometries.
%\end{divers}

\setlength{\parskip}{\smallskipamount}

\tableofcontents

\setlength{\parskip}{\medskipamount}

\section{Introduction}

\subsection{Overview} \label{overview}

Every hyperbolic group $\Ga$  has a canonical action 
on its boundary at infinity $\D \Ga$;
with respect to any visual metric on $\D \Ga$, this action
is by uniformly quasi-Moebius homeomorphisms.
 This structure has a central role in the proofs 
of  Mostow's
rigidity theorem and numerous other results in the same vein,
which are based on the analytic theory of quasiconformal homeomorphisms
of the boundary (see the survey papers \cite{GPa, BP3, K} and their references). 
With the aim of extending these rigidity results to a larger class of 
hyperbolic groups, one may hope to apply the work
of J. Heinonen and P. Koskela \cite{HeK} and subsequent authors
(\emph{e.g.} \cite{cheeger,Ty,khshanmuga,keith}), which has generalized much
of the classical quasiconformal theory to the setting of Loewner
spaces\footnote{
We will use the 
shorthand {\em Loewner space} for a compact metric space that is
 Ahlfors $p$-regular and $p$-Loewner
for some $p>1$, in the sense of \cite{HeK}.},  a certain class of metric
measure spaces \cite{HeK}.
Unfortunately, among the currently known
examples of Loewner spaces, the only ones which  arise
as boundaries of hyperbolic groups 
are the  boundaries of rank-one
symmetric spaces, and Fuchsian buildings 
\cite{HeK,BP1}.   One of the goals of this paper 
is to take a step
toward improving this situation, by finding infinitely
many new examples
of hyperbolic groups whose boundaries satisfy the
{\em Combinatorial Loewner Property}, a
combinatorial variant of the Loewner property which is conjecturally
equivalent 
to the property of being
quasi-Moebius homeomorphic to a Loewner space.  
In addition to this,
using similar techniques, we prove the Cannon conjecture for hyperbolic Coxeter groups,
show that 
the $\ell_p$-equivalence relation  studied in  
\cite{G,E,B2} has a particularly simple form in the case of 
hyperbolic Coxeter groups, and prove that the standard square 
Sierpinski carpet in the plane and the standard cubical Menger curve
in $\hr ^3$ satisfy the Combinatorial Loewner Property.

\bigskip

\subsection{Statement of results} \label{statement}

We now present some of the ideas of the paper,  illustrating them
with non-technical statements.
More general results, as well as detailed discussion justifying the statements
made here, may be found in the body of the paper.

\medskip

\subsection*{Combinatorial modulus}
Let $Z$ be a compact metric space. For every $k \in \hn $, 
let $G_k$ be the incidence graph of a
ball cover $\{B(x_i,2^{-k})\}_{i\in I}$,
where $\{x_i\}_{i\in I}\subset Z$ is a maximal $2^{-k}$-separated subset.
Given $p \ge 1$ and a curve family $\mathcal F$ 
in $Z$, we denote by $\Mod _p (\mathcal F , G_k )$ the $G_k$-combinatorial
$p$-modulus of $\mathcal F$ (see Subsection \ref{first} for the definition);
also for any pair of subsets $A,B\subset Z$, we let
$\Mod _p(A,B,G_k)=\Mod _p(\mathcal F,G_k)$ where
$\mathcal F$ is the collection of paths joining 
$A$ and $B$.

In our study
of combinatorial modulus, we will assume that $Z$ is
 \emph{approximately self-similar} (see 
Definition \ref{homothety}).   Examples of approximately self-similar metric spaces
include many classical fractals such as 
the square Sierpinski carpet or the cubical Menger sponge,  boundaries 
of hyperbolic groups equipped with their visual metrics,
metric spaces associated with finite subdivision rules, 
and metric $2$-spheres arising from expanding Thurston map.

\medskip
One of our principal goals is to find criteria for 
the {\em Combinatorial Loewner Property} (CLP).   Roughly speaking,
a  doubling space $Z$ satisfies the CLP if there is a $p\in (1,\infty)$
such that for any pair $A,B\subset Z$ of disjoint
nontrivial continua, 
the $p$-modulus $\Mod_p(A,B,G_k)$ is controlled by the relative distance
$$\Delta (A, B) = \frac {\dist (A, B) } {\min \{\diam A, \diam B \}}\;;$$
see Section \ref{comb} for the definition.
As indicated in the overview, 
our interest in the CLP stems from the fact that a Loewner space 
satisfies the CLP, and that  
the converse is conjecturally true for compact approximately
self-similar spaces.  Thus -- modulo the conjecture -- this paper would provide
new examples of hyperbolic groups to which the recently 
developed quasiconformal theory would be applicable.
Should the conjecture turn out to be false, the CLP would be 
of independent interest, since it shares many of the features of the
Loewner property, \emph{e.g.} quasi-Moebius invariance; see Section \ref{comb}
for
more discussion.

\subsection*{Coxeter group boundaries satisfying the CLP}

Recall that a group $\Gamma$ is a \emph{Coxeter group} if it admits
a presentation of the form
$$ \Gamma = \langle s_i \in S ~ \vert ~ s_i ^2 =1, (s_i s_j)^{m_{ij}} =1
~\textrm{for}~ i \neq j \rangle , $$
with $\vert S \vert < + \infty $, and with $m_{ij} \in \{2,3, ..., + \infty \} $. 
A subgroup is {\em special} if it is generated by a subset of the
generating set $S$.   A subgroup is   \emph{parabolic}
if it is conjugate to a special subgroup.
Now suppose  in addition  that $\Gamma$ 
is   hyperbolic, and  $\partial \Gamma$ is its boundary
at infinity.  Then a non-empty limit set $\partial P \subset 
\partial \Gamma$ of a parabolic subgroup $P \leqq \Gamma$ is
called a \emph{parabolic limit set}.

\medskip

In Theorem \ref{CLP} we give a sufficient condition 
(of combinatorial flavour) for the boundary of a 
hyperbolic Coxeter group to satisfy the CLP. 
A special case is the following result, which shows that the 
CLP holds when the parabolic limit sets form a 
combinatorially simple collection of subsets.

\begin{theorem} [Corollary \ref{simple}] \label{theosimple}
Let $\Gamma$ be a hyperbolic Coxeter group whose boundary is connected
and such that $\Confdim (\partial \Gamma ) > 1$. Assume that for every 
proper, connected, parabolic limit set $\partial P \subset \partial
\Gamma$, one has 
$$\Confdim (\partial P) < \Confdim (\partial \Gamma). $$
Suppose furthermore that for every pair $\partial P$, $\partial Q$ of
distinct, proper, connected, parabolic limit sets, the subset
$\partial P \cap \partial Q$ is totally disconnected or empty.
Then $\partial \Gamma$ 
satisfies the CLP.
\end{theorem}

In the statement $\Confdim (Z)$ denotes the Ahlfors regular conformal
dimension of $Z$ \emph{i.e.} the infimum of Hausdorff dimensions of Ahlfors regular
metric spaces quasi-Moebius homeomorphic to $Z$ (see \cite{MT} for more
details about conformal dimension). 
\medskip

To illustrate Theorem \ref{theosimple} with some simple examples,
consider a Coxeter group $\Gamma$  with
a Coxeter presentation $$\langle s_1,\ldots,s_4\mid s_i^2=1,\;(s_is_j)^{m_{ij}}=1
\;\mbox{for}\;i\neq j\rangle\,,$$
where
the order $m_{ij}$ is finite for all $i\neq j$,
and for every $j \in \{1, ..., 4 \}$ one has 
$\sum _{i \neq j} \frac{1}{m_{ij} } <  1$.   For these examples, the 
proper connected parabolic limit sets are circles, and hence have
conformal dimension $1$, while  $\partial\Gamma$ is homeomorphic
to the Sierpinski carpet and therefore has conformal dimension $>1$
by a result
of J. Mackay \cite{Mac}.  Theorem \ref{theosimple} therefore applies,
and $\partial\Gamma$ has the CLP.

\medskip

Applying similar techniques in a simplified setting, we prove:
\begin{theorem}[Theorem
\ref{hs}]
The square
Sierpinski carpet and the cubical Menger sponge satisfy the CLP.
\end{theorem}

\medskip

\subsection*{The Cannon conjecture for Coxeter groups} 
We obtain a 
proof of  Cannon's conjecture in the special case
of Coxeter groups:

\begin{theorem} [Theorem \ref{cannon}] \label{introcannon}
Let $\Gamma$ be a hyperbolic Coxeter group
whose boundary is homeomorphic to the $2$-sphere. Then there is
a properly discontinuous, cocompact, and isometric action of
$\Gamma$   on $\hh ^3$.
\end{theorem}
\noindent
This result was essentially
known (see the discussion at 
the end of Section \ref{cannonsection}).  Our view is that the
principal value of the proof is that it  illustrates  the feasibility of the
asymptotic approach (using the ideal boundary
and modulus estimates), and it may suggest ideas for attacking the 
general case.
It also gives a new proof of the Andreev's theorem on realizability
of polyhedra in $\hh^3$, in the case when the prescribed
dihedral angles are submultiples of $\pi$.

\medskip

\subsection*{$\ell _p $-equivalence relations} 
Let $\Gamma$ be a hyperbolic group and let $p \ge 1$.
The first $\ell _p $-cohomology group of $\Gamma$
induces on $\partial \Gamma$ an equivalence relation
 --  the \emph{$\ell _p $-equivalence relation} -- which 
is invariant under  
quasi-isometries of
 $\Gamma$ \cite{G,E,BP2,B2}. A natural problem
is to determine its cosets.  The existence of a non trivial coset,
\emph{i.e.} a coset different from a point and from the whole space,
was shown to be an obstruction to the Loewner property
in  \cite{BP2}.  Inspired by this, we prove 
the analogous statement for the CLP in  Corollary \ref{confdim}.
Moreover, in consequence of some 
of our previous results, we get:

\begin{corollary} [Corollaries \ref{coset} and \ref{Confdim}]
Assume $\Gamma$ is a hyperbolic Coxeter group,
let $p \ge 1$, and denote by $\sim _p$ the 
$\ell _p $-equivalence relation on $\partial \Gamma$.
Then:
\begin{enumerate}
\item Each coset of $\sim _p $
is either a point or a connected parabolic limit set.
\item If $\partial \Gamma$ is connected, 
and $ \sim _p$ admits a coset different from a point and the whole 
set $\partial \Gamma$,
then $ \sim _p$ admits a coset $F$ with $\Confdim (F) = \Confdim (\partial \Gamma)$.
\end{enumerate}

\end{corollary}

\smallskip

\subsection*{Beyond the CLP}
At present, our understanding of Coxeter groups is still quite limited. 
We have only been able to show that certain 
very special groups have boundaries which satisfy the
CLP.   
While this is consistent with our expectation that
the CLP should be a highly non-generic property,
we only have a few examples which 
are known not to have the CLP,  apart from groups whose boundaries
could not have the CLP for topological reasons, see 
Proposition \ref{linear} and the remark at the end of Section \ref{equivalence}.  

It would be desirable to have effective criteria for showing that
a group boundary does not have the CLP, as well as new examples
of such groups.
In addition, when a group boundary does not have the CLP, we expect
that alternative structure will be present instead, such as
a quasi-Moebius invariant equivalence relation.

\bigskip

\subsection{Discussion of the proofs} \label{presentation}
We now give an indication of the ideas that go into some of the 
proofs.

\subsection*{Combinatorial modulus on approximately self-similar spaces}
Let $Z$ be a compact, approximately self-similar space.
For $d_0 > 0$, let $\mathcal F _0 $ be the family of curves $\gamma \subset
Z $ with $\diam (\gamma) \ge d_0 $.
The approximate self-similarity of $Z$ allows one to compare moduli of
curve families at different locations and scales with the modulus of
$\mathcal{F}_0$, and this observation leads readily to a 
submultiplicativity relation between combinatorial
moduli at different scales: 

\begin{proposition} [Proposition \ref{submult}]
Let $Z$ be an arcwise connected approximately self-similar metric space.  
Let $p \ge 1$ and set $M_k := \Mod _p (\mathcal{F} _0 , G_k) $. 
Then, for $d_0$ sufficiently small, 
there exists a constant $C>0$ such that for every pair of integers
$k, \ell$ one has : 
\begin{equation}
\label{eqn-submultiplicative}
M_{k+\ell} \le C \cdot M_k \cdot M_{\ell}\,.
\end{equation}
In addition when $p$ belongs to a compact subset of $[1, +\infty)$
the constant $C$ may be chosen independent of $p$.
\end{proposition}

Iterating (\ref{eqn-submultiplicative}) yields
$M_{nk}\leq C^{n-1}M_k^n$, which implies that $\lim_{j\ra\infty}M_j= 0$
if $M_k<C^{-1}$ for some $k$.  Therefore if we define the critical exponent
to be
$$
Q_M=\inf\{p\in [1,\infty) ~;~ \lim_{k\ra\infty}\,M_k=0\}\,
$$
then $M_k=\Mod_{Q_M}(\mathcal{F}_0,G_k)\geq C^{-1}$ for all $k$.
In fact, $Q_M$ is the Ahlfors regular conformal dimension of $Z$,
\cite{KeK}.

When $Z$ is the standard Sierpinski carpet constructed from the 
unit square, one can exploit the reflectional 
symmetry to get additional control on the modulus.  Using it,
one shows
that for any nonconstant curve $\eta:[0,1]\ra Z$ and any $\eps>0$, 
if $\mathcal{U}_{\eps}(\eta)$
denotes the $\eps$ neighborhood of $\eta$ in the $C^0$ topology, then
the $p$-modulus  $\Mod_p(\mathcal{U}_{\eps}(\eta),G_k)$ is 
uniformly comparable
to $M_k$, independent of $k$.  In other words, the modulus of the 
curves near an arbitrary curve is comparable to the modulus of 
all curves.   From this, and using the planarity of the carpet,
one can prove a supermultiplicativity inequality as well:
\begin{equation}
\label{eqn-supermultiplicative}
M_{k+\ell} \ge C' \cdot M_k \cdot M_{\ell}\,,
\end{equation}
where $C'\in (0,\infty)$ may be chosen in terms of an upper bound
on $p$.  Reasoning as above, it follows that if $M_k>C'^{-1}$
for any $k$, then $\lim_{k\ra\infty}M_k=\infty$; this implies that
at the critical exponent $Q_M$, the sequence $\{M_k\}$ is bounded 
away from zero and infinity.  From this, one concludes that
the statement of the CLP holds for pairs of balls.  By 
imitating an argument from \cite{BK3}, one shows that 
the CLP holds provided it holds for pairs of balls
(Proposition \ref{ball}), and therefore
the Sierpinski carpet satisfies the CLP (see Theorem \ref{hs}).

In the case of the  Menger curve, instead of (\ref{eqn-supermultiplicative}),
one obtains an estimate for
$M_k$ in terms of the moduli $\{M_j\}_{j<k}$, which is sufficient to
verify the CLP.

\bigskip

\subsection*{Dynamics of curves and crossing}
Our strategy for understanding the boundary of a
hyperbolic Coxeter group $\Gamma$ is inspired by the
analysis of the 
Sierpinski carpet, although the story is more 
complicated.

For the purpose of this paper, much of the dynamics of the
$\Gamma$-action on $\partial \Gamma$ is encoded in the 
parabolic limit sets. For example we get:

\begin{theorem} [Corollary \ref{1.3}]  
Consider a $\Gamma$-invariant equivalence relation 
on $\partial \Gamma$ whose cosets are connected.
Then:
\begin{itemize}
\item  The closure of
each coset is either a point or a parabolic limit set.
\item  If a nontrivial coset $F$ is path-connected, and $P$ is the 
parabolic subgroup with $\overline F = \D P$, then for every $\eps>0$
and every path $\eta:[0,1]\ra\D P$, there is a path $\eta':[0,1]\ra
F$  such that
$$
d(\eta,\eta')=\max_{t\in [0,1]}d(\eta(t),\eta'(t))<\eps\,.
$$ 
\end{itemize}
\end{theorem}

As an illustration, let $\gamma$ be a nontrivial curve in
$\partial \Gamma$, and let $\sim$ be the  smallest 
equivalence relation 
on $\partial \Gamma$ such that for every $g\in \Ga$,
the  curve $g\gamma$ lies in a single coset; in other words,
two points $x,y\in \partial\Gamma$ lie in the same coset if there
is a finite chain $g_1\ga,\ldots,g_k\ga$ joining $x$
to $y$. Thanks to the previous 
theorem, the coset closures are either points
or parabolic limit sets. In particular, if $\gamma$ is not contained in 
any proper parabolic limit set, then any path $\eta:[0,1]\ra \D\Ga$
is a uniform limit of paths lying in the coset of $\gamma$. 

A key ingredient in the analysis 
of the 
combinatorial modulus on $\partial \Gamma$
is a quantitative version of this
phenomenon, which  is established in Proposition \ref{1.8}.

\bigskip

\subsection*{The proof of the Cannon conjecture for Coxeter groups} 
By \cite{S},
if $\Gamma$ is a  hyperbolic group
and  $\partial \Gamma$ is
quasi-Moebius homeomorphic to the Euclidean $2$-sphere, 
then $\Gamma$ admits a properly discontinuous, cocompact,
isometric action on $\hh ^3 $. 
Also, as a consequence of the uniformization
theorem established in \cite{BK1}, we obtain:

\begin{corollary} [Corollary \ref{corunif}] \label{uniformization}
Suppose $Z$ is an approximately self-similar metric space 
homeomorphic
to the $2$-sphere. Assume that for $d_0 >0$ small enough, there exists a constant
$C = C(d_0) \ge 1$ such that for every $k \in \hn $ one has 
\begin{equation}
\label{0.1}
\Mod _2 (\mathcal F _0 , G_k ) \le C .
\end{equation}
Then $Z$ is quasi-Moebius homeomorphic to the Euclidean $2$-sphere.
\end{corollary}
\noindent
Thus, we are reduced to verifying the hypotheses of 
the above corollary when 
$\Gamma$ is  
a Coxeter group.   We note that an alternate reduction to the same assertion
can be deduced using \cite{CFP}.

One of the main results of the paper is the existence of a finite number of
``elementary curves families'', whose moduli govern the modulus of every
(thick enough) curve family in $\partial \Gamma$ (see Theorem \ref{3.4} and Corollary 
\ref{comparable}). Each elementary curve family is associated to a conjugacy
class of an infinite parabolic subgroup.

In consequence, to obtain the bound (\ref{0.1}), it is enough to establish that every connected 
parabolic limit set $\partial P$ enjoys the following property:
there exists a non constant continuous curve $\eta \subset \partial P$,
such that letting $\mathcal U _{\epsilon} (\eta)$ be the
$\epsilon$-neighborhood of $\eta$ in the ${C}^0 $-topology, the
modulus $\Mod _2 (\mathcal U _{\epsilon} (\eta), G_k)$ is bounded
independently of $k$, for $\epsilon >0$ small enough.
 
To do so, two cases are distinguished: either $\partial P$ is a circular
limit set \emph{i.e.} it is homeomorphic to the circle, or it
is not.
 
In the second case one can find two crossing curves $\eta_1, \eta_2
\subset \partial P$. Since $\partial \Gamma$ is a planar set, one gets that 
$\min _{i=1,2} \Mod _2 (\mathcal U _{\epsilon} (\eta _i), G_k)$
is bounded independently of $k$, for $\epsilon$ small enough.
Note that crossing type arguments in relation with the 
combinatorial $2$-modulus, appear frequently in the papers
\cite{Ca}, \cite{CW} and \cite{CFP} (not to mention the whole body
of literature on $2$-dimensional quasiconformal geometry).

Let $r >0$, and denote by $\mathcal F _1$ 
the subfamily of $\mathcal F _0 $ consisting of
the curves $\gamma \in \mathcal F _0 $ which do
not belong to the $r$-neighborhood 
$N_r (\partial P)$ of any
circular limit set $\partial P$. 
At this stage one knows that for $r$ small enough,
$ \Mod _2 (\mathcal F _1 , G_k )$ is bounded independently
of $k$.  To bound the modulus of $\mathcal F_0\setminus \mathcal F_1$,
we proceed as follows.  Consider a 
curve $\ga\in  \mathcal F_0$ contained in $N_r(\partial P)$, where
$\partial P$ is a circular parabolic limit set.  The idea is to 
break $\ga$ into pieces $\ga_1,\ldots,\ga_i$, such that for each
$j\in \{1,\ldots,i\}$, the maximal distance $\max\{d(x,\partial P)
~;~  x\in \ga_j\}$
is comparable to $\diam(\ga_j)$.  Then for each $j$, applying a suitable group element
$g\in \Gamma$, we can arrange that both $g\ga_j$ and $g(\partial P)$ have
roughly unit diameter.  Since $g\,\ga_j$ lies close to $g(\partial P)$, but not
too close, it cannot lie very close to a circular limit set;
it follows that $g\ga_j$ belongs to a curve family with controlled
modulus.  We then apply $g^{-1}$ to the corresponding admissible function,
and renormalize it suitably; by summing the collection of functions which 
arise in this fashion from all such configurations, we arrive at an admissible function for
all such curves $\ga$.  The fact that the conformal dimension of $S^1$
is strictly less than $2$ allows us to bound the $2$-mass of this admissible function, and 
this yields the desired bound (\ref{0.1}).

We note that in the body of the paper, the
argument in the preceding paragraph appears in
Theorem \ref{theocircular}, where it is formulated in greater
generality.   It is also used in the proofs of Theorem \ref{CLP} and Corollary \ref{Confdim}.

\subsection{Organization of the paper}

Combinatorial modulus, Combinatorial Loewner Property, and 
approximately self-similar metric spaces are presented
in Sections \ref{modulus} and \ref{hyp}.  
In Section \ref{menger} the Combinatorial Loewner Property is established
for the square Sierpinski carpet and the cubical Menger sponge.
Section \ref{preliminary} 
contains preliminary results about various 
dynamical aspects of the action of a hyperbolic Coxeter group on its 
boundary. 
Section \ref{cox} is the heart of the 
paper, 
it focusses on the combinatorial modulus on boundaries of hyperbolic 
Coxeter groups. 
Section \ref{cannonsection} discusses a proof of the Cannon's conjecture 
in the Coxeter group case.
An application to Coxeter groups with Sierpinski carpet boundary is given.
Section \ref{ex} establishes a sufficient condition for a Coxeter boundary 
to satisfy the Combinatorial Loewner Property.
Examples are presented in Section \ref{examples}. 
Section \ref{equivalence} 
discusses applications to
$\ell _p$-equivalence relations.

\subsection{ Suggestions to the reader}

Readers who are concerned only with the Coxeter group
case of the Cannon's conjecture may read Section 
\ref{preliminary}, Subsection \ref{first},
Section \ref{hyp} until Corollary \ref{corunif},
and Sections \ref{cox} and \ref{cannonsection}.

\medskip
\subsection*{Acknowledgements} We would like to thank Mario Bonk 
and Frederic Haglund 
for several helpful discussions on the subject of 
this article. Special thanks to Peter Ha\"issinsky for his interest
and for several comments on a first version of the paper.  We thank the referee
for his/her numerous and very valuable suggestions. 
M. B. was partially supported by ANR grant ``Cannon'', and B. K.  
by NSF grant DMS-0701515.

\subsection*{Notation and conventions}  Any curve 
$\gamma : [a, b] \to X$ is
assumed to be continuous. Often we do not distinguish between $\gamma$
and its image $\gamma ([a, b])$. 
Two real valued functions $f, g$ defined on a space $X$ are said to be 
\emph{comparable},
and then we write $f \asymp g$, if there exists a constant $C>0$ such that
$C^{-1} f \le g \le C f$.  We write $f \lesssim g$ if there is a constant 
$C >0$ such that $f \le C g$.
A \emph{continuum} is a non-empty compact connected topological space, it 
is \emph{non-degenerate}
if it contains more than one point.
For an open ball $B = B(x,r)$ in a metric space and for $\lambda >0$, we 
denote by $\lambda B$
the ball $B(x, \lambda r)$. The radius of a ball $B$ is denoted by $r(B)$.
The open $r$-neighborhood of a subset $E$ is denoted by $N_r (E)$.
A subset $E$ of a set $F$ is called a \emph{proper subset of $F$} if
$E \subsetneqq F$. 
 A subset $E$ of a geodesic metric space $X$ is \emph{convex} if for every
$x, y \in E$, all the geodesic segments in $X$ between $x$ and $y$ are contained in $E$.
Given a graph $G$, we define $G^0$ and $G^1$ to be the set of vertices and of
(non-oriented) open edges respectively. The inclusion of a subgroup $H$ in a group
$\Gamma$ is denoted by $H \leqq \Gamma$.

\section{Combinatorial modulus} \label{modulus}

This section develops the theory of combinatorial modulus in a general setting.
The Combinatorial Loewner Property and related topics are discussed. 
For the classical
notions of geometric function theory used in this paper, we refer to
\cite{He,Vai,BK1}.

Versions of combinatorial modulus
have been considered 
by several authors in connection with Cannon's conjecture
on groups with $2$-sphere boundary
(see \emph{e.g.} \cite{Ca,CW,CFP,BK1,H}), and in a more general
context  \cite{Pa2,Ty}. 

\subsection{Definitions and first properties} \label{first}

Let $(Z, d)$ be a compact metric space, let $k \in \hn$, and let 
$\kappa \geq 1$. A finite graph
$G_k$ is called a \emph{$\kappa$-approximation of $Z$ on scale $k$},
 if it is
the incidence graph of a
covering of $Z$,  
such that for every $ v \in G_k ^0$ there
exists $z_v \in Z$ with
$$ B(z_v, \kappa ^{-1} 2^{-k}) \subset v \subset B(z_v, \kappa 2^{-k}), $$
and for $v, w \in G_k ^0$ with  $v \neq w $ :
$$  B(z_v, \kappa ^{-1} 2^{-k}) \cap B(z_w, \kappa ^{-1} 2^{-k}) = 
\emptyset. $$
Note that we identify every vertex $v$ of $G_k$ with the 
corresponding subset in $Z$.
A collection of graphs $\{ G_k \} _{k \in \hn}$ is called a
\emph{$\kappa$-approximation} of $Z$, if for each $k \in \hn$
the graph $G_k$ is a
$\kappa$-approximation of $Z$ on scale $k$.

\medskip

Let $\gamma \subset Z$ be a  curve (or a subset) and let $\rho :
G_k ^0 \to \hr_+ $ be any function. The \emph{$\rho$-length} of $\gamma$ is
$$ L_{\rho} (\gamma) = \sum_{v \cap \gamma \neq \emptyset} \rho (v).$$
For $p \geq 1$ the \emph{$p$-mass} of $\rho$ is
$$ M_p (\rho) = \sum_{v \in G_k ^0} \rho (v)^p.$$
Let $\mathcal{F}$ be a non-void family of  curves in $Z$, we define
its \emph{$G_k$-combinatorial $p$-modulus} by
$$\Mod _p (\mathcal{F}, G_k) = \inf _{\rho} M_p(\rho),$$
where the infimum is over all \emph{$\mathcal{F}$-admissible functions} 
\emph{i.e.} functions  $\rho : G_k ^0 \to \hr_+ $ which satisfy 
$ L_{\rho} (\gamma) \geq 1$ for every $\gamma \in \mathcal{F}$.
If $\mathcal{F} = \emptyset$ we set $\Mod _p (\mathcal{F}, G_k) = 0$.
Observe that admissible functions with minimal
$p$-mass are smaller than or equal to $1$. 

\medskip

We denote  by $\mathcal F (A, B)$ the family of curves joining 
two subsets $A$ and
$B$ of $Z$ and by $\Mod _p (A, B, G_k)$ its  $G_k$-combinatorial
$p$-modulus. 
The following properties are routine to verify.

\begin{proposition}\label{3.1}
\mbox{}
\begin{enumerate}
\item
If $\mathcal F_1 \subset \mathcal F_2$ then 
$\Mod _p (\mathcal{F}_1, G_k) \leq \Mod _p (\mathcal{F}_2, G_k)$.
\item 
Let $\mathcal F_1 ,..., \mathcal F_n$ be curves families then
$$ \Mod _p (\cup \mathcal{F}_j, G_k) \leq \sum \Mod _p (\mathcal{F}_j, G_k).$$
\item
Let $\mathcal F_1, \mathcal F_2$ be two curve families. 
Suppose  that each curve in $\mathcal F_1$ admits
a subcurve in $\mathcal F_2$, then 
$\Mod _p (\mathcal{F}_1, G_k) \leq \Mod _p (\mathcal{F}_2, G_k)$.
\end{enumerate}
\end{proposition}

Recall that a metric space $Z$ is called a \emph{doubling metric space} if 
there is a constant $n\in\mathbb{N}$ such that every ball $B$ can be covered by
at most $n$ balls of radius $\frac{r(B)}{2}$. 
For a doubling metric space the combinatorial modulus does not 
depend on the choice
of the graph approximation up to a multiplicative constant. 
More precisely we have :

\begin{proposition}\label{3.2}
Assume that $Z$ is a doubling metric space.
Then for every $\kappa,\kappa' \geq 1$ and every $p \geq 1$
there exists a constant 
$D \geq 1$ such that 
for any $k \in \hn$ and for any graphs $G_k, G'_k$ which
are respectively $\kappa$ and $\kappa'$-approximations of $Z$ on scale $k$,
one has 
$$D^{-1 }\Mod _p (\cdot , G_k) \leq \Mod _p (\cdot , G'_k) \leq 
D \Mod _p (\cdot , G_k) .$$
\end{proposition}

\begin{proof}
The doubling property allows one to bound
the maximal number $N$ of pieces of $G_k ^0$ which overlap a given 
piece of $G'_{k}{}^0$, in terms of $n,\kappa,\kappa'$.

Let $\mathcal F$ be a family of curves in $Z$. For every  
$\mathcal F$-admissible 
function $\rho : G_k ^0 \to \hr_+$ define $\rho ' 
: G'_k {}^0 \to \hr_+$ by
$$\rho '(v') = \max \{\rho(v)~ ; ~ v \in G_k ^0,~  
v \cap v' \neq \emptyset \}.$$
 For $v' \in G'_{k}{}^0$ one has
$$N \rho '(v') \ge \sum _{v \cap v' \neq \emptyset} \rho (v).$$
Hence  $N \rho'$ is a $\mathcal F$-admissible 
function for the graph $G'_k$. 
In addition, letting $N'$ be the maximal number of pieces of $G'_k {}^0$ 
which overlap a given 
piece of $G_{k} ^0$, we have
$$ M_p (\rho ') \leq  \sum _{v' \in G'_k {}^0} \sum _{v \cap v' \neq
  \emptyset} \rho (v) ^p \le N'  M_p (\rho),$$
so we get
$$\Mod _p (\mathcal F , G'_k) \leq N^p N' \Mod _p (\mathcal F , G_k).  $$
\end{proof}

The following lemma is sometimes useful to understand the asymptotic behaviour 
of a minimal admissible function when $k$ tends to $+\infty$.

\begin{lemma}\label{min}
Let $\mathcal F$ be a curve family in $Z$, let $G_k$ be a 
$\kappa$-approximation of $Z$ on scale $k$, and let
$\rho : G_k ^0 \to \hr_+$ be a $\mathcal F$-admissible 
function with minimal $p$-mass. 
For $v \in  G_k ^0$ define  
$$\mathcal {F} _v = \{ \gamma \in \mathcal F~;~ \gamma \cap v \neq 
\emptyset \}.$$
Then one has $\rho (v) \le 
\Mod _p (\mathcal {F} _v  , G_k)^{1/p}$.
\end{lemma}

\begin{proof}
Let $\rho _v : G_k ^0 \to \hr_+$ be a minimal $\mathcal {F} _v$-admissible 
function,  
and let $ \tilde{\rho}$
be the function on $G_k ^0$ defined by 
$ \tilde{\rho} (w) = \max \{\rho (w), \rho _v (w) \}$
for $w \neq v$, and by $ \tilde{\rho} (v) = \rho _v (v)$.
Clearly $\tilde{\rho}$ is a $\mathcal F$-admissible function,
thus:
$$
\Mod _p (\mathcal F , G_k) 
    \le  M_p (\tilde{\rho}) 
    \le \sum _{w \neq v} \rho (w) ^p + \sum _w \rho _v (w) ^p ,$$
which implies that   
$$   \Mod _p (\mathcal F , G_k)  \le 
\Mod _p (\mathcal F , G_k) - \rho (v) ^p + 
\Mod _p (\mathcal {F} _v  , G_k).$$ 
The statement follows.
\end{proof}

\subsection{The Combinatorial Loewner Property (CLP)} \label{comb}

In this subsection we define a combinatorial analog of the
Loewner property introduced by J. Heinonen
and P. Koskela in \cite{HeK}. This notion appears in \cite{K} Section 7.
We will show that the CLP has a number of features in common
with the Loewner property.  
Examples of spaces satisfying
the combinatorial Loewner property will be given in Sections
\ref{examples} and \ref{menger}. 

\medskip

We assume that $Z$ is a compact arcwise connected doubling metric
space, in particular Proposition \ref{3.2} holds. Let 
$\{ G_k \} _{k \in \hn}$ be a $\kappa$-approximation of $Z$.
Recall that the \emph{relative distance} between two disjoint 
non-degenerate continua $A, B \subset Z$ is 
$$\Delta (A, B) = \frac {\dist (A, B) } {\min \{\diam A, \diam B \}}.$$

\begin{definition}
\label{def-clp}
Suppose  $p >1$.   Then $Z$ satisfies the \emph{Combinatorial
$p$-Loewner Property} if there exist 
two positive increasing functions  $\phi , \psi $ on $(0, +\infty)$ 
with $\lim _{t \to 0} \psi (t) = 0$, such that for 
all disjoint non-degenerate continua $A, B \subset Z$
and for all $k$ with $2^{-k} \le \min \{\diam A, \diam B \} $ 
one has :
$$ \phi (\Delta (A, B) ^{-1}) \le \Mod _p (A, B, G_k)
\le \psi (\Delta (A, B) ^{-1}) .$$
We say that $Z$ satisfies the \emph{Combinatorial
Loewner Property} if it satisfies the Combinatorial
$p$-Loewner Property for some $p>1$.
\end{definition}

Note that thanks to Proposition \ref{3.2}, the definition 
does not depend on the choice of the $\kappa$-approximation.

\medskip

The following properties have been established by J. Heinonen and P. Koskela for Loewner
spaces (see \cite{HeK} Theorem 3.13 and Remark 3.19). Their proof generalizes verbatim
to the spaces which satisfy the combinatorial Loewner property.

\begin{proposition}\label{linear}
Assume that $Z$ satisfies the CLP, then :

\begin{enumerate}
\item
It is \emph{linearly connected}, in other words there exists a 
constant $C \ge 1$
such that any two points $z_1 , z_2 \in Z$ 
can be joined by a path of diameter less or equal to
$C d(z_1 , z_2)$,

\item 
It has \emph{no local cut point}, in other words no connected open subset is disconnected by removing
a point.
\end{enumerate}
\end{proposition}

We also have:
\begin{theorem} \label{theoloewner}

\mbox{}
\begin{enumerate}
\item 
If $Z$ is a compact Ahlfors $p$-regular, $p$-Loewner metric space,
then $Z$ satisfies the combinatorial $p$-Loewner property.
\item If $Z'$ is quasi-Moebius  homeomorphic to a compact space $Z$
satisfying the CLP, then $Z'$ also satisfies the CLP (with the same exponent).
\end{enumerate}
\end{theorem}

The proof of (1) involves transferring admissible functions
on the metric measure space to admissible functions on an associated
discrete approximation, and vice-versa.  The arguments are straightforward
imitations of those appearing in \cite{HeK,BK1,H}, so we omit them.

The proof of (2) is similar in spirit, except that the admissible functions
are transferred between two discrete approximations.  
It involves some of the techniques 
that will be used
frequently in the sequel.

\begin{proof}[Proof of (2)]
We start by some general observations. 
Let $f : Z \to Z'$ be a quasi-Moebius  homeomorphism. 
Since $Z$ and $Z'$
are bounded, $f$ is a quasi-symmetric homeomorphism (see \cite{Vai}).
The doubling property and linear connectedness are 
preserved by quasi-symmetric homeomorphisms. Since 
$Z$ is a doubling linearly connected space (see Prop. \ref{linear}), $Z'$ is also. 
Given a $\kappa '$-approximation $G_\ell '$ of $Z'$, the preimages
$f ^{-1} (v')$ of the pieces $v' \in G'_\ell {}^0$ form a covering
of $Z$ that we denote by $\mathcal U _\ell$.  It enjoys the following properties :
there exist constants $\lambda \ge 1$ and $N \in \hn$, depending 
only on $\kappa '$, on the doubling constants of $Z, Z'$,  and on 
the quasi-symmetric parameters of $f$, such that 
\begin{itemize}
\item [(i)] For every $u \in \mathcal U _\ell$, there is a ball $B_u \subset Z$ with
$\frac{1}{\lambda} B_u \subset u \subset B_u$.
\item [(ii)] For every $z \in Z$, the number of balls $B_u$ containing $z$
is bounded by $N$.
\item [(iii)] If $2B_u $ and $2B_v$ intersect,
their radii satisfy $r(B_u) \le \lambda r(B_v)$.
\end{itemize}
For disjoint continua $E, F \subset Z$ one defines in an obvious way their 
\emph{$p$-modulus
relative to the covering $\mathcal U _\ell$}, denoted by $\Mod _p (E,
F,\mathcal U _\ell)$, 
so that  $\Mod _p (E, F,\mathcal U _\ell) = \Mod _p (f(E), f(F), G_\ell ' )$.

Assuming that $Z$ satisfies the combinatorial $p$-Loewner property, we will
compare $\Mod _p (E, F,\mathcal U _\ell)$ and $\Mod _p (E, F, G_k)$ for $k
\gg \ell$.
The statement (2) will follow since, for a quasi-Moebius homeomorphism $f$,
the relative distances $\Delta (E, F)$ and $\Delta (f(E), f(F))$
are quantitatively related (see \cite{BK1} Lemma 3.2). We begin by the :

\smallskip 

\emph{Left hand side CLP inequalities. }   We are looking for 
a positive increasing function $\phi '$ on $(0, +\infty)$, such that
for all disjoint non-degenerate continua $E,F \subset Z$, 
and all $\ell$
with $2 ^{- \ell} \le \min \{\diam f(E),  \diam f(F) \}$ : 
\begin{equation} \label{lhs}
\phi ' (\Delta (f(E), f(F)) ^{-1}) \le \Mod _p (f(E), f(F), G'
_\ell). 
\end{equation}
We may
assume that $2^{-\ell}$ is small in comparison with 
$\dist (f(E), f(F))$. Indeed, in the contrary, there is a family of curves joining $f(E)$ and
$f(F)$, which lies in the union of a controlled amount of pieces of $G_\ell '$
(recall that $Z'$ is a doubling linearly connected space).
Its $G_\ell '$-modulus is bounded from below in terms of the number of these
pieces. Thus, via the correspondence induced by $f$, we may
assume in addition that  
\begin{itemize}
\item [(iv)] None of the $2B_u$ ($u \in \mathcal U _\ell$) intersects both continua
  $E, F$.
\end{itemize}
 We will show the existence of a constant $C \ge 0$, such that for
every  $E,F, \ell $ as
above and for all $k \gg \ell$, one has
\begin{equation} \label{lhs2}
\Mod _p (E, F, G_k) \le C \Mod _p (E, F,\mathcal U _\ell).
\end{equation}
As explained in the first part of the proof, inequality (\ref{lhs}) will
follow. 
Let $\varphi 
: \mathcal U _\ell \to \hr _+ $ be a 
$\mathcal F (E,F)$-admissible function. For $k \gg \ell$
we wish to define a $\mathcal F (E,F)$-admissible function
of $G_k ^0$ whose $p$-mass is controlled by above by $M _p (\varphi) $.
To this aim, for every $u \in \mathcal U_\ell$, consider a minimal 
$\mathcal F (\overline{B_u}, Z \setminus \frac{3}{2}B_u)$-admissible
function $\rho _u : G_k ^0 \to \hr _+ $, and let 
$\rho _k : G_k ^0 \to \hr _+$ be defined by :
$$ \forall v \in G_k ^0, ~~ \rho _k (v) = \sum _{u \in \mathcal U} \varphi (u)
\rho _u (v) .$$
For every $\gamma \in \mathcal F (E,F)$, one has with (iv) and the definition 
of $\rho _u$ 
$$1 \le \sum _{u \cap \gamma \neq \emptyset} \varphi (u)
\le \sum _{u \cap \gamma \neq \emptyset} \varphi (u)
\Big( \sum _{v \cap \gamma \cap 2B_u \neq \emptyset} \rho _u (v) \Big) .$$
Using property (ii) we obtain that $N \rho_k$ is $\mathcal F (E,F)$-admissible.

To estimate $p$-masses, observe that $\rho _u$ being minimal, it is
supported on the set of  $v \in G_k ^0$ such that $v \subset 2B_u$.
Moreover its $p$-mass is smaller than $\psi (1/2)$ since
$Z$ satisfies the CLP. In combination with (ii), these properties show
that the $p$-mass of $\rho _k$ is 
$$\sum _{v \in G_k ^0} \Big( \sum _{u \in \mathcal U, v \subset 2B_u} \varphi (u)
\rho _u (v) \Big)^p \le N ^{p-1}\sum _{v \in G_k ^0} \sum _{u \in \mathcal U} 
\varphi (u) ^p \rho _u (v) ^p ,$$
which is less than  
$N ^{p-1} \psi (1/2) M_p (\varphi)$. 
 Therefore inequality \ref{lhs2} holds with $C = N^{2p-1}\psi
(1/2)$.
It remains to establish the :

\medskip

\emph{Right hand side CLP inequalities. }  We are looking for
a positive increasing function $\psi '$ on $(0, +\infty)$,
with  $\lim _{t \to 0} \psi' (t) = 0$, such
that for all disjoint non-degenerate continua $E,F \subset Z$ and all $\ell$
with $2 ^{- \ell} \le \min \{\diam f(E),  \diam f(F) \}$, one has
\begin{equation*} 
\Mod _p (f(E), f(F), G'
_\ell) \le \psi ' (\Delta (f(E), f(F)) ^{-1}). 
\end{equation*}
As explained in the first part of the proof, it is enough to find 
a constant $D \ge 0$, such that for all $E, F, \ell$ as above, and for all $k \gg
\ell$ :
\begin{equation} \label{rhs2}
\Mod _p (E, F,\mathcal U _\ell) \le D \Mod _p (E, F, G_k).
\end{equation}  
It requires the following general observation, that will also be used
frequently in the sequel :

\begin{lemma} \label{dual} let $\mathcal {F}$ be a family of curves
in a general compact metric space $Z$ and let $M >0$. 
Then $\Mod _p (\mathcal{F}, G_k) \geq M$ if and only if for every
function $\rho : G_k ^0 \to \hr _+$ there exists a curve $\gamma \in \mathcal{F} $
with
$$ L_{\rho} (\gamma) \le \Big ( \frac{M_p (\rho)}{M} \Big ) ^{1/p}.$$
\end{lemma}

\begin{proof} Just notice that the modulus can be written as
$$\Mod _p (\mathcal{F}, G_k) = \inf _{\rho} \frac{M_p (\rho)}{L_{\rho}
  (\mathcal F ) ^p}
~~\textrm{with}~~L_{\rho} (\mathcal F )
= \inf _{\gamma \in \mathcal F } L_{\rho} (\gamma),$$
where $\rho$ is any positive function on $G_k ^0$. 
\end{proof}

Let $k \gg \ell$ and let $\rho _k : G_k ^0 \to \hr_+ $ be any function. According to
Lemma \ref{dual}, to establish inequality (\ref{rhs2}), 
it is enough to construct a curve $\gamma \in \mathcal F (E,F)$ 
whose $\rho_k$-length is controlled by above by 
$$\Big( \frac{M_p (\rho _k)}{\Mod _p (E, F,\mathcal U_\ell)} \Big) ^{1/p}.$$ 
Recall that $2^{-\ell} \le \min \{\diam f(E), \diam f(F)
\}$. Increasing $\lambda$ if necessary, it yields :
\begin{itemize}  
\item [(v)] Whenever $B_u$ intersects $E$, the diameter of 
$2B_u \cap E$ is larger than $\frac {1}{\lambda} r(B_u)$ (and the same holds for
$F$ too).
\end{itemize}
Let $\Lambda \ge 1$ be a constant (that will be specified later on), and let 
$\varphi : \mathcal U_\ell \to \hr _+$ be the function :
$$ \forall u \in \mathcal U_\ell, ~~ \varphi (u) = \Big( \sum _{v \cap 2 \Lambda B_u \neq
  \emptyset} \rho _k (v) ^p \Big) ^{1/p} .$$
The obvious generalization of property (ii) implies that 
\begin{equation}
\label{comparison} 
M_p (\rho _k) \asymp M_p (\varphi) \,.
\end{equation} 
Lemma \ref{dual} shows that there is a curve $\delta \in \mathcal F (E,F)$
such that
\begin{equation}
\label{length}
L_{\varphi} (\delta) \le \Big( \frac{M_p (\varphi)}{\Mod _p (E, F,\mathcal
  U_\ell)} \Big) ^{1/p} \,.
\end{equation}
Let $u_i \in \mathcal U_\ell$ so that $\delta$ enters successively $u_1, ..., u_n$,
and set $B_i := B_{u_i}$ for simplicity. 
%Trimming and reindexing the sequence 
%$\{u_i \}_{i=1} ^n $ if necessary, we may assume that
%\begin{itemize}  
%\item [(vi)] $B_1 \cap E \neq \emptyset $ and $B_i \cap E = \emptyset $ for
%  $i>1$.
%\end{itemize}
We will now use the following lemma whose proof is similar to the one 
of Lemma 3.17 in \cite{HeK}. 

\begin{lemma}\label{clpball} Let $Z$ be a compact metric space satisfying the CLP. Then for every $\alpha \in
(0,1)$ there exist constants
$\Lambda  \ge 1$ and $m >0$, such that for every ball $B \subset Z$ and every disjoint
continua $E_1, E_2 \subset B$ with $\diam E_i \ge \alpha r(B)$, the
$G_k$-combinatorial $p$-modulus of the family
$$\{\eta \in \mathcal F (E_1, E_2)~;~ \eta \subset \Lambda B \}$$
is greater than $m$, for every $k$ with $2^{-k} \le \min \{\diam E_1, \diam E_2 \}$.
\end{lemma}
 
This lemma in combination with Lemma \ref{dual} and properties (iii),
(v), allows one to construct by induction on $s \in \{1, ..., n-1\}$ a curve 
$\gamma _s \subset \cup _{i=1} ^s 2 \Lambda B_i$, joining $E$ to $B_{s+1}$, whose
$\rho_k$-length is bounded linearly by above by $\sum _{i=1} ^s \varphi
(u_i)$. Indeed this follows from letting $\alpha =
\frac {1}{2 \lambda}$, $\Lambda = \Lambda (\alpha)$  and $B = 2B_s$
in the statement of the above lemma.  
A step futher gives a curve $\gamma \in \mathcal F (E,F)$ whose
$\rho_k$-length is bounded linearly by above by 
$L_{\varphi} (\delta) $. Thanks to the estimates (\ref{comparison}) and
(\ref{length}), 
the curve $\gamma$ enjoys the
expected properties. The theorem follows. 
\end{proof}

\subsection{Combinatorial ball Loewner condition}

Our next result is a combinatorial version of
Proposition 3.1 in \cite{BK3}. It asserts  that a space that satisfies a combinatorial
Loewner type condition for pairs of balls satisfies the combinatorial Loewner condition for
all pairs of continua. It is a main tool to exhibit examples of spaces satisfying the
combinatorial Loewner property.

\begin{proposition} \label{ball} Let $p \ge 1$. Assume that for every
$A >0$ there exist constants $m = m(A) >0$
and $L = L(A) >0$ such that if $r > 0$ and $B_1, B_2 \subset Z$
are $r$-balls with $\dist(B_1, B_2) \le A r$, then
for every $k \ge 0$ with $2^{-k} \le r$ 
the $G_k$-combinatorial $p$-modulus of the family
$$ \{\gamma \in \mathcal{F} (B_1, B_2)~;~ \diam (\gamma) \le Lr \} $$
is greater than $m$.
Then there exists a positive increasing function $\phi$ of $(0, +\infty) $
such that for 
every disjoint non-degenerate continua $E_1, E_2 \subset Z$
and for every $k$ with $2^{-k} \le \min \{\diam E_1, \diam E_2 \} $, one has :
$$ \phi (\Delta (E_1,  E_2) ^{-1}) \le \Mod _p (E_1, E_2, G_k).$$
\end{proposition}

Its proof is a rather straighforward discretization of the proof 
of Proposition
 3.1 in \cite{BK3}. For the sake of completeness and because similar ideas will be used in
Section \ref{ex}, we will give the details of the proof. Since it is the most
technical part of paper, readers may skip it at the first reading. 
The proof requires 
the following lemma which is the analog of Lemma 3.7 in \cite{BK3}.

\begin{lemma} \label{path} Let $Z$ as in Proposition \ref{ball} and 
suppose 
 $0< \lambda <1/8$. 
There exist constants $\Lambda = \Lambda (\lambda)$
 and $C = C(\lambda)$ with the following property. Let $\rho : G_k ^0 \to \hr _+$ be any function, 
let $B = B(z,r)$ be a ball of radius $0 < r < \diam Z$, 
and let $F_1, F_2 \subset Z$ be two continua with $F_i \cap \frac{1}{4} B \neq \emptyset$ and $F_i
\setminus B \neq \emptyset$ for $i= 1,2$. Then for $2^{-k+2} \le \lambda r$ 
there exist disjoint balls $B_i$, $i= 1,2$ 
and a path $\sigma \subset Z$ such that :
\begin{itemize}
\item[(i)] $B_1$ and $B_2$ are disjoint, they are centered on $F_i$ and of radius $\lambda r$, 
\item [(ii)] $B_i \subset \frac{7}{8} B$ and 
$$ \sum _{v \subset B_i} \rho (v) ^p  \le  8 \lambda \sum _{v \subset B} \rho (v) ^p ,$$
\item [(iii)] the path $\sigma$ joins $\frac{1}{4} B_1$ to $\frac{1}{4} B_2$, it is contained
in $\Lambda B$ and it has $\rho$-length at most 
$$ C \cdot \Big( \sum _{v \subset \Lambda B } \rho (v) ^p \Big ) ^{1/p} .$$
\end{itemize}
\end{lemma}

\begin{proof}
 We can find a subcontinuum $E_1 \subset F_1$ which is contained in $ \overline{B}(z, \frac{3r}{8}) \setminus
B(z, \frac{r}{4})$
and which joins $Z \setminus B(z, \frac{3r}{8})$ to $\overline{B}(z, \frac{r}{4})$. Similarly we can find a subcontinuum
$E_2 \subset F_2$ which is contained in $ \overline{B} (z, \frac{3r}{4})\setminus B(z, \frac{5r}{8})$
and which joins $Z \setminus B(z, \frac{3r}{4})$ to $\overline{B}(z, \frac{5r}{8})$. 
Then we have $\diam E_i \ge r/8$ for $i=1,2$ and $\dist (E_1, E_2) \ge r/4$.

\smallskip

Since $\diam E_i \ge r/8$,
for every $\lambda$ with $0< \lambda <1/8$ there exist at least
$\frac{1}{8 \lambda}$ pairwise disjoint balls centered on $E_i$ and of radius $\lambda r$.
So at least one of them - called $B_i$ - satisfies item (ii) of the statement. The condition (i)
is clearly satisfied by the pair of balls $B_1, B_2$.

\smallskip

Item (iii) follows from the hypothesis on the modulus of curves joining 
the pair of balls $\frac{1}{4} B_1, 
\frac{1}{4} B_2$ and from Lemma \ref{dual}. Note that the radius of $\frac{1}{4}B_i$ is equal to 
$\frac{\lambda r}{4}$ therefore our hypotheses requires that $2 ^{-k} \le \frac{\lambda r}{4}$ 
\emph{i.e.} $2^{-k+2} \le \lambda r$.
\end{proof}

\begin{proof} [Proof of Proposition \ref{ball}.]
Let $\lambda \in \hr $ subject to the conditions $0< \lambda <1/8$ and $2 \cdot (8 \lambda) ^{1/p} <1$.

Suppose $E_1, E_2 \subset Z$ are disjoint non-degenerate continua and let $\rho : G_k ^0 \to \hr _+$.
According to Lemma \ref{dual} we are looking
for a curve $\gamma$ joining $E_1$ to $E_2$ whose $\rho$-length is at most 
$$\Big( \frac{M_p (\rho)}{M} \Big) ^{1/p},$$
where $M>0$ depends only on the relative distance between $E_1$ and $E_2$.

\smallskip

Pick $p_i \in E_i$ such that $d(p_1, p_2) = \dist (E_1, E_2)$. Set
$$ r_0 := \frac{1}{2} \min \{ d(p_1, p_2), \diam E_1, \diam E_2 \} >0.$$
Let $B_i = B(p_i, r_0)$ for $i=1,2$. Then $B_1 \cap B_2 = \emptyset$ and $E_i \setminus B_i \neq
\emptyset$ for $i=1,2$. In addition :
$$\dist (\frac{1}{4}B_1, \frac{1}{4}B_2) \le \Big ( \frac{4d(p_1, p_2)}{r_0 } \Big ) \frac{r_0}{4} \le 
t \frac{r_0}{4} , $$
where $t := 8 \max \{ 1, \Delta ( E_1, E_2) \}$.
By our hypotheses one can find a path
$\sigma$ joining $\frac{1}{4}B_1$ to $\frac{1}{4}B_2$ and whose $\rho$-length is at most 
$( \frac{M_p (\rho)}{m} ) ^{1/p}$
where $m = m(t)$ is the constant appearing in the statement of Proposition \ref{ball}.

\smallskip

By using Lemma \ref{path}
inductively,  we will construct for successive integers $n$ a family
of balls $\mathcal{B} _n$, a collection of continua $\Omega _n$, and a collection of paths $\Sigma _n$, such that :
\begin{itemize}
\item [(1)] $\mathcal{B} _0 = \{B_1, B_2\}$, $\Omega _0 = \{E_1, E_2 \}$ and  $\Sigma _0 = \{ \sigma \}$ 
are defined previously .
\item [(2)] For $n \ge 1$ the family $\mathcal{B} _n$ consists of $2^{n+1}$ 
disjoint balls of radius $\lambda ^n r_0$.
Each ball $B \in \mathcal{B} _n$ is centered on a continuum 
$\omega \in \Omega _n$ with $\omega \setminus B \neq \emptyset$.
The collection $\Sigma _n$ consists of $2^{n}$ paths, for each element $\sigma \in \Sigma _n$ 
there are exactly 
two elements  $B_1, B_2 \in \mathcal{B} _n$ such that $\sigma$ joins $\frac{1}{4}B_1$ to $\frac{1}{4}B_2$.
\end{itemize}
The induction proceeds as follows.
According to item (2) for every ball $B \in \mathcal{B} _n$
there exist $\omega \in \Omega _n$ and $\sigma \in \Sigma _n$ 
such that the pair $\{F_1, F_2 \} := \{\omega,  \sigma \}$ satisfies 
$F_i \cap \frac{1}{4} B \neq \emptyset$ and $F_i
\setminus B \neq \emptyset$ for $i= 1,2$. Thus applying Lemma \ref{path} we get two disjoint balls
$B_1, B_2$ 
and a path $\sigma$ joining $\frac{1}{4}B_1$ to $\frac{1}{4}B_2$. 
The balls $B_i$ are of radius $\lambda ^{n+1} r_0$, they are centered respectively on $\omega $ and $\sigma$.
The definitions of $\mathcal{B} _{n+1}$ and $\Sigma _{n+1}$ are therefore clear.
We define $\Omega _{n+1}$ to be $ \Omega _n \cup \Sigma _n$.

\smallskip

With Lemma \ref{path} and by construction the following additional properties are satisfied :
\begin{itemize}
\item [(3)] For every ball $B \in \mathcal{B} _n$ there exists a ball $B' \in \mathcal{B} _{n-1}$
such that $B \subset \frac{7}{8} B'$ and 
$$\sum _{v \subset B} \rho (v) ^p  \le  8 \lambda \sum _{v \subset B'} \rho (v) ^p .$$
\item [(4)] For every path $\sigma \in \Sigma _n$ 
there exists a ball $B' \in \mathcal{B} _{n-1}$
such that $\sigma$ lies in $\Lambda B'$ and has $\rho$-length at most 
$$ C \cdot \Big( \sum _{v \subset \Lambda B' } \rho (v) ^p \Big ) ^{1/p} ,$$
where $\Lambda = \Lambda (\lambda)$ and $C = C(\lambda)$ .
\item [(5)] At each stage $n$ one can index the elements of $ \bigcup _{\ell \le n} \Sigma _{\ell} $
by $\sigma _1, ... , \sigma _m $ and the elements of $\mathcal B _n$ by $B_1, ... , B_{m+1}$,
in order that $B_1$ meets $E_1$ and $\sigma _1$, $B_i$ meets $\sigma _{i-1}$
and $\sigma _i$ for $2 \le i \le m$  , and $B_{m+1}$ meets $\sigma _m$ and $E_2$.
\end{itemize}
We iterate this procedure as long as $n \in \hn$ satisfies $2 ^{-k +2} \le \lambda ^n r_0$. 
Let $N$ be the largest integer satisfying this condition. 

\smallskip

It remains to connect the paths $\sigma _i$ described in item (5) to obtain a curve joining $E_1$
to $E_2$. For this purpose observe that 
the hypotheses of Proposition \ref{ball} assert that for every pair $z_1, z_2$ of points in $Z$
there exists a path joining $B (z_1, \frac{1}{4}d(z_1, z_2))$ to $B (z_2, \frac{1}{4}d(z_1, z_2))$
whose diameter is comparable to $d(z_1, z_2)$. This property implies that every pair of points 
$z_1, z_2$ are connected by a path whose diameter is comparable to $d(z_1, z_2)$, (see 
Lemma 3.4
in \cite{BK3} for more details).
Therefore increasing $\Lambda$ if necessary and using item (5) with $n=N$ we exhibit a family 
$\Theta$ consisting of $2^{N+1}$ curves such that :
\begin{itemize}
\item [(6)] For every $\theta \in \Theta$ the union of the subsets $v \in G_k ^0$ with $v \cap
\theta \neq \emptyset$ is contained in a ball of the form $\Lambda B$ with $B \in \mathcal {B} _N $ .
\item[(7)] The following subset contains a curve $\gamma$ joining $E_1$ to $E_2$ :
$$ \bigcup _{\theta \in \Theta} \theta \cup \bigcup _{n \le N} \bigcup _{\sigma \in \Sigma _{n} } \sigma.$$
\end{itemize}
Observe that for $B \in \mathcal {B} _N$ the number of $v \in G_k ^0$ with $v \subset \Lambda B$ is bounded
in terms of $\Lambda, \kappa$ and the doubling constant of $Z$. Hence increasing $C$ if necessary
we obtain with property (6) that 
\begin{itemize}
\item [(8)] Every $\theta \in \Theta $ has $\rho$-length at most
$$ C \cdot \Big( \sum _{v \subset \Lambda B } \rho (v) ^p \Big ) ^{1/p} ,$$
where $B \in \mathcal {B} _N$ is the ball attached to $\theta$ in item (6).
\end{itemize}

\smallskip

We now compute the $\rho$-length of the curve $\gamma$ defined in (7). At first with properties (2) and 
(3) one obtains that for $0 \le s \le n \le N $ and for every ball $B \in \mathcal{B}_n $ there exists
a ball $B' \in \mathcal{B}_{n-s}  $ such that $(\frac{1}{8 \lambda})^s B \subset B'$. Let
$s = s(\lambda)$ be the smallest integer such that $(\frac{1}{8 \lambda})^s \ge \Lambda$. 

With properties (4) and (3)
we get that for $s < n \le N$ the $\rho$-length of every curve $\sigma \in \Sigma _n$ is less
than 
$$C \cdot \Big ( (8 \lambda) ^{n-s-1} M_p (\rho) \Big ) ^{1/p} .$$
Similarly with properties (8) and (3) we get that the $\rho$-length of every curve $\theta \in \Theta$
is less than
$$ C \cdot \Big ( (8 \lambda) ^{N-s} M_p (\rho) \Big ) ^{1/p} .$$
For $0 < n \le s$ each $\sigma \in \Sigma _n$ has $\rho$-length at most $ C (M_p (\rho)) ^{1/p}$.
Recall that $\sigma \in \Sigma _0$ has $\rho$-length at most $( \frac{M_p (\rho)}{m} ) ^{1/p}$.
Finally : 
$$ L_{\rho} (\gamma) \le M_p (\rho) ^{1/p} \Big ( \frac{1}{m ^{1/p}} + C \sum _{n = 1} ^{s} 2^n 
+ C \sum _{n = s+1} ^{N+1} 2^n (8 \lambda) ^{(n-s-1)/p} \Big ) .$$
Therefore letting $D := C \cdot (8\lambda)^{(-s-1)/p}$ and $a := 2 \cdot (8 \lambda) ^{1/p}$ we get
$$ L_{\rho} (\gamma) \le M_p (\rho) ^{1/p} \Big ( \frac{1}{m ^{1/p}} + D \sum _{n = 1} ^{N+1} a^n \Big),$$
and $$M := \Big ( \frac{1}{m ^{1/p}} +  D \sum _{n = 1} ^{N+1} a^n \Big ) ^{-p} $$
satisfies the desired properties since $a<1$ by assumption.
\end{proof}

\subsection*{Remarks}  1) We could have chosen to define the CLP 
using the more
general notion of $\kappa$-approximations from \cite{BK1}, so as to make the
quasi-Moebius invariance automatic from the definition.  However, this
would simply make it harder to verify in examples, forcing one
to prove the equivalence of the two definitions anyway.

\medskip

 2) We emphase that, unlike the CLP, the classical
Loewner property is not, strictly speaking, invariant 
by quasi-Moebius homeomorphisms.
Indeed assume that $Z$ and $Z'$ are quasi-Moebius equivalent metric spaces,
and that $Z$ is a Loewner space. Then by a theorem of J. Tyson \cite{Ty},
$Z'$ is a Loewner space if and only if the Hausdorff dimensions of $Z$ and
$Z'$ are equal.

For that reason it is in general a difficult problem to decide
whether a given metric space is quasi-Moebius equivalent to a 
Loewner space. The analogous problem for the CLP is in principle 
easier.

\section{Self-similarity}\label{hyp}

This section derives from self-similarity several general principles that will be useful
in the sequel.  In particular self-similarity permits the reduction of the right hand
side of the CLP inequality to a simpler condition (Proposition \ref{rd}). It also
implies a certain asymptotic behaviour for the
combinatorial modulus (Proposition \ref{submult} and Corollary
\ref{Q_M}). 

\subsection{Approximately self-similar spaces}

The following definition appears in \cite{K} Section 3. 

\begin{definition} \label {homothety} A compact metric space $(Z, d)$ is called \emph{approximately
self-similar} if there is a constant $L_0 \ge 1$ such that if $B(z,r) \subset Z$ is a ball 
of radius $0 < r \le \diam (Z)$, then there is an open subset $U \subset Z$ which is $L_0$-bi-Lipschitz
homeomorphic to the rescaled ball $(B(z,r), \frac{1}{r} d)$.
\end{definition}

Observe that approximately
self-similar metric spaces are doubling metric spaces. 
Examples include some classical fractal spaces like the square Sierpinski
carpet and the cubical  Menger sponge (their definitions
are recalled in Section \ref{menger}).
A further source of examples comes from expanding Thurston maps,
\cite{bonkmeyer}.  It follows readily from 
\cite[Theorem 1.2]{bonkmeyer} that with respect to any visual metric as in \cite{bonkmeyer},
the $2$-sphere is approximately self-similar.

\medskip 

 We now present examples arising from group actions on boundaries (Proposition
\ref{propself}). For this purpose, we first 
recall some standard notions and results 
(see \emph{e.g.} \cite{GO, BH, KapB} 
for more details). Let  $X$ be a
geodesic proper metric space. Assume that $X$ is \emph{hyperbolic} (in the sense of Gromov),
\emph{i.e.} there exists a
constant $\delta _X$, called the \emph{triangle fineness constant of $X$},
such that for every geodesic triangle 
$[x, y] \cup [y,z] \cup [z,x] \subset
X$  and for every $p \in [x,y]$, one has 
\begin{equation} \label{triangle}
\dist (p , [y,z] \cup [z,x]) \le \delta _X.  
\end{equation}
Denote by $\partial X$ the boundary at infinity of $X$. It carries a \emph{visual metric},
\emph{i.e.} a metric $\delta$ for which there are constants $a >1$, $C \ge 1$ such that
for every $z, z' \in \partial X$, one has
\begin{equation} \label{visualbis}
C ^{-1} a ^{- \ell} \le \delta (z,z') \le C a ^{- \ell},
\end{equation}
where $\ell$ denotes the distance from $x_0$ (an origin in $X$) to 
a geodesic $(z, z') \subset X$.
Suppose in addition that $\Gamma$ is a (Gromov) hyperbolic group 
that acts on $X$ by isometries, properly discontinuously and
cocompactly. 
Then $\Gamma$ acts canonically on $(\partial X, \delta)$ by uniformly quasi-Moebius 
homeomorphisms. 
Moreover, assuming that 
$\partial \Gamma$ and $\partial X$ are both equipped with visual metrics, 
the orbit map $g \in \Gamma \mapsto gx_0 \in X$
extends to a canonical quasi-Moebius homeomorphism $\partial \Gamma
\to \partial X$. We set :

\begin{definition} \label{self} Let $\Gamma$ be a hyperbolic group.
A metric $d$ on $\partial \Gamma$ is called 
a \emph{self-similar metric}, if there exists a hyperbolic geodesic metric space
$X$, on which $\Gamma$ acts by isometries, properly discontinuously and cocompactly, 
such that $d$ is the preimage of a visual metric $\delta$ on $\partial X$ by 
the canonical homeomorphism
$\partial \Gamma \to \partial X$.
\end{definition}

\begin{proposition} \label{propself} The space $\partial \Gamma$ equipped with a self-similar 
metric $d$ is approximately
self-similar, the partial bi-Lipschitz maps being restrictions 
of group elements. 
Moreover $\Gamma$ acts on $(\partial \Gamma, d)$ by 
(non-uniformly) bi-Lipschitz homeomorphisms, and $(\partial \Gamma, d)$
is linearly connected as soon as it is connected.
\end{proposition}

\begin{proof}
 Let $X$ and $\delta$ be as in Definition \ref{self}, let $x_0 \in X$ be an origin
and let $\delta _X$ be the triangle fineness constant of $X$. Let $a$ and $C$ be
the constants associated to $\delta$ and $x_0$ by the relations (\ref{visualbis}). 
We denote the
distance in $X$ between $x$ and $y$ by $\vert x - y \vert$. We are looking
for a constant $L_0$ such that for every ball $B(z,r) \subset (\partial X,
\delta)$
with $0<r<\diam \partial X$, there is $g \in \Gamma$ whose restriction to
$B(z,r)$ is a $L_0$-bi-Lipschitz homeomorphism from the rescaled ball
$(B(z,r), \frac{1}{r} \delta)$ onto an open subset of $\partial X$.

\smallskip

Let $B(z,r)$ as above and let $y \in [x_0, z) \subset X$ be such that
$$\vert x_0 - y
\vert = - \log _a 2rC - \delta _X.$$ 
At first we claim that for every $z_1, z_2  
\in B(z,r)$ and every $p \in (z_1, z_2)$, one has 
$$\dist (y, [x_0, p]) \le 3 \delta _X .$$
To see it, notice that $\delta (z_1, z_2) \le 2r$. Hence from inequalities
(\ref{visualbis}) one gets : $\dist (x_0, (z_1,z_2)) \ge - \log _a 2rC$.
Let $y_1 \in [x_0, z_1)$ be such that
$$\vert x_0 - y_1 \vert = - \log _a 2rC - \delta
_X.$$
One has : $\dist (y_1, (z_1,z_2)) \ge d(x_0, (z_1,z_2)) - \vert x_0 - y_1 \vert
\ge \delta _X .$ For every $p \in (z_1,z_2)$ consider the geodesic triangle 
$[x_0, z_1) \cup (z_1, p] \cup [p, x_0]$. With the inegality (\ref{triangle})
we get : $\dist (y_1, [x_0, p]) \le \delta _X$. 
For $z_2 = z$ and $p = z_2$ we obtain in particular that $\vert y - y_1 \vert 
\le 2 \delta _X$. The claim follows.

\smallskip

Let $D$ be the diameter of the quotient space $X / \Gamma$, and let $g \in
\Gamma$ be such that $\vert g^{-1} x_0 - p \vert \le D$. For any $z_1, z_2
\in B(z,r)$, we deduce from the claim that 
$$
\arrowvert \dist (g^{-1} x_0 , (z_1, z_2))-\dist (x_0, (z_1,z_2))-\log _a 2rC
 ~\arrowvert  \le D + 6\delta _X .$$
One has :
$\dist (g^{-1} x_0 , (z_1, z_2)) = \dist (x_0 , (g z_1, g z_2))$
since $\Gamma$ acts on $X$ by isometries.
With inegalities (\ref{visualbis}) it leads to 
$$L_0 ^{-1} \frac{\delta (z_1,z_2)}{r} \le \delta (g z_1, g z_2)
\le L_0  \frac{\delta (z_1,z_2)}{r},$$
where $L_0 = \frac{C^3}{2}a^{D+6\delta _X}$. The first part of the statement
follows. 

\smallskip

The fact that $\Gamma$ acts on $(\partial X, \delta)$ by (non-uniformly)
bi-Lipschitz homeomorphisms can be proved by similar arguments. The
connectedness statement is a result from \cite{BK2}.
\end{proof}

\subsection{Modulus on approximately self-similar spaces }

For the rest of the section $Z$ denotes an arcwise connected,
approximately self-similar metric space.  Let $\{ G_k \} _{k \in \hn}$ 
be a $\kappa$-approximation of $Z$. We fix a positive constant $d_0$ that is small compared to 
the diameter of $Z$  and to the constant $L_0$ of the Definition \ref{homothety}. 
Denote by $\mathcal F _0$ the family of curves $\gamma \subset Z$ with $\diam
(\gamma) \ge d_0$.
 A first result concerns the right hand side CLP inequality  :

\begin{proposition}\label{rd} Let $p >1$ and suppose 
there exists a constant $C \geq 0 $
such that for every $k \in \hn$ one has $ \Mod _p (\mathcal{F} _0 , G_k) 
\le C$. 
Then there exists a positive increasing function $\psi$ of
$(0, +\infty) $ with $\lim _{t \to 0} \psi (t) = 0$, such that for every disjoint
non-degenerate continua $A, B \subset Z$ and every integer $k$ 
satisfying  $2 ^{-k} \le \min \{\diam (A), \diam (B)\}$, one has 
$$ \Mod _p (A, B, G_k) \le \psi (\Delta (A, B) ^{-1}).$$

\end{proposition}

\begin{proof}
Define for $t > 0$ : 
$$\psi (t) = \sup \Mod _p (A, B, G_k) ,$$
where the supremum is over all disjoint continua $A, B$
with $\Delta (A, B) \ge 1/t$  and over all integers $k$
such that $2 ^{-k} \le \min \{\diam (A), \diam (B)\} $.
>From the monotonicity of the modulus (Proposition \ref{3.1}.1) and from 
our hypotheses,
the function $t \mapsto \psi(t)$ is non-decreasing with 
non-negative real 
values. It is enough to prove that $\psi(t)$ tends to $0$ when $t$
tends to $0$.

Let $A$ and $B$ be disjoint non-degenerate continua, 
assume that $d := \diam (A)$ is smaller than $\diam (B)$ and let $n$
be the largest integer with $2^n \cdot d \le \dist (A, B)$.
Pick $z_0 \in A$, we get :
$$ A \subset B(z_0 , d) ~\textrm{and}~B \subset Z \setminus  
B(z_0 , 2^n \cdot d).$$
For $i \in \{1, ... , n-1 \}$ define $B_i = B(z_0, 2^i \cdot d)$.
There exists a constant $C_1$ depending only on $\kappa$ 
and on the geometry of $Z$
such that for every $k \in \hn$ with $2 ^{-k} \le d$ and every 
$i \in \{1, ... , n-1 \}$
one has
$$ \Mod _p (B_i, Z \setminus B_{i + 1} , G_k) \le C_1.$$
Indeed applying the self-similarity property (Definition \ref{homothety}) to $B_{i+1}$ we may inflate 
$\mathcal F (B_i, Z \setminus B_{i+ 1} )$ to a family of curves
of essentially unit diameters.
Hence the above inequalities follow from the monotonicity of the modulus, 
from our hypotheses,
and from Proposition \ref{3.2}.

\smallskip

Choose for every $i \in \{1, ... , n-1 \}$ a minimal
$\mathcal F (B_i, Z \setminus B_{i+ 1} )$-admissible function
$\rho _i : G_k ^0 \to \hr _+ $, and let
$\rho = \frac{1}{n-1} \sum _{i =1} ^{n-1} \rho _i$.
This is a $\mathcal F (A, B)$-admissible function since every
curve joining $A$ and $B$ joins $B_i$ and $Z \setminus B_{i + 1}$ too.
For $2^{-k} \le d$ the minimality of the $\rho _i$'s shows that 
their supports are essentially
disjoint, thus
$$M_p (\rho) \le 
\frac{C_2}{(n-1)^{p} } \sum _{i = 1} ^{n-1} M_p (\rho _i)
\le \frac{C_1 C_2}{(n-1)^{p-1} },$$
where $C_2$ depends only on $\kappa$.
Since $p >1$ and $\Delta (A, B) \le 2^{n +1}$ we get that for $t$ small enough 
$\psi (t) \le C_3 (\log 1/t)^{-1}$ where $C_3$ depends only on $C_1, C_2, p$.
\end{proof}

The following uniformization criterion is a consequence of 
an uniformization theorem established in \cite{BK1} (see also \cite{CFP} 
for related results).

\begin{corollary} \label{corunif}
Suppose in addition that $Z$ is 
homeomorphic to the $2$-sphere. 
Assume that there exists a constant
$C \ge 1$ such that for every $k \in \hn $ one has 
$$\Mod _2 (\mathcal F _0 , G_k ) \le C \,.$$
Then $Z$ is quasi-Moebius homeomorphic to the Euclidean $2$-sphere.
\end{corollary}

\begin{proof} The statement is a consequence of the previous
proposition in combination with \cite{BK1} Th. 10.4.
This theorem supposes that $Z$ is doubling and linearly
locally contractible. 

The doubling property follows  
from the fact that $Z$ is approximately self-similar.

Recall that a metric space $Z$ is \emph{linearly locally contractible}
if there exists a constant $\lambda \ge 1$, such that every ball $B \subset Z$
with $0< r(B) < \frac{\diam Z}{\lambda}$, is contractible in $\lambda B$.
Approximately
self-similar manifolds enjoy this property.
Indeed, if $Z$ is a compact manifold, then there is a $R_0 >0$ and a positive function
$\Phi$ of $(0, R_0)$  with $\lim _{t \to 0} \Phi (t) = 0$, such that every ball
$B \subset Z$ of radius $r$ with $0< r < R_0 $, is contractible in $\Phi (r) B$.
Let $\lambda \ge 1$ and let $B(z,r) \subset Z$ be a ball with $\lambda r <
\diam Z$.
Applying self-similarity to the ball $B(z, \lambda r)$, one obtains a open
subset $U \subset Z$ and a bilipschitz map $f : B(z, \lambda r) \to U$
such that for every pair $x,y \in B(z, \lambda r)$ :
$$ L_0 ^{-1} \frac{d(x,y)}{\lambda r}  \le d(f(x), f(y)) \le 
 L_0 \frac{d(x,y)}{\lambda r} .$$
Thus if $\lambda$ satisfies $\frac{L_0}{\lambda} \le L_0 ^{-1} $,
one has
$$f(B(z,r)) \subset B(f(z),\frac{L_0}{\lambda}) \subset  B(f(z), L_0 ^{-1})
\subset f(B(z, \lambda r)).$$ 
Choosing $\lambda$ subject to the conditions 
$\Phi (\frac{L_0}{\lambda}) < L_0 ^{-1}$ and 
$\frac{L_0}{\lambda} \le L_0 ^{-1} $, we get that $f(B(z,r))$ is
contractible in  $f(B(z, \lambda r))$. Therefore so is $B(z,r)$ in $B(z,
\lambda r)$.

Assumptions 
of Theorem 10.4 of \cite{BK1} also
require the graphs $G_k$ to be (essentially) homeomorphic
to $1$-skeletons of  triangulations of $S^2$. Corollary 6.8 of
\cite{BK1} ensures existence of such graphs. 

We remark finally that when $Z$ is the boundary of a hyperbolic group, the corollary
can also be deduced from \cite{CFP} Th. 1.5 and 8.2.
\end{proof} 

We now study the behaviour of $\Mod _p (\mathcal{F} _0 , G_k) $ 
when $k$ tends to $+ \infty$, depending on $p \ge 1$. For this purpose we establish a 
submultiplicative inequality :

\begin{proposition}\label{submult}
Let $p \ge 1$ and write $M_k := \Mod _p (\mathcal{F} _0 , G_k) $ for simplicity. 
There exists a constant $C>0$ such that for every pair of integers
$k, \ell$ one has : $M_{k+\ell} \le C \cdot M_k \cdot M_{\ell}$.
In addition, when $p$ belongs to a compact subset of $[1, +\infty)$,
the constant $C$ may be choosen independent of $p$.
\end{proposition}

\begin{proof} Let $\rho _k : G_k ^0 \to \hr _+ $ be a minimal
$\mathcal F _0$-admissible function on scale $k$.
By definition
each $v \in G_k ^0$ is 
roughly a ball of $Z$ of radius
$2^{-k}$. For every $v \in G_k ^0$, let $B_v$ be a ball containing
$v$ and whose radius is approximately $2^{-k}$.

\smallskip

 We inflate every ball $2B_v$ 
to essentially unit diameter as in Definition \ref{homothety}.
Let $g_v$ be the corresponding partial bi-Lipschitz map. Define 
$G_{k+\ell} \cap 2B_v$ to be the incidence graph of the covering
of $2B_v$ by the   subsets $w \in G^0 _{k+\ell}$ with $w \cap 2B_v \neq
\emptyset$. Using the map $g_v$, one can consider the graph
$G_{k+\ell} \cap 2B_v$ as a $(2\kappa L_0)$-approximation
of $g_v (2B_v)$ on scale $\ell$. We pull back to $G_{k+\ell} \cap 2B_v$ 
a normalized minimal $\mathcal F _0$-admissible function
on scale $\ell$, in order to get for every $v \in G_k ^0$
a function $\rho _v : (G_{k+\ell} \cap 2B_v)^0 \to \hr _+ $
with the following properties 
\begin{itemize}
\item[(i)] its $p$-mass is bounded by above by $D \cdot M_{\ell}
\cdot \rho _k (v) ^p$,
where $D$ is independent of $k$ and $v$,
\item[(ii)] every curve $\gamma \subset 2B_v$ whose diameter is larger than
$r(B_v)$ picks up a $(\rho_v)$-length larger than
$\rho _k (v) $, where $r(B_v)$ is the radius of $B_v$.
\end{itemize}
Define a function $\rho_{k+\ell}$ on $G_{k+\ell}^0 $ by
$$ \rho_{k+\ell} (w) = \max \rho _v (w), $$
where the maximun is over all $v \in G_{k} ^0 $ with
$w \cap 2B_v \neq \emptyset$.
Its $p$-mass is linearly bounded by above by
$$\sum _{v \in G_k ^0} M_p (\rho _v) ,$$
which in turn is linearly bounded by $M_p (\rho_k) \cdot M_{\ell} 
= M_k \cdot M_{\ell} $ (see item (i)).

\smallskip

It remains to prove that $\rho_{k+\ell}$ is a $\mathcal{F} _0$-admissible
function -- up to a multiplicative constant independent of the scale.
For $\gamma \in \mathcal{F} _0$ we have
$$ 1 \le \sum _{v \cap \gamma \neq \emptyset} \rho _k (v) 
\le \sum _{v \cap \gamma \neq \emptyset} 
\sum _{\gamma \cap 2B_v \cap w \neq \emptyset} \rho _v (w). $$
Indeed the last inequality follows from item (ii) since the relation
$v \cap \gamma \neq \emptyset$ implies that there exists a subcurve
of $\gamma$ of diameter greater than $r(B_v)$ contained in $2B_v$.

Thus the $\rho_{k+\ell}$-length of $\gamma$ is larger than $1/N$
where $N$ is maximal number of elements $v \in G_k ^0$ such that
$2B_v$ intersects a given piece $w \in G_{k + \ell} ^0$ .
Therefore $N \cdot \rho_{k+\ell}$ is $\mathcal{F} _0$-admissible.
\end{proof}

For every $k \in \hn$, observe that  
$p \mapsto \Mod _p (\mathcal{F} _0 , G_k) $ is a non-increasing
continuous function on $[1, +\infty)$ (monotonicity comes 
from the fact that minimal admissible
functions are smaller than or equal to $1$).  
We define a critical exponent associated to the 
curve family $\mathcal{F} _0$ by
\begin{equation}\label{defQ_M} 
Q_M := \inf \{p \in [1, +\infty) ~;
\lim _{k \to +\infty} \Mod _p (\mathcal{F} _0 , G_k) = 0 \}.
\end{equation}
With Propositions \ref{rd}, \ref{submult} and Lemma \ref{min} one gets

\begin{corollary} \label{Q_M}
\mbox{}
\begin{itemize}
\item [(1)] For $p>Q_M$ one has 
$\lim _{k \to +\infty} \Mod _p (\mathcal{F} _0 , G_k) = 0 $,
\item [(2)] For $1 \le p \le Q_M$ the sequence  $\{ \Mod _p (\mathcal{F} _0 , G_k) \}
_{k \ge 0} $ admits a positive lower bound,
\item [(3)] If in addition $Z$ is linearly connected, then for 
$1 \le p < Q_M$ the sequence  $\{ \Mod _p (\mathcal{F} _0 , G_k) \} _{k \ge 0} $ is unbounded.
\end{itemize}
In particular when $Z$ satisfies the combinatorial $p$-Loewner 
property one has $p = Q_M$.
\end{corollary}

\begin{proof} 
Part (1) is a consequense of the definition. Part (2) comes from
Proposition \ref{submult} and from the fact that if a positive sequence
$\{M_k \} _{k \in \hn} $ satisfies $M_{k+\ell} \le C \cdot M_k \cdot
M_{\ell}$, then one has :
$$\lim _{k \to +\infty} M_k = 0 ~~\Longleftrightarrow~~ \exists \ell \in \hn ~~\textrm{with}~~
M_{\ell} < C^{-1} .$$ 

\smallskip

To establish
part (3) suppose by contradiction that the sequence is bounded
for some $Q$ with $1<Q<Q_M$. Then the conclusion of 
Proposition \ref{rd} holds for the exponent $Q$. 
We will prove that for $p > Q$
the sequence $\{ \Mod _p (\mathcal{F} _0 , G_k) \}
_{k \ge 0} $ tends to $0$, contradicting
the definition of $Q_M$.

\smallskip

 Consider a finite covering of $Z$ by balls of radius
$\frac{d_0}{4}$. Every $\gamma \in \mathcal F _0$ satisfies 
$\diam (\gamma) \ge d_0$, thus meets at least two
disjoint balls of the covering. Therefore,  
by Proposition
\ref{3.1}.2, it is enough to show that for every disjoint
non-degenerate continua $A$ and $B$ of $Z$ the sequence $\{ \Mod _p (A, B, G_k) \}
_{k \ge 0} $ tends to $0$.
To do so, we will establish that the $\mathcal {F} (A, B)$-admissible functions 
$\rho _k : G_k ^0 \to \hr _+$ with minimal $Q$-mass satisfy
$$\lim _{k \to +\infty} \Vert \rho _k \Vert _{\infty} =0,
~~\textrm{where}~~ \Vert \rho _k \Vert _{\infty}
:= \sup _{v \in G_k ^0}  \rho _k (v).$$
For each $v \in G_k ^0$ pick a continuum $E_v$
containing $v$
 and whose diameter $ d_v $ is comparable to  
$\diam (v) \asymp 2^{-k}$.
The existence of $E_v$ follows from the assumption that $Z$ is linear connected.
Any curve in $\mathcal F (A,B)$ passing
through $v$ possesses a subcurve in each family $ \mathcal F (A, E_v)$
and $ \mathcal F (B, E_v)$. Hence with Lemma \ref{min} and 
Proposition \ref{3.1} we get 
$$ \rho _k (v) \le \min \{\Mod _Q (A, E_v, G_k), \Mod _Q (B, E_v, G_k) \}
^{1/Q}.$$
For $2^{-k}$ small enough compared with $\dist (A, B)$ one has
$$\max \{\Delta (A, E_v), \Delta (B, E_v) \} \ge \dist (A, B) /2 d_v .$$
Therefore Proposition \ref{rd} applied with exponent $Q$ shows
that $$ \rho _k (v) \le \psi (2 d_v / \dist (A, B))^{1/Q} ,$$
and so $\Vert \rho _k \Vert _{\infty} $ 
tends to $0$ when $k$ tends to $+\infty$.

\smallskip

Finally the CLP assertion follows from Proposition \ref{linear}.1.
\end{proof} 

\subsection*{Remarks }
 
1) In \cite{H2} the authors have proved that if $Z$
satisfies the combinatorial $Q$-Loewner property then for every disjoint
non-degenerate continua $A, B \subset Z$ one has for $1 \le p < Q$ :
$$\lim _{k \to +\infty} \Mod _p (A, B, G_k) = + \infty . $$

2) In \cite{KeK} S. Keith and the second (named) author
establish that $Q_M$ is equal to the \emph{Ahlfors regular conformal dimension} of $Z$
\emph{i.e.} the infimum of Hausdorff dimensions of Ahlfors regular
metric spaces quasi-Moebius homeomorphic to $Z$.  
 A proof of this equality is also part of the 
forthcoming Phd Thesis \cite{Car}. 
As a consequence $Q_M$ is a quasi-Moebius
invariant of $Z$.

\medskip
3) M. Barlow and R. Bass have established submultiplicative and supermultiplicative
inequalities for the combinatorial $2$-modulus on some self-similar space like the 
square Sierpinski carpet.
Their method relies on the analysis of random walks on graph approximations
(see \cite{Ba} Lemma 3.2).

\medskip
4)
One may formulate a variant of Proposition \ref{submult} for doubling spaces
in general, by modifying the definition of  $M_k$.   For 
constants $C_1, C_2$, and each $k$,  let $M_k$ be the supremum,
as $x$ ranges over $Z$, of the $p$-modulus of the family 
of curves of diameter at least $C_12^{-k}$ lying in the ball
$B(x,C_22^{-k})$.  Then for suitably chosen $C_1$, $C_2$, one obtains a
submultiplicative inequality $M_{k+\ell}\leq C\cdot M_k \cdot M_\ell$.

\section{The Combinatorial Loewner Property for the standard 
Sierpinski carpet and Menger curve} \label{menger}

The \emph{square Sierpinski carpet} is the continuum 
$\hs\subset [0,1]^2\subset \hr^2$ 
constructed as follows : start with the unit square
in the plane, subdivide it into nine equal subsquares, remove the middle 
open square,
and then repeat this procedure inductively on the remaining squares.
The \emph{cubical Menger sponge} is
the continuum  
$$
\hm=\bigcap_{i=1}^3\;\pi_i^{-1}(\hs)\,, 
$$
where $\pi_i:\hr^3\ra \hr^2$ is the map which forgets the $i^{th}$
coordinate.
Both $\hs$ and $\hm$ are endowed with the induced Euclidean metric.
This section establishes : 

\begin{theorem} \label{hs}
$\hs$ and $\hm$ satisfy the Combinatorial
Loewner Property.
\end{theorem}

We notice that similar reasoning applies
to other families of self-similar examples, 
see Subsection \ref{subsec-otherexamples}.

Many ideas appearing in this paper have their origin in the analysis
of the combinatorial modulus on $\hs$ and $\hm$. Futhermore, 
most of the ideas present in the proof of Theorem \ref{hs}
will be reused to study the boundaries of 
Coxeter groups in the next five Sections.
We think that this exposition may be useful to
the reader to understand the Coxeter group case. 

We will treat  $\hs$ and $\hm$ separately. Due to planarity, the Sierpinski
carpet case is  simpler. One can view the Menger sponge case as 
an intermediate stage
between the Sierpinski carpet and  Coxeter group boundaries.

\subsection{Proof of Theorem \ref{hs} for the Sierpinski carpet}
The proof is divided in several steps.

\subsection*{Curves and crossing}
We start by analysing the dynamics of the curves in $\hs$.
Consider
the isometric action $D_4\acts [0,1]^2$ of the dihedral group $D_4$,
and let $\De\subset [0,1]^2$ be the triangular fundamental domain $\{(x,y)\in[0,1]^2
~;~ y \leq x\leq \frac12\}$.   We identify $\De$ with the orbit space $[0,1]^2/D_4$, and hence
we have the orbit map $\pi:[0,1]^2 \ra \De$.

\smallskip

Pick $k \in \hn$ and consider the tiling $\T_k$ of $\hs$ by homothetic copies of itself 
of side length $3^{-k}$, so that each $T \in \T_k$ is contained in a square of the form 
$[\frac{m-1}{3^k},\frac{m}{3^k}]\times [\frac{n-1}{3^k},\frac{n}{3^k}]$ with $1\leq m,n\leq 3^k$.
For every $T\in\T_k$,  let $\phi_T:T \ra \hs$ be the homeomorphism
which is a composition of a translation and scaling.  
We have a (continuous) folding map $\pi_k: \hs \ra \De$ whose restriction to each tile
$T\in\T_k$ is the composition
$$T\stackrel{\phi_T}{\ra} \hs \subset [0,1]^2 \stackrel{\pi}{\ra}\De.$$
%We also have a collection of similarities $\De \ra [0,1]^2$, 
%which are compositions
%$\De\stackrel{g}{\ra} [0,1]^2 \stackrel{\phi_T^{-1}}{\ra}T$, for $g\in D_4$, $T\in\T_k$.
\begin{lemma}
\label{lem-squarecoxeter}
Suppose $k\in \hn$ and  $\gamma \subset \hs $ is a curve.  
Then one of the following holds:
\begin{itemize}
\item [(1)] The curve $\pi_k (\gamma)$ intersects all three (closed) sides of $\De$,  and for every 
(parametrized) curve
$\eta \subset \hs$ the subset
$\pi _k ^{-1} (\pi_k(\gamma)) \subset \hs$ 
contains a curve which approximates $\eta$ to within error
$2 \cdot 3^{-k}$ with respect to the
$C^0$ distance. 
\item [(2)] The curve $\gamma$ is disjoint from $\pi_k^{-1}(\si)$ for one of the (closed) sides $\si$ of $\De$, and 
$\diam (\gamma) \le 2 \cdot 3^{-k}$. 
\end{itemize}
\end{lemma}

\begin{proof}  
Suppose the path $\pi_k(\gamma)$ intersects all three sides of $\De$, and consider
elements $g_1,g_2\in D_4$ such  $g_1\De$ and $g_2\De$ meet along an edge $\si$.  Then 
the paths $g_1\circ\pi_k (\gamma)$ and $g_2\circ \pi_k (\gamma)$ both touch
$\si$, and differ by reflection across $\si$. 
Therefore $g_1\circ\pi _k (\gamma)\cup g_2\circ\pi _k(\gamma)$ is path connected.
It follows that  $\delta := \cup_{h\in D_4}\; h \circ \pi _k (\gamma)$ is
path connected, and that it touches every sides of $[0,1]^2$. 

Let $\eta \subset \hs$ be a curve, and let $T_1 ,...,T_n$ be a sequence of tiles
in $\T_k$ such that every two consecutive tiles share a side,
and such that $\eta$ enters successively the $T_i$'s. Then $\cup _{i =1} ^n \phi _{T_i} ^{-1} (\delta)$
is a subset of $\pi _k ^{-1} (\pi_k(\gamma))$.
Moreover for  every $i = 1, ..., n-1$ the subsets  
$\phi _{T_i} ^{-1} (\delta)$ and $\phi _{T_{i+1}} ^{-1} (\delta)$ both meet the side $T_i \cap T_{i+1}$
and differ 
by reflection along $T_i \cap T_{i+1}$.  Hence $\phi _{T_i} ^{-1} (\delta) \cup \phi _{T_{i+1}} ^{-1} (\delta)$
is path connected. Therefore $\cup _{i =1} ^n \phi _{T_i} ^{-1} (\delta)$ contains a curve 
which approximates $\eta$ to within $2 \cdot 3^{-k}$ with respect to the
$C^0$ distance.

To prove the second part of the alternative 
we just notice that an Euclidean reflection group generated
by a pair of reflections along sides of $\De$ fixes one of the vertices of $\De$. 
\end{proof}

\subsection*{Bounds for the combinatorial modulus}

We now analyse the combinatorial modulus using the previous lemma. 
Let $G_k$ be the
incidence graph of the tiling $\mathcal T _k $. The collection $\{ G_k \} _{k \in \hn} $ is a 
$\kappa$-approximation of $\hs$ (\footnote{More rigourously there is a subsequence of 
$\{ G_k \} _{k \in \hn} $ which is a $\kappa$-approximation of $\hs$. Indeed any
interval of the form $[\frac{1}{2^k}, \frac{3}{2^k}]$ contains a power of
$\frac{1}{3}$.}). 
As usual we identify every
vertex of $G_k$ with the corresponding subset in $\hs$.

\smallskip

Let $d_0$ be a fixed (small) positive constant and let
$\mathcal F _0$ be the family of curves $\gamma \subset \hs$ with $\diam (\gamma) \ge d_0$.
For $\epsilon > 0$ and for any (parametrized) curve $\eta$ in $\hs$, 
let $\mathcal U _{\epsilon} (\eta)$ be the set of curves whose $C^0$
distance to $\eta$ is smaller than $\epsilon$.

\begin{lemma} \label{hypsimpleCLP}
For every $p \ge 1$ and for
every $\epsilon >0$, there exists a constant $C = C(p, \epsilon)$
such that for every non constant curve $\eta \subset \hs$ and for every $k \in \hn$ :
$$\Mod _p (\mathcal {F} _0 , G_k) \le C \Mod _p (\mathcal {U} _{\epsilon}
(\eta), G_k) .$$
Futhermore, when $p$ belongs to a compact subset of $[1, +\infty)$,
the constant $C$ may be choosen independent of $p$.
\end{lemma}

\begin{proof}
Pick $\ell \in \hn$ large enough so that $\epsilon, d_0 \ge 2 \cdot 3^{-\ell}$.
Then, by the previous lemma, for every pair of curves $\eta \subset Z$ and 
$\gamma \in \mathcal F _0$ the subset
$\pi _\ell ^{-1} (\pi_\ell (\gamma)) \subset \hs$ contains a curve of 
$\mathcal {U} _{\epsilon} (\eta)$.

Given $T, T' \in \mathcal T _\ell$ and $h \in D_4$, we have an isometry
$g_{T,T'\!,h} : T \to T'$ which is the composition of a translation and
$\phi _T ^{-1} \circ h \circ \phi _T$. Observe that for every curve $\gamma \subset
\hs$ :
\begin{equation} \label{equahyp}
\pi _\ell ^{-1} (\pi_\ell (\gamma)) = \cup \{g_{T,T'\!,h} (\gamma \cap T) ~;~
T,T' \in \mathcal T _\ell ~, h \in D_4 \}.
\end{equation}
Let $k \ge \ell$. For every $\mathcal {U} _{\epsilon}(\eta) $-admissible function
$\rho : G^0 _k \to \hr_+$  define a function 
$\rho ': G_k ^0 \to \hr_+$ as follows. For $v \in G_k ^0$, let  $T \in \mathcal T _\ell$
be so that $v \subset T$, and set : 
$$\rho' (v) = \sum _{T',h} \rho (g_{T,T'\!,h} (v)),$$
where $T' \in \mathcal T_\ell $ and $h \in D_4$. 
We claim that $\rho '$ is 
$\mathcal F _0$-admissible. Indeed let $\gamma
\in \mathcal F _0$, and let $\theta \in  
\mathcal U _{\epsilon} (\eta)$ be so that $\theta \subset 
\pi _\ell ^{-1} (\pi_\ell (\gamma))$. One has with (\ref{equahyp}) :
$$L _{\rho '} (\gamma) = \sum _{T,T'\!,h} ~\sum _{\substack{v \subset T \\ v \cap \gamma \neq \emptyset}}
\rho (g_{T,T'\!,h} (v)) = L_\rho (\pi _\ell ^{-1} (\pi_\ell (\gamma)) \ge 
L _\rho (\theta) \ge 1 .$$
The claim follows. Moreover one has
$$M_p (\rho ') = \vert \mathcal T _\ell \vert \cdot \vert D_4 \vert \cdot M_p (\rho).$$
Thus 
$$\Mod _p (\mathcal F _0 , G _k)
\le \vert \mathcal T _\ell \vert \cdot \vert D_4 \vert \cdot \Mod _p (\mathcal U _{\epsilon} (\eta) , G_k).$$
Therefore the statement holds for $k \ge \ell$. Since $\ell$ depends only on $\epsilon, d_0$
the lemma follows. 
\end{proof}

\subsection*{Multiplicativity}

We now study the asymptotic behaviour of the modulus $\Mod _p (\mathcal F _0 , G _k)$
using the planarity of $\hs$. We have

\begin{lemma}\label{sup} Let $p \ge 1$ and write $M_k := \Mod _p (\mathcal F _0 , G _k)$ for simplicity.
There exists a constants $C\ge1$ such that for every pair of integers $k, \ell$ one has
$$C^{-1} \cdot M_k \cdot M_\ell \le M_{k+\ell} \le C \cdot M_k \cdot M_\ell .$$
\end{lemma}

\begin{proof}
The right hand side inequality follows from Proposition \ref{submult}. 

To prove the 
opposite inequality, let $\mathcal F_\hs$ be the family of of curves in $\hs$ joining
the sides $\{0\} \times [0,1]$ and $\{1\} \times [0,1]$. Define $m_k := 
\Mod _p (\mathcal F _\hs , G _k)$. More generally, given a tile $T \in \cup _{k\in \hn}
\mathcal T _k$, let $\mathcal F _T$ be the family of curves in $T$ which join
the left and right sides of $T$. For $k, \ell \in \hn$ and $T \in \mathcal T _k$ 
one has: $\Mod _p (\mathcal F _T , G _{k+\ell}) = m_\ell$. We claim that
$m_k \asymp M_k$. Indeed we have $\mathcal F _\hs \subset \mathcal F _0$, thus $m_k \le M_k$.
Pick $T := ([\frac{1}{3},\frac{2}{3}] \times [0,\frac{1}{3}]) 
\cap \hs \in \mathcal T_1$
and $\eta := [0,1] \times \{0\} \subset \hs$. Every curve in $\mathcal U_{\frac{1}{3}} (\eta)$
admits a subcurve in $\mathcal F _T$. Therefore, in combination with Lemma 
\ref{hypsimpleCLP}, we obtain :
$$m_k = \Mod _p (\mathcal F _T , G _{k+1}) \ge  \Mod _p (\mathcal U_{\frac{1}{3}} (\eta),
 G _{k+1}) \gtrsim M_{k+1} \asymp M_k.$$
The claim follows. 

\smallskip

We shall prove that $m_{k + \ell} \gtrsim m_k m_\ell$. 
To do so let $\rho_{k+\ell}:G_{k+l}^0\ra \R_+$ be a function. We wish to construct a curve $\gamma
\in \mathcal F _\hs$ with controlled $\rho_{k+\ell}$-length.  Define 
$\rho_{k+\ell}':G_{k+\ell}^0\ra \R_+$
as follows.  For each $v \in G_{k+\ell}^0$, let $T\in \mathcal T_k$ be so that $v \subset T$,
and 
%, and $v \let $G_{k+l}^{T,0}=\{W\in G_{k+1}^0\mid W\subset T\}$.
set $\rho_{k+\ell}' (v)$ to be the sum of the  
$\rho_{k+\ell}\circ \phi (v)$ where $\phi: T \ra T'$ ranges over
isometries of $T$ with tiles $T'$ intersecting $T$. 

Pick $T \in \mathcal T_k$.  By Lemma \ref{dual} there is a curve $\ga_T\in \F_T$ such that
\begin{equation}\label{equasuper}
L_{ \chi_T \rho_{k+\ell}'}(\ga_T)\leq\left(  \frac{M_p( \chi _T \rho_{k+\ell}')}{m_\ell}  
\right)^{\frac1p}\,,
\end{equation}
where $\chi _T: G_{k+l}^0\ra \R_+$ is the charateristic function of the set of $v \in G_{k+l}^0$
which are contained in $T$.   
Let $\ga_T'$ be the union of the images of $\ga_T$ under all isometries 
$T\ra T'$, where $T'\in \mathcal T _k$ intersects $T$.   Then $\ga_T'$
is path connected, and
\begin{equation}\label{equasuper2}
L_{ \chi _T \rho_{k+\ell}}(\ga_T') = L_{ \chi _T \rho_{k+\ell}'}(\ga_T).
\end{equation}
Note that if $T_1,T_2\in \mathcal T_k$ intersect, then $\ga_{T_1}'$ intersects $\ga_{T_2}'$,
since $\ga_{T_2}'$ contains a curve in $T_1$ joining the top and bottom, which must
intersect $\ga_{T_1}$ by planarity. 

Define $\rho_k:G_k^0\ra \R_+$ by letting $(\rho_k(T))^p$ be the 
$p$-mass of $\chi _T \rho_{k+\ell}'$.
By Lemma \ref{dual}, 
there is a vertex path $T_1,\ldots,T_n$ in $G_k$ such that $T_1$ intersects $\{0\}\times [0,1]$,
$T_n$ intersects $\{1\}\times [0,1]$, and whose $\rho_k$-length satisfies
\begin{equation}\label{equasuper3}
\sum_{i}\rho_k(T_i)\leq \left(\frac{M_p(\rho_k)}{m_k}\right)^\frac1p
=\left(\frac{M_p(\rho_{k+\ell}')}{m_k}\right)^\frac1p
\end{equation}
Then the union $\cup_i \ga_{T_i}'$ is a path connected subset of $\S$ joining the 
left and right sides. Therefore it contains a curve $\gamma \in \mathcal F _\hs$ whose
$\rho_{k+\ell}$-length is at most (thanks to (\ref{equasuper2}), (\ref{equasuper}) and
(\ref{equasuper3})) :
$$
\sum_i L_{\chi _{T_i} \rho_{k+\ell}}(\ga_{T_i}') =
\sum_i L_{\chi _{T_i} \rho_{k+\ell}'}(\ga_{T_i})
\leq 
\sum_i \left(  \frac{M_p(\chi _{T _i} \rho_{k+\ell}')}{m_\ell}  \right)^{\frac1p}
$$
$$
=\sum_i   \frac{\rho_k(T_i)}{m_\ell ^{1/p}} 
\leq \left( \frac{M_p(\rho_{k+\ell}')}{m_km_\ell} \right)^\frac1p 
\asymp \left( \frac{M_p(\rho_{k+\ell})}{m_km_\ell} \right)^\frac1p .
$$
With  Lemma \ref{dual}, inequality $m_{k + \ell} \gtrsim m_k m_\ell$ follows.
\end{proof}

\begin{proof}[Proof of Theorem \ref{hs} for the Sierpinski carpet, concluded]
Let $Q_M$ be the critical exponent defined in (\ref{defQ_M}). We
establish the CLP for $p=Q_M$.
Lemma \ref{sup} and the proof of Corollary \ref{Q_M}.2 imply that
$M_k \asymp 1$ for $p = Q_M$. Moreover one has $Q_M >1$. Indeed 
for every $N \in \hn$, there is a disjoint collection of paths $\ga_1,\ldots, \ga_N \in\F _0$.
It follows that with $p=1$ and $k$ large enough :  $\Mod _1(\F_0 ,G_k)\geq N$. Since 
$\Mod _p(\F _0,G_k)$ is bounded for $p=Q_M$, we obtain $Q_M >1$.
Therefore Proposition \ref{rd} shows that the right hand side CLP inequality
holds.  On the other hand, the lower bound $M_k \gtrsim 1$, Lemma \ref{hypsimpleCLP}
and self-similarity
imply that the assumptions of Proposition \ref{ball} are satisfied. Thus the left
hand side CLP inequality holds.
\end{proof}

\subsection{A criterion for a self-similar space to satisfy the 
CLP,  and the proof of Theorem \ref{hs} for the Menger sponge} 
Observe that $\hm$ satisfies the obvious analog of Lemma 
\ref{lem-squarecoxeter} with $D_4$ replaced by the symmetry 
group of the cube. Therefore Lemma \ref{hypsimpleCLP} is also valid 
for $\hm$. Moreover, as in the last paragraph, one has (for $p=1$) :
$$ \lim _{k \to +\infty} \Mod _1(\F_0 ,G_k) = +\infty.$$
The following general result implies that $\hm$ admits the CLP. 
Notation is the same as in Lemma \ref{hypsimpleCLP}.

\begin{proposition} \label{simpleCLP}
Let $Z$ be a linearly connected, approximately self-similar
metric space. For $p=1$ assume that $\{\Mod _p(\F_0 ,G_k)\}_{k \in \hn}$ 
is
unbounded. Suppose
that for every $p \ge 1$, for every non constant curve $\eta \subset Z$ and for
every $\epsilon >0$, there exists a constant $C = C(p, \eta, \epsilon)$
such that for every $k \in \hn$ :
$$\Mod _p (\mathcal {F} _0 , G_k) \le C \Mod _p (\mathcal {U} _{\epsilon}
(\eta), G_k) .$$
Suppose futhermore that when $p$ belongs to a compact subset of $[1, +\infty)$
the constant $C$ may be choosen independent of $p$.
Then $Z$ satisfies the CLP.
\end{proposition}

The proof uses some arguments that will be discussed in greater generality
in Section \ref{ex}. We refer to this section for more details. 

In the following two lemmata
the hypotheses of  Proposition \ref{simpleCLP} are assumed  to hold. The first lemma 
is a replacement of the supermultiplicativity that was established in Lemma
\ref{sup} using planarity.

\begin{lemma}
Write $M_k := \Mod _p (\mathcal F _0 , G_k)$ and $L_k := M_k ^{-1/p} $ for 
simplicity. 
There exist constants $C \ge 1$ and $b \in (0,1)$ such
that for every $k, \ell \in \hn$ one has :
$$ L_{k+ \ell} \le 
C \cdot L _k \cdot \sum _{n=0} ^\ell L_{\ell - n} ~b ^n .$$
Morever when $p$ belongs to a compact subset of $[1, +\infty)$ the constants $C$ and $b$
may be choosen independent of $p$.
\end{lemma}

\begin{proof}
The assumptions on $Z$ allow one to check that the proofs
of Lemmata  \ref{gen1}, \ref{gen2} and of Proposition \ref{super} apply
verbatim to $Z$, with $M_k$ instead of $m_k$.
Therefore the inequality holds.
\end{proof}

\begin{lemma} \label{behaviour}
For $p=Q_M$ the sequence $\{M_k \} _{k \in \hn} $ admits a
  positive lower bound, and it doesn't tend to $+ \infty$.
\end{lemma}

\begin{proof}
The first part of the statement follows from Corollary \ref{Q_M}.2.
The previous lemma shows that the set of  $p$ such
that  $\{M_k \} _{k \in \hn} $ tends to $+ \infty$ is an open subset
of $[1, + \infty)$ --- see the proof of  Corollary \ref{Q_m} for more details.
Therefore the second part of the statement holds.
\end{proof}

\begin{proof} [Proof of Proposition \ref{simpleCLP}]
It is enough to prove that for $p=Q_M$ one has $M_k \asymp 1$.
The $p$-combinatorial Loewner property will then 
come from Lemma \ref{gen1} and from  Propositions \ref{ball} and \ref{rd}. Note that 
Proposition \ref{rd} requires that $p >1$. At this point we use the hypothesis on
$\{\Mod _1(\F_0 ,G_k)\}_{k \in \hn}$ to ensure that $Q_M > 1$. 

Thanks to Lemma \ref{behaviour}  the sequence $\{M_k \} _{k \in \hn} $ admits a
positive lower bound. To obtain the upper bound, one argues exactly as in the second part of
the proof of Theorem \ref{CLP}, using the fact that $\{M_k \} _{k \in \hn} $
does not tend to $+ \infty$.
\end{proof}

\subsection{Other examples}
\label{subsec-otherexamples}

Similar reasoning establishes the CLP for the following examples:

1) \emph{Higher dimensional Menger spaces.}  For $n\geq 2$, 
the Sierpinski
carpet  construction may be generalized
by iterating the subdivision of the 
unit cube $[0,1]^n$ into $3^n$ subcubes of side length $\frac{1}{3}$,
and removing the central open subcube.  Denote by $\hs_n$ the resulting 
space. Then for $\ell \geq n$, we get
an analog of the Menger space by letting
$$
\hm_{\ell, n}=\cap\{\pi_I^{-1}(\hs_n)~;~ I\subset 
\{1,\ldots,\ell\}, |I|=n\}\,,
$$
where $\pi_I:[0,1]^{\ell} \ra [0,1]^n$ is the projection map which retains the
coordinates indexed by elements of $I$.  These constructions can be further
generalized by subdividing into $k^n$ subcubes instead of $3^n$, where
$k$ is an odd integer, or by removing a symmetric pattern of subcubes
at each stage, instead of just the central cube, etc.

\medskip

2) \emph{Higher dimensional snowspheres, cf. \cite{meyer}.}
Fix $n\geq 2$.  Let $Q_k$ denote the $n$-dimensional
polyhedron obtained from the
boundary of the cube $[0,3^{-k}]^{n+1}$ by removing the interior of 
one $n$-dimensional face.  
Construct a sequence $\{P_k\}_{k\in \mathbb{Z}_+}$, 
where $P_k$ is a
metric polyhedron consisting of Euclidean $n$-cubes of side length
$3^{-k}$, as follows.
Let $P_0$ be the boundary of the unit cube $[0,1]^n$, and inductively construct $P_k$ from $P_{k-1}$
by subdividing each $n$-cube face of $P_{k-1}$
into $3^n$ subcubes, removing the 
central open subcube, and gluing on a copy of $Q_k$ along the
boundary.   If we endow $P_k$ with the path metric, then the
sequence $\{P_k\}$ Gromov-Hausdorff converges to a self-similar
space $Z$.  As with Menger spaces, there are further generalizations
of these constructions.

\medskip

3) \emph{Pontryagin manifolds.}
We modify slightly the above construction to obtain some Pontryagin manifolds.
Fix $n \ge 2$, and let $T ^n$ be the standard $n$-torus obtained by
identification of the opposite faces of the unit cube $[0,1]^n$.
Let $Q_k$ denote the polyhedron obtained as follows. Tesselate $T^n$
by $3^n$ equal subcubes, remove the interior of one of them, and
normalize the metric so that the side length of every subcube
is $3^{-k}$. Define $P_0 = T^n$ and inductively construct $P_k$ from $P_{k-1}$
by subdividing each $n$-cube face of $P_{k-1}$
into $3^n$ subcubes, removing the 
central open subcube, and gluing on a copy of $Q_k$ along the
boundary. Then, as above, $\{P_k\}$ Gromov-Hausdorff converges to a self-similar
space $Z$.

\section{Hyperbolic Coxeter groups}\label{preliminary}

This section exhibits some specific dynamical properties of the action 
of a (Gromov) hyperbolic Coxeter group on its boundary.  Results of
this section will serve to study the combinatorial modulus in Section \ref{cox} 
and the $\ell _p$-equivalence relations in Section \ref{equivalence}.

\subsection{Definitions and first properties }

We start by recalling some standard definitions, 
see \cite{D} for more details.
Let $(\Gamma, S)$ be a Coxeter  group (see section \ref{statement}
for the definition). 
A \emph{special subgroup} of $\Gamma$
is a subgroup generated by a non-empty subset $I$ of $S$; we shall denote 
it by
$\Gamma _I$. A \emph{parabolic subgroup} of
$\Gamma$ is a conjugate of a special subgroup. 
Let $G$ be the Cayley graph of $(\Gamma, S)$.
We define $G^0$ and $G^1$ to be the set of vertices and of
(non-oriented) open edges respectively. Each edge carries a \emph{type} 
which is an element
of $S$. The distance between two vertices
$x, y$ of $G$ is denoted by $\vert x - y \vert$,  and the distance
from the identity
to $x$ by $\vert x \vert$. An edge with
endpoints $x,y$ is noted $(x,y)$.
For $s \in S$ the \emph{wall associated to $s$} is the subset 
$M_s \subset G^1$
of $s$-invariant 
(open)
edges. The graph $G \setminus M_s$
consists of two disjoint  convex closed subsets of $G$, denoted by 
$H_{-}(M_s)$ 
and $H_{+}(M_s)$,
and called the \emph{half-spaces bounded by $M_s$}.
They satisfy the relations
$$H_{-}(M_s)^0 = \{ g \in \Gamma~ ;~ \vert sg \vert = \vert g \vert +1 \}, \ \ 
  H_{+}(M_s)^0 = \{ g \in \Gamma~ ;~ \vert sg \vert = \vert g \vert -1 \},$$
and $s$ permutes $H_{-}(M_s)$ and $H_{+}(M_s)$. A \emph{wall} of $G$ is
a subset of the form $g(M_s)$, with $g \in \Gamma$ and $s \in S$. The 
involution $gsg^{-1}$ is called the \emph{reflection along the wall $g(M_s)$}.
The set of walls forms a partition of $G^1$. Each wall $M$ divides
$G$ in two disjoint  convex closed subsets, called the 
\emph{half-spaces} bounded by $M$ and denoted by $H_{-}(M)$ 
and $H_{+}(M)$, with the convention that $e \in H_{-}(M)$, where $e$
denotes the identity.
 
\medskip

Assume now that $\Gamma$ is a hyperbolic Coxeter group
(see \emph{e.g.} \cite{GO,BH, KapB}
for  hyperbolic groups and related topics).
We denote by $\partial \Gamma$ its boundary at infinity equipped with
a visual metric $d$ (see (\ref{visualbis})).
For a subset $E$ of 
$G$ we denote by $\partial E$ its limit set in $\partial \Gamma$.
A non-empty limit set of a parabolic subgroup of $\Gamma$
will be called a \emph{parabolic limit set}. If in addition it is
a topological circle we shall call it a \emph{circular limit set}.
 We notice that :

\begin{proposition} \label{fls}
For every $\epsilon > 0$ there is only a finite number of parabolic limit
sets of diameters greater than $\epsilon$.
\end{proposition}

\begin{proof}
Let $P = g \Gamma _I g^{-1}$ be a parabolic subgroup with non empty
limit set. One has
$\partial P = \partial (g \Gamma _I)$. Since $\Gamma$
is a Coxeter group, the inclusion map of the Cayley graph
of $(\Gamma _I, I)$ into $G$ is a totally geodesic isometric embedding.
By using the visual property 
(\ref{visualbis}) and the fact that $\Gamma$ acts by isometries on $G$,
we obtain that the diameter of $\partial (g \Gamma _I)$
is comparable to $a ^{-\dist (e , ~g \Gamma _I)}$. 
Only finitely many cosets $g \Gamma _I$  intersect
a given ball of $G$.  The statement follows. 
\end{proof}

The limit set of a wall $M$ is either of empty interior or equal
to $\partial \Gamma$. Indeed this property is well-known for limit
sets of subgroups, and the stabilizer of $M$ in $\Gamma$
acts cocompactly on $M$. Using the convexity 
of the half-spaces one easily sees that 
$$\partial  H_{-}(M) \cup \partial H_{+}(M) =  \partial \Gamma ,\ \
\partial  H_{-}(M) \cap \partial H_{+}(M) = \partial M.$$
In consequence, in $\partial \Gamma$,  the fixed point set of a reflection is
the limit set of its wall.

The following property asserts in particular that the limit sets of half-spaces 
form a basis of neighborhoods in $\partial \Gamma$.

\begin{proposition} \label{basis} There exists a constant $\lambda \ge 1$  
such that for every 
$z \in \partial \Gamma$
and every $0< r \le \diam (\partial \Gamma) $ 
there exists a half-space $H$ of $G$ with
$$ B(z, \lambda ^{-1} r) \subset \partial H  \subset  
B(z, \lambda r).$$
\end{proposition}

 The proof relies on the following lemma. In its statement
$\delta _G$ denotes the triangle fineness constant of $G$ (see (\ref{triangle})).

\begin{lemma} There exists a constant $L>0$ with the following property.
Let $x, x_1, y_1, y$ be arbitrary $4$ points of $G$ lying in this order on a geodesic line.
If $\vert x_1 - y_1 \vert \ge L$ then there exists a wall $M$ passing between $x_1$ and $y_1$ 
such that 
$$ \dist (x, M) \ge \vert x - x_1 \vert - 3 \delta _G ~~\textrm{and}~~ 
\dist (y, M) \ge \vert y - y_1 \vert - 3 \delta _G .$$
\end{lemma}

\begin{proof}
Observe that different edges of a geodesic segment give rise to different walls.
Thus there are as many walls passing between $x_1$ and $y_1$ as edges in $[ x_1, y_1 ]$. 
Let $(a,b)$ be such an edge and let $M$ be the corresponding wall.
Let $p \in G^0$ be a vertex  adjacent to $M$ realizing
$\dist (x, M)$. Observe that the geodesic segment $[a, p]$ 
lies in a $\delta _G$-neighborhood
of $M$. Indeed, $a$ and $p$ are adjacent to $M$ thus the Hausdorff distance between $[a,p]$ 
and its image by the reflection along $M$ 
is smaller than $2 \delta _G$, moreover $[a,p]$ and its image
are separated by $M$. 
Therefore we have 
\begin{equation} \label{basis1}
\dist (x , [a,p]) \ge \vert
x - p \vert - \delta _G.
\end{equation}
Assume that $\vert x - p \vert < \vert x - x_1 \vert  - 3 \delta _G$. 
Consider the triangle $[x,a] \cup [a,p] \cup [p,x]$ and the point $q \in 
[p,x]$ with $\vert q - p \vert = 2 \delta_G$. The inequality (\ref{basis1})
implies that the nearest point projection of $q$ on $[x,a] \cup [a,p]$ lies in $[x,a]$.
Thus one gets that 
$$
\dist ( p , [x, x_1]) \le 3 \delta _G .
$$
Pick $r \in [x, x_1]$
with $\vert r - p \vert \le 3 \delta _G$. Consider the triangle 
$[r,a] \cup [a,p] \cup [p,r]$, and the nearest point projection
of $x_1$ on $[a,p] \cup [p,r]$. We obtain that 
$$\dist ( x_1, [a,p]) \le 4 \delta _G .$$ 
Since $[a,p]$ 
lies in a $\delta _G$-neighborhood
of $M$, the distance between $M$ and  $x_1$ is smaller than 
$5 \delta _G$. A similar argument applies to $y$ and $y_1$.
Let $N$ be the number of walls $M \subset G$ such that $\dist (e, M) \le 5 \delta _G$.
The property follows letting $L = 2N +1$. 
\end{proof}

\begin{proof}[Proof of Proposition \ref{basis} ]
Since $d$ is a visual metric, there exist constants $a >1$, $C \ge 1$ such that
for every $z, z' \in \partial \Gamma$ one has
\begin{equation} \label{visual}
C ^{-1} a ^{- \ell} \le d(z,z') \le C a ^{- \ell},
\end{equation}
where $\ell$ denotes the distance from $e$ (the identity of $\Gamma$) to 
a geodesic $(z,z') \subset G$. Let $L$ be the constant of the previous
lemma.
Pick $x_1, y_1, y$ in this order 
on the geodesic ray $[e, z) \subset G$, with 
$$\vert x_1 \vert = -\log _a r ,~~ \vert y_1 \vert =  -\log _a r +L,
~~ \vert y \vert = -\log _a r + L + 6 \delta _G .$$
Let $M$ be a wall between $x_1$ and $y_1$
as in the lemma (applied to $e, x_1, y_1, y$). Set
$H := H_+ (M)$ for simplicity. By convexity, for every $z' \in \partial H$ the 
geodesic $(z, z')$ lies in $H$. Therefore $\dist (e, (z,z'))$ is larger
than $-\log _a r - 3 \delta _G$. From (\ref{visual}) one obtains that
$\partial H \subset B(z, C a^{3 \delta _G}r)$.

To establish the other inclusion,
let $z' \in \partial \Gamma$ with $\dist (e, (z,z')) \ge \vert y \vert 
+ 3 \delta _G$. By the same kind of metric argument as in the first part of
the proof
of Proposition \ref{propself}, one sees that $[e, z')$ intersects $B(y, 3 \delta _G)$.
>From the lemma we have $B(y, 3 \delta _G) \subset H$.
Therefore the geodesic ray $[e, z')$ enters $H$.
Since the complementary half-space $H_- (M)$ is convex, the ray $[e, z')$ remains in $H$;
thus $z' \in \partial H$. With (\ref{visual})
it follows that $B(z, C^{-1}a^{-L-9 \delta _G} r) \subset \partial H$.
   
\end{proof}

\subsection{Invariant subsets and parabolic subgroups}

We study now the parabolic subgroups of $\Gamma$ in connection with the
action of $\Gamma$ on its boundary.

\begin{definition}
Let $M$ be a wall of $G$ and let $F$ be a subset of $\partial \Gamma$. We say that 
$M$ \emph{cuts} $F$ if $F$ meets both subsets 
$\partial  H_{-}(M)$
and $\partial H_{+}(M)$. 
\end{definition}

\begin{theorem}\label{1.1} Let $F$ be a  subset
of $\partial \Gamma$ containing
at least two distinct points.  Assume that 
$F$ is invariant under each reflection whose wall 
cuts $F$, and 
let $P$ be the subgroup of $\Gamma$  generated by 
these reflections. Then $P$ 
is a parabolic subgroup of $\Gamma$, the closure of $F$ is the 
limit set of $P$, and $P$ is the stabilizer of $F$
in $\Gamma$.
\end{theorem}

To prove the theorem we introduce the following notion which  will also
be useful in the sequel:

\begin{definition} \label{1.6} Let $F$ be a  subset of 
$\partial \Gamma$.
Assume it contains at least two distinct points. If $\overline{F}
\neq \partial \Gamma$
the \emph{convex hull} $\mathcal C _F$ of $F$ is the intersection of all 
the half-spaces $H$ in $G$ such
that $F$ belongs to the interior of $\partial H$ in $\partial \Gamma$. 
It is a convex subgraph of $G$. If $\overline{F}
= \partial \Gamma$ we set $\mathcal C _F = G$.
\end{definition}

\begin{lemma} \label{convex}
The limit set $\partial \mathcal C _F$ is the closure of $F$. 
In addition every
wall of $G$ which intersects $\mathcal C _F $ cuts $F$. 
\end{lemma}

\begin{proof}
Since $\mathcal C_F $
contains every geodesic in $G$ with both endpoints in $F$,
one has $F \subset \partial \mathcal C _F$,
and thus $\overline F\subset \partial\mathcal{C}_F$. To prove the converse inclusion
pick a point $z \in \partial \Gamma \setminus \overline F$. Since $\overline F$ is a 
closed set,
Proposition \ref{basis} insures the existence of a wall $M$, with associated half-spaces $H_+$ and $H_-$, 
such that
$$z \in (\partial H_+ \setminus \partial M) ~~ \textrm{and} ~~ \partial H_+
\subset (\partial \Gamma \setminus F).$$
 
It follows that  $z \notin \partial H_-$ and that 
$F \subset (\partial \Gamma \setminus \partial  H_+) = \textrm{int} (\partial H_-)$, 
therefore $z \notin \partial \mathcal C _F$.

Let $M$ be a wall intersecting $\mathcal C _F $ and let $(x,y)$ be an edge in the intersection.
Since $\mathcal C _F $ is a closed subset of $G$ it contains $x$ and $y$.
Hence both half-spaces $H_+ , H_-$ bounded by $M$ intersect $\mathcal C _F $.
It follows that $F$ is contained neither in $\textrm{int} (\partial H_+) = \partial \Gamma \setminus \partial H_-$, 
nor in $\textrm{int} (\partial H_-) = \partial \Gamma \setminus \partial H_+$. Thus 
$M$ cuts $F$.
\end{proof}

\begin{proof}[Proof of Theorem \ref{1.1} ]
 Consider the convex hull $\mathcal C _F$ of $F$ in $G$. 
   
The canonical action of $\Gamma$ on $G
\cup \partial \Gamma$ possesses the following property: 
$$ \forall g \in \Gamma,~  g (\mathcal C _{F}) = \mathcal C 
_{g(F)}. $$
Therefore, up to a translation of $F$ by a group element and a conjugation
of $P$, we may assume that $e$ belongs to 
$\mathcal C _F$.

\smallskip

The graph $\mathcal C _F$ is a closed non-empty 
$P$-invariant
convex subset of $G$. At fist we claim that $P$ acts 
simply transitively on it. Since $\Gamma$ acts freely on
$G^0$, $P$ acts freely on $\mathcal C _F ^0$.
To establish the transitivity,
let $x \in \mathcal C _F ^0$ and let
$x_{0} = e, x_1, ..., x_n = x $ be the successive vertices of a geodesic 
segment joining $e$ to $x$ in $G$. Since $\mathcal{C}_F$ is convex, the
edge $(x_k, x_{k+1})$ belongs to $\mathcal{C}_F$ and thus  
the wall $M_k$ between $x_k$ and
$x_{k+1}$ intersects $\mathcal C _F$. According to Lemma \ref{convex}, it 
cuts $F$.
Therefore the reflection along the wall $M_k$ belongs to $P$. In consequence
every $x_i$ belongs to $P$. The transitivity 
follows.

\smallskip

The equality $P = \mathcal C _F ^0 $ and Lemma \ref{convex} imply 
that $\partial P = \overline{F}$. Moreover the stabilizer 
of $F$ in $\Gamma$ stabilizes $\mathcal C _F ^0 $ too, so it coincides with $P$.

\smallskip

Finally we claim that $P$ is equal to the special subgroup $\Gamma _I$ 
where $I$ is the set 
of elements $s \in S$ such that
the walls $M_s$ intersect $\mathcal C _F$. 
The inclusion $\Gamma _I < P$ comes
from Lemma \ref{convex}.
To establish the converse one, let $x \in P$.
We shall prove by induction on 
$\vert x \vert$ that $x$ belongs to $\Gamma _I$. 
When $\vert x \vert =0$ this is obvious.
Let $n \in \mathbb N ^*$, 
assume that the elements $y \in P$ with $\vert y \vert \le n-1$ belong in
$\Gamma _I$, 
and consider $x \in P$ with $\vert x \vert =n$.
Pick a geodesic 
segment $\gamma \subset G$ joining $e$ to $x$, and denote by  
$x_{0} = e, x_1, ..., x_n = x $ its successive vertices. 
Since $x \in \mathcal C _F ^0$, the
convexity of $\mathcal C _F $ implies that $\gamma$ is contained in $\mathcal
C _F $. Therefore the reflection $x_1$ belongs to
$\Gamma _I$. By the induction assumption we get that $x_1 ^{-1} x$ belongs to $\Gamma
_I$, thus $x$ does too. 
\end{proof}

A first corollary concerns a special class of equivalence relations on $\partial \Gamma$.
Examples of such equivalence relations will be considered in Section \ref{equivalence}.

\begin{corollary}\label{1.3}
Consider a $\Gamma$-invariant 
equivalence relation 
on $\partial \Gamma$ whose cosets are connected.
Then:
\begin{enumerate}
\item  The closure of
each coset is either a point or a parabolic limit set.
\item  If a nontrivial coset $F$ is path-connected, and $P$ is the 
parabolic subgroup with $\overline{F}=\D P$, then for every $\eps>0$
and every path $\eta:[0,1]\ra\D P$, there is a path $\eta':[0,1]\ra
F$  such that
$$
d(\eta,\eta')=\max_{t\in [0,1]}d(\eta(t),\eta'(t))<\eps\,.
$$ 
\end{enumerate}
\end{corollary}

\begin{proof}

(1).  If $F$ is a coset, and a
wall $M$  cuts $F$, then the limit set $\D M$ intersects $F$, because
$F$ is connected.  Since the reflection in $M$ is the identity
map on $\D M$, the coset and its image intersect, so they are equal. 
Therefore the assertion follows from
Theorem \ref{1.1}.

(2).  Suppose $F$ is path-connected,  $\overline F=\D P$, $\eta\subset\D P$
is a path, and $\eps>0$.  By Proposition \ref{basis}, we can find a finite
collection of half-spaces $H_1,\ldots,H_k\subset G$, such that the
limit sets $\D H_1,\ldots,\D H_k\subset \D\Gamma$ each have diameter $<\eps$,
and their interiors cover the image of $\eta$.  We can then choose
$0=t_0<t_1<\ldots<t_n=1$ such that for each $i$, the pair $\eta(t_{i-1})$, $\eta(t_i)$
is contained in $\D H_{j_i}$ for some $j_i\in\{1,\ldots,n\}$,  and 
$\diam(\eta([t_{i-1},t_i]))<\eps/2$.  Since $F$ is dense in $\D P$, for each
$j\in \{1,\ldots,n\}$ we
may choose $s_j\in F$ close enough to $\eta(t_j)$, such that for all $i$,
the pair $s_{i-1}$, $s_i$ also lies in $\D H_{j_i}$.  By the path 
connectedness of $F$, we may join $s_{i-1}$ to $s_i$ by a path $\bar \ga_i
\subset F$; reflecting the portion of $\bar\ga_i$ lying outside
$\D H_{j_i}$ into $\D H_{j_i}$ using the reflection whose wall bounds
$H_{j_i}$, we get a path $\ga_i\subset \D H_{j_i}$ joining $s_{i-1}$ to 
$s_i$.  Concatenating the $\ga_i$'s, we obtain the desired path $\eta'$. 
\end{proof}

The same proof as for part (1) of the previous corollary gives :

\begin{corollary} Let $\Phi$ be a quasiconvex subgroup
of $\Gamma$ with connected limit set.  Assume 
that for all $g \in \Gamma$, the intersection 
$g(\partial \Phi) \,\cap\,\partial \Phi$ is either 
empty, or equal to $\partial \Phi$.
Then  $\Phi$ is virtually a parabolic subgroup of $\Gamma$. In other
words there exists a parabolic subgroup $P$ of $\Gamma$
such that $\Phi \cap P$ is of finite index in both subgroups
$\Phi$ and $P$.
\end{corollary}

\begin{corollary} \label{1.5} Each connected component of $\partial \Gamma$ 
containing more than one point is a parabolic limit set.
\end{corollary}

\subsection{Shadowing curves} \label{shadow}

In this paragraph we present a sort of a quantitative version of Theorem \ref{1.1}.

\medskip

 Let $F$ be a subset in $\partial \Gamma$.
As before we denote by  $\mathcal C _{F}$ its convex hull in $G$
(see Definition \ref{1.6}). The canonical action of $\Gamma$ on $G
\cup \partial \Gamma$ possesses the following property : 
$$ \forall g \in \Gamma,~  g (\mathcal C _{F}) = \mathcal C 
_{g(F)}. $$
Therefore, up to a translation of $F$ by a group element, 
we are in the situation where $\mathcal C _{F}$ 
contains $e$ the identity
of $\Gamma$.

\begin{definition} Let $\gamma$ be a non-constant curve in $\partial \Gamma$,
let $I$ be a non-empty subset of $S$ and let $L \ge 0$. We say that
$\gamma$ is a \emph{$(L, I)$-curve} if 
$e \in \mathcal C _{\gamma}$,
and if for every $s \in I$
there exists an edge $a_s$ of type $s$ in $\mathcal C _{\gamma}$
with $\dist (e, a_s ) \le L$.
\end{definition}

\begin{proposition} \label{1.8} 
Let $\epsilon > 0$ and let $L , I$  be as in the above
definition. Let $P \leqq \Gamma$ be a conjugate of $\Gamma _ I$,
and let $\eta$ be a (parametrized) curve contained in $\partial P$.
There exists a finite subset $E \subset \Gamma$ such that for any $(L,I)$-curve $\gamma$
the subset $\bigcup _{g \in E} g \gamma \subset \partial \Gamma$ contains
a curve which approximates
$\eta$  to within $\epsilon$ with respect to the ${C} ^0$ distance.
\end{proposition}

\begin{proof} 
 Since $P = h \Gamma _I h^{-1}$ for some $h \in \Gamma$, we have 
$\partial P = h(\partial \Gamma _I)$. Moreover $h$ is a bi-Lipschitz
homeomorphism of $\partial \Gamma$. Therefore it is enough to establish
the proposition for $P = \Gamma _I$.

\smallskip

\emph{First step} :  we show that for every $(L,I)$-curve $\gamma$
the subset 
$$\bigcup _{ \{g \in \Gamma ~; ~\vert g \vert \le L\} } g \gamma$$ 
of $\partial \Gamma$ contains a curve passing through every 
 $\partial M_s$ with $s \in I$. 
For this purpose pick $s \in I$ , let $a_s \subset 
\mathcal C _{\gamma}$ be an edge of type $s$ with 
$\dist (e, a_s ) \le L$, and let $M$ be the wall containing $a_s$.
Let $g \in H_- (M)$ be such that $a_s =  (g, gs)$. The geodesic
segment $[e, g]$ belongs to $\mathcal C _{\gamma}$ (by convexity).
Denote by 
$$g_0 =e,~ g_1 = \sigma_1,~ g_2 = \sigma_1 \sigma_2,~ ...~ ,~ 
g_n = g = \sigma_1 ... \sigma_n , $$
the successive vertices of the segment $[e, g]$ with $\sigma_i \in S$.
According to Lemma \ref{convex} the wall $M_i$ passing between $g_i$ and 
$g_{i+1}$ cuts $\gamma$.
Thus for every $i \in
\{0, ..., n-1\}$ the curves $\gamma$ and $g_i \sigma _{i+1} g_i ^{-1} \gamma$ 
intersect. 
One has $g_i \sigma _{i+1} g_i ^{-1} = g_i g_{i+1} ^{-1}$. 
It follows that 
the subset $\gamma \cup g_1 ^{-1} \gamma \cup ... \cup g_n ^{-1} \gamma$
is an arcwise connected set. It intersects the limit set of
the wall $g^{-1} (M) = M_s$.
Thus the subset 
$$\bigcup _{ \{g \in \Gamma ~; ~\vert g \vert \le L \} } g \gamma $$
is an arcwise connected set which intersects the limit set of every wall 
$M_s$ with $s \in I$.

\smallskip

\emph{Second step} : Consider  a collection $H_1, ... , H_k$ of
half-spaces in $G$ all of them intersecting the special subgroup
$\Gamma _I$ \emph{properly} \emph{i.e.} $\Gamma _I \cap H_i  $
is distinct from $\emptyset$ and
$\Gamma _I $.
We will show that there exists a finite 
subset $E_0$ of $\Gamma$ such that for any $(L,I)$-curve $\gamma$,
one can find in the subset $\bigcup _{g \in E_0 } g \gamma \subset \partial \Gamma$
a curve passing through every $\partial H_1, ... , \partial H_k$.
To do so, pick for each $i \in \{1, ... , k\}$ an element  $p_i \in \Gamma _I$ 
adjacent to the wall which bounds $H_i$. Let $c$ be a path in $G$ 
which joins successively 
$p_1, ... , p_k$ and whose vertices $c_1, ... , c_n$ 
lie in $\Gamma _I$. 
Define
$$ E_0 = \{c_i g \in \Gamma ~; ~\vert g \vert \le
L,~ 1 \le i \le n  \}, $$
and let $\theta$ be a curve made of translates of $\gamma$ passing through 
every $\partial M_s$ with $s \in I$ (as constructed in step 1). The subset 
$\bigcup _{1 \le i \le n} c_i \theta$ meets the limit set of any wall
$c_i (M_s)$  with $i \in \{1, ... , n\}$ and $ s \in I$.
 In particular it
intersects $\partial H_1, ... , \partial H_k$.
In addition this is an arcwise connected set. Indeed write $c_{i + 1} = c_i s
= \sigma c_i$, with $s \in I$ and $\sigma = c_i s c_i ^{-1}$. 
Then $c_{i + 1} \theta =  \sigma c_i \theta$. The curve $\theta$
intersects $\partial M_s$, thus  $c_{i} \theta$
intersects $c_i (\partial M_s)$. The intersection set is pointwise invariant by 
the reflection $\sigma$ and thus it belongs to $c_{i + 1} \theta$ too.

\smallskip

\emph{Last step} : We finally prove the proposition. By Proposition \ref{basis}
there exists a collection of half-spaces $H'_1, ... , H'_{k+1}$
of $G$ such that the union of their limit sets is a neighborhood
of $\eta$ contained in the $\epsilon/2$-neighborhood of $\eta$.
Reordering if necessary we may assume that the curve $\eta$ enters
successively $\partial H'_1, ... , \partial H'_k$.
Pick a collection of half-spaces $H_1, ... , H_k$  each of them intersecting
 $\Gamma _I$ properly, and such that for every 
$i \in \{1, ... , k\}$ one has
$$ \partial H_i \subset \partial H'_i \cap \partial H'_{i+1}.
$$
Existence of the $H_i$'s follows from Proposition \ref{basis}. 
According to step 2 there exists
a subset $E_0$ of $\Gamma$ such that for every $(L,I)$-curve $\gamma$,
one can find in the subset $\bigcup _{g \in E_0 } g \gamma$
a curve passing through every $\partial H_1, ... , \partial H_k$. 
Denote by $\theta$ such a curve. The part of $\theta$ between 
$\partial H_{i-1}$ and $\partial H_i$ may exit from $\partial H'_i$.
If it happens reflect the outside part of $\theta$ along the wall
which bounds $H'_i$. The resulting curve 
can be parametrized to approximate $\eta$  to within $\epsilon$ with respect
to the 
${C} ^0$ distance.
Let $\sigma_i$
be the reflection along the wall of $H'_i$. The following subset of $\Gamma$ 
$$ E = E_0 \cup \bigcup _{i=2} ^k \sigma_i (E_0) $$
satisfies the property we were looking for. 
\end{proof}

We now establish the abundance of $(L,I)$-curves. Denote
by $N_r (E)$ the open $r$-neighborhood in $\partial \Gamma$
of a subset $E \subset \partial \Gamma$. 

\begin{proposition} \label{1.9} Let $I \subset S$ and let
$P \leqq \Gamma$ be a conjugate of $\Gamma _I$. 
For all $r > 0 $ there
exist $ L \ge 0 $ and $\delta > 0 $ such that 
every curve $\gamma \subset \partial \Gamma$ satisfying the following
conditions is a $(L,I)$-curve :
\begin{itemize}
\item [(i)] Its convex hull contains $e$,
\item [(ii)] $\gamma \subset N_{\delta} (\partial P)$, 
\item [(iii)] $\gamma \nsubseteq N_{r} (\partial Q)$ for any parabolic 
$Q \lneqq P$ with connected limit set.
\end{itemize}

\end{proposition}

\begin{proof}
Suppose by contradiction  that for every $ L \ge 0 $ and $\delta > 0 $
there exists a curve $\gamma$ which satisfies property (i), (ii) and (iii), and
which is not a $(L,I)$-curve. Choose $L = n$, $\delta = 1/n$
and pick $\gamma _n$ a corresponding curve. 
We may assume,
by extracting a subsequence if necessary, that 
there exists an element $s \in I$ such that
for every $n \ge 1$ no edge of 
$\mathcal C_{\gamma _n} \cap B(e , n)$  is of type $s$. 
 We can also suppose that the sequence of compact subsets 
$ \{ \gamma _n \}_{n   \ge 1}$
converges with respect to the Hausdorff distance to a non-degenerate continuum
$\mathcal L \subset \partial P$. Indeed this follows from the compactness of
the set of continua in a compact metric space, equipped with the Hausdorff 
distance
(see \cite{Mun} p. 281).  
In addition, by a standard diagonal argument, 
we may assume the sequence $\{\mathcal C_{\gamma _n} \}_{n   \ge 1}$
converges to a convex subset $\mathcal C \subset G$
on every compact subset of $G$. 
With item (i) one has $e \in \mathcal C$. Moreover one sees easily that 
$\mathcal L \subset \partial \mathcal C$.
The fact that no edge of $\mathcal C $ is of type $s$ implies 
that $\mathcal L$ is contained in the limit set of the special
subgroup generated by $S \setminus \{s \}$. 
Intersections of parabolic subgroups are again parabolic subgroups 
(see \cite{D} Lemma 5.3.6); 
thus $\mathcal L$ is contained in the limit set of a proper parabolic subgroup
of $P$. By Corollary \ref{1.5} the connected component which contains
$\mathcal L$ is the limit set of a proper parabolic subgroup $Q$ 
of $P$.
So we get a contradiction with the hypothesis (iii).
\end{proof}

\subsection*{Remarks and questions~: } 1) Theorem \ref{1.1} admits a partial converse.
Indeed let $P = g \Gamma _I g^{-1}$ be a parabolic subgroup of $\Gamma$ and let $M$ be a wall 
of $G$ such that $\partial P$ meets both open subsets 
$\partial  H_{-}(M) \setminus \partial M$
and $\partial H_{+}(M) \setminus \partial M$. Then a convexity argument
shows that $M$ admits an edge whose end-points lie in $g \Gamma _I$. 
Thus letting $\sigma$ be the reflection along $M$, there exist $h, h' \in \Gamma _I$
such that $\sigma g h = g h'$. It follows that  $\sigma = g h' h^{-1} g^{-1} \in P$,
and so $\partial P$ is $\sigma$-invariant.

\medskip 

2) Let $\Gamma$ be an arbitrary hyperbolic group. The limit set of the intersection of 
two quasiconvex subgroups is the intersection of their limit sets (see \cite {G} p. 164). 
Hence for any quasiconvex subgroup $\Phi \leqq \Gamma$
the following properties are equivalent :
\begin{itemize}
\item [--]
 for every $g \in \Gamma, \ \ g \partial \Phi \cap \partial \Phi = \partial \Phi \ \ \textrm{or}
\ \ \emptyset $,
\item [--] for every $g \in \Gamma$, either $g \Phi g^{-1} \cap \Phi$
is finite or is of finite index in both subgroups $g \Phi g^{-1}$ and
$\Phi$.
\end{itemize}

3) Given a subset $E \subset \partial \Gamma$ there exists a unique smallest parabolic
limit set $\partial P$ containing $E$, moreover if  $E$ is connected and non-reduced to a point 
then $\partial P$ is so.
Indeed this follows from the fact that parabolic subgroups are stable by intersection (see 
\cite{D} Lemma 5.3.6),
from the property of intersections of limit sets recalled in Remark 2 above, and from 
Corollary
\ref{1.5}.

\medskip

4)
Let $\Gamma$ be a hyperbolic group and consider a closed
$\Gamma$-invariant equivalence relation on $\partial \Gamma$ 
whose cosets are continua. 
What one can say about such equivalence relations ? 
In particular for which groups $\Gamma$ do the nontrivial cosets 
arise as the limit sets
of a finite collection of conjugacy classes of quasiconvex
subgroups ?

Note that the quotient space of $\partial \Gamma$ by a  
$\Gamma$-invariant closed equivalence relation $\sim$ is a compact metrizable  
space on which $\Gamma$ acts as a convergence group, \emph{i.e.}
$\Gamma$ acts properly discontinuously on the set of triples of distinct
points of $\partial \Gamma / \!\sim$ (see \cite{Bow}).

\section{Combinatorial modulus and Coxeter groups}\label{cox}

In this section we study the combinatorial modulus
on boundaries of hyperbolic Coxeter groups.
Consider a hyperbolic Coxeter group $\Gamma$ with connected boundary, and let
$Z$ be the metric space $\partial \Gamma$ equipped with a self-similar
metric (see Definition \ref{self}).
We fix in the sequel some $\kappa$-approximation 
$\{ G_k \} _{k \in \hn}$ of $Z$.

\medskip

 Let $d_0$ be a fixed (small) positive constant. As before $\mathcal
F _0$ denotes the set of curves $\gamma \subset Z$ with $\diam (\gamma) \ge d_0$. 
We wish to establish a kind of filtration of $\mathcal F _0$  
by ``elementary curve families'' (see the discussion after Corollary \ref{comparable}). 
To this aim we introduce the following families of curves. 
Let $\partial P \subset Z$ be a parabolic limit set and let $\delta, r > 0$.
Denote by $ \mathcal F _{ \delta, r } (\partial P)$
the family of curves $\gamma$ in $Z$ satisfying the following conditions:
\begin{itemize}
\item [--] $\gamma \subset N_{\delta} (\partial P)$ and  $ \diam \gamma \ge d_0$, 
\item [--] $\gamma \nsubseteq N_{r} (\partial Q)$ for any connected parabolic limit set 
$\partial Q \subsetneqq \partial P$ .
\end{itemize}
Notice that for $\delta$ small enough  there is only a finite number of
$\partial P$ such that $\mathcal F _{ \delta, r } (\partial P) \neq
\emptyset$ (see Proposition \ref{fls}).

For $\epsilon > 0$ and for any (parametrized) curve $\eta$ in $Z$, 
let $\mathcal U _{\epsilon} (\eta)$ be the set of curves whose $C^0$
distance to $\eta$ is smaller than $\epsilon$.

\begin{theorem}\label{3.4}
There exists a positive increasing function $\delta _0$ on
$(0, +\infty)$, that depends on $d_0$ only, and which possesses
the following property. Let $p \ge 1$, $\epsilon$, $r >0$, let $\partial P \subset Z$
be a parabolic limit set and let $\eta$
be any curve in $\partial P$. 
There exists
$C = C(p, \epsilon, \eta, r, d_0)>0$ such that 
for every $\delta \le \delta _0 (r)$ and for every $k \in \hn $ : 
$$\Mod _p (\mathcal F _{ \delta, r } (\partial P) , G_k)
\le C \Mod _p (\mathcal U _{\epsilon} (\eta) , G_k).$$
In addition when $p$ belongs to a compact subset of $[1, +\infty)$
the constant $C$ may be choosen independent of $p$. 
\end{theorem}

In combination with  Proposition \ref{3.1}.1, this leads to :
 
\begin{corollary} \label{comparable}
Let  $\eta \in \mathcal F_0$ 
and let $\partial P$ be the smallest parabolic limit set containing $\eta$. 
Let $r > 0$ be small enough in order
that $\eta \nsubseteq \overline{N} _{r} (\partial Q)$ for any parabolic limit set 
$\partial Q \subsetneqq \partial P$.  
Let $\delta \le \delta _0 (r)$ be as in the previous theorem. 
Let $\epsilon > 0$ be small enough in order that
$\mathcal U _{\epsilon} (\eta)
\subset \mathcal F _{ \delta, r } (\partial P)$. 
Then for every $k \in \hn$ one has :
$$\Mod _p (\mathcal U _{\epsilon} (\eta) , G_k) 
\le \Mod _p (\mathcal F _{ \delta, r } (\partial P) , G_k)
\le C \Mod _p (\mathcal U _{\epsilon} (\eta) , G_k),$$
where $C$ is the constant defined in the previous theorem.
\end{corollary}

 Before going into the proof of the theorem, 
we discuss the meaning of these
results and the roles and the 
interdependence of the parameters $r, \delta, \epsilon$. Corollary
\ref{comparable} means that for every curve $\eta \in \mathcal F_0$ a
sufficiently small neighborhood $\mathcal U _{\epsilon} (\eta)$ 
behaves -- for the
combinatorial modulus -- like 
one of curve families $\mathcal F _{ \delta, r } (\partial P)$ with $\delta
\le \delta _0 (r)$ . 
In particular the behaviour of 
$\Mod _p (\mathcal U _{\epsilon} (\eta) , G_k)$
is independent of $\epsilon$ (small enough depending on $\eta$).
Theorem \ref{3.4} states that for all $r>0$ and all $\delta \le 
\delta _0 (r)$ the modulus of $\mathcal F _{ \delta, r } (\partial
P)$
is controlled by above by the modulus of any neighborhood $\mathcal U
_{\epsilon} (\eta)$ of any curve $\eta \subset \partial P$.
To complete the picture we notice that, given $r_0 >0$, the curve family 
$\mathcal F _0$ can be express as a finite union
of curve families $\mathcal F _{ \delta _i, r_i } (\partial P_i)$
with $r_i \le r_0$ and $\delta _i\le \delta _0 (r_i)$. To do so, one first defines 
the \emph{height} of a connected parabolic limit set $\partial P$
to be the largest integer $n \ge 0$ such that there exists a chain 
of connected parabolic limit sets of the form
$$
\partial Q_0  \subsetneqq \partial Q_1 \subsetneqq \dots \subsetneqq
\partial Q_n = \partial P.
$$
Then one considers the finite collection
of the connected parabolic limit sets $\partial P_i$ of diameter 
larger than $\frac{d_0}{2}$, and defines the corresponding parameters 
$(r_i, \delta _i)$
by induction on the height
of $\partial P_i$.

\medskip
 
We now give the

\begin{proof}[Proof of Theorem \ref{3.4}]
Instead of considering all curves with diameter larger than $d_0$,
we may restrict ourselves to those whose convex hulls contain $e$
(see the discussion at the beginning of Subsection \ref{shadow},
and notice that 
$\sup _{\gamma \in \mathcal F _0} \dist (e, \mathcal C_\gamma) < +\infty$).
In the sequel of the proof we make this restriction and 
we keep the same notation for the restriction of 
$ \mathcal F _{ \delta, r } (\partial P)$. We shall also use --  without
further mention -- the metric equivalence between $Z$ and $\partial \Gamma$. 

\smallskip

Let $I \subset S$ such that $P$ is a conjugate of $\Gamma _I$.
Thanks to Proposition \ref{1.9} there exists $L$ and $\delta  > 0$
depending only on $r$, 
such that every element of $ \mathcal F _{ \delta, r } (\partial P)$
is a $(L, I)$-curve. By Proposition \ref{1.8}
there exists a finite subset $E \subset \Gamma$ 
such that for every $\gamma \in \mathcal F _{ \delta, r } (\partial P)$
the subset $\bigcup _{g \in E} g \gamma \subset Z$ contains a curve of 
$\mathcal U _{\epsilon} (\eta)$. 
Let $G' _k$ be the incidence graph of the covering of $Z$ by the elements
of $\bigcup _{g \in E} g^{-1} G^0 _k $. According to Proposition 
\ref{3.2}, and  since $\Gamma$ acts on $Z$ by bi-Lipschitz 
homeomorphisms (Proposition \ref{propself}), the moduli $\Mod_p (\cdot , G'_k )$
and $\Mod_p (\cdot , G_k )$ are comparable.
For all $\mathcal U _{\epsilon} (\eta)$-admissible function 
$\rho : G^0 _k \to \hr_+$  define a function 
$\rho ': G'_k {}^0 \to \hr_+$ by
$$ \rho '(g^{-1}v) = \rho (v), $$ 
where $v \in G^0 _k$ and $g \in E$.  We claim that $\rho '$ is 
$\mathcal F _{ \delta, r } (\partial P)$-admissible. Indeed let $\gamma
\in \mathcal F _{ \delta, r } (\partial P)$, and let $\theta \in  
\mathcal U _{\epsilon} (\eta)$ so that $\theta \subset 
\bigcup _{g \in E} g \gamma $. One has :
$$\ell _{\rho '} (\gamma) = \sum _{g \in E} \ell _\rho (g \gamma) 
\ge \ell _\rho (\theta) \ge 1.$$
The claim follows.
Thus 
$$\Mod _p (\mathcal F _{ \delta, r } (\partial P) , G' _k)
\le \vert E \vert \sum _{v \in G^0 _k } \rho (v)^p ,$$
and so
$$\Mod _p (\mathcal F _{ \delta, r } (\partial P) , G' _k)
\le \vert E \vert \Mod _p (\mathcal U _{\epsilon} (\eta) , G_k).$$ 
\end{proof}

We now present a companion result to Theorem \ref{3.4}. Its statement requires some 
notation. 
For a parabolic subgroup $P \leqq \Gamma$ let $\mathcal {N} _r (\partial P)$ be the family
of curves $\gamma \subset N_r (\partial P)$ with $\diam \gamma \ge d_0$.
Let $\mathcal L$ be a collection of parabolic limit sets. 
For $r>0$ set 
$$\mathcal {F} _r (\mathcal L) := \mathcal {F} _0 \setminus \bigcup _{\partial Q \in \mathcal L} 
\mathcal {N} _r (\partial Q). $$
We denote by $\Confdim (\partial P)$ the Ahlfors regular conformal dimension of $\partial P$ (its 
definition is recalled in
Remark 2 at the end of Section \ref{hyp}), and we set 
$$\Confdim (\mathcal L) := \max _{\partial Q \in \mathcal L} \Confdim (\partial Q). $$

The following property ``controls'' the modulus of neighborhoods of some
limit sets by the modulus of the complementary subsets.

\begin{theorem} \label{theocircular} Let $\mathcal L$ be a $\Gamma$-invariant collection of 
connected proper parabolic limit sets. 
Assume that there exists  $\partial P \in \mathcal L$ such that for every $\partial Q \in \mathcal L$ with 
$\partial Q \neq \partial P$, 
the set $\partial P \cap \partial Q$ is totally disconnected or empty.
Let $p > \Confdim (\partial P)$. 
There exist constants $C>0$ and $a \in (0,1)$ such that for $r>0$ small enough and 
for every $k \in \hn $ one has :
$$\Mod _p (\mathcal {N} _r (\partial P) , G _k) \le C \cdot
\sum _{\ell = 0}^{k} \Mod _p (\mathcal {F} _r (\mathcal L) , G _{k - \ell} )~ a^{\ell}.$$
\end{theorem}

\begin{corollary} \label{corcircular} Let $\mathcal L$ be a $\Gamma$-invariant collection of 
connected proper parabolic limit sets. 
Assume that for every $\partial P, \partial Q \in \mathcal L$ with $\partial Q \neq \partial P$, 
the set $\partial P \cap \partial Q$ is totally disconnected or empty.
Let $p > \Confdim (\mathcal L)$. 
There exist constants $C>0$ and $a \in (0,1)$ such that for $r>0$ small enough and 
for every $k \in \hn $ :
$$\Mod _p (\mathcal {F} _0 , G _k) \le C \cdot
\sum _{\ell = 0}^{k} \Mod _p (\mathcal  {F} _r (\mathcal L) , G _{k - \ell} )~ a^{\ell}.$$
\end{corollary}

\begin{proof} [Proof of the Corollary \ref{corcircular}]
 Notice that $\mathcal F _0 = \mathcal {F} _r (\mathcal L) \cup \bigcup
_{\partial P \in \mathcal L} \mathcal {N} _r (\partial P)$. Let $r>0$ be such
that $d_0 - 2r >0$.  Only finitely many parabolic limit sets 
have diameter at least  $d_0 - 2r$ (see Proposition \ref{fls}).
Hence, all but a finite number of the curve families 
$\mathcal {N} _r (\partial P)$
are empty. Theorem \ref{theocircular} applied to the non-empty ones, in
combination
with the basic properties of the combinatorial modulus, completes the proof.
\end{proof}

As an example the collection $\mathcal L$ of all circular limit sets satisfies the above hypotheses with
$\Confdim (\mathcal L) = 1$.

\medskip

The proof of the theorem relies on some lemmata. 
The first one shows that two connected parabolic limit sets whose intersection
is at totally disconnected or empty,
have to ``move away one from another linearly''.

\begin{lemma}\label{C_0}
There is a constant $C_0 \ge 1$ with the following property. Let $ \partial P_1, \partial P_2$
be two connected parabolic limit sets such that either $ \partial P_1 = \partial P_2$
or $\partial P_1 \cap \partial P_2$ is totally disconnected or empty. 
Assume that there exist $r>0$ and a curve $\gamma 
\subset N_r(\partial P_1) \cap N_r(\partial P_2)$ with
$\diam (\gamma)\geq C_0r$. Then $\partial P_1 = \partial P_2$.
\end{lemma}

\begin{proof} Applying self-similarity property \ref{homothety} we may assume that 
$\gamma$ has unit diameter.  Then $\partial P_1$ and $\partial P_2$ 
have diameter bounded
away from zero, and hence $P_1$ and $P_2$ belong to a 
finite collection of parabolic
subgroups with connected limit sets (Proposition \ref{fls}).  But if $P_1 \neq P_2$, then 
$\partial P_1 \cap \partial P_2$ cannot contain an arc, because the 
intersection is totally disconnected.  Hence by choosing $C_0$ large enough, 
the lemma follows. 
\end{proof}

>From now on we consider a connected parabolic limit set $\partial P \subset Z$.  

\begin{lemma} Let $r>0$, $C \ge 1$ and let $\gamma \subset N_r (\partial P)$
be a curve with $\diam (\gamma) > C r$.
 For every $z \in (\gamma \setminus \partial P)$ there exists  a subcurve
$\gamma'$ of $\gamma$ such that letting $d:= \diam (\gamma')$, one has :
$$ \gamma' \subset B(z, 4d) \cap N_{d/C}(\partial P)
~\textrm{and}~ \gamma' \nsubseteq N_{d/(8C)}(\partial P).$$
\end{lemma}

\begin{proof} Either there exists a point $z_1$ on $\gamma$ with 
$d(z,z_1)  < 2 C \dist (z, \partial P)$ and $\dist (z_1 , \partial P)
= 2 \dist (z,\partial P)$, or there does not. 

If not then the part of $\gamma$ contained 
 in the ball centered
on $z$ of radius $ 2 C \dist (z, \partial P)$ lies in the
$( 2 \dist (z,\partial P))$-neighborhood of $\partial P$.  
The assumption on $\diam (\gamma)$ implies that $\gamma$
exits the ball $B(z, 2 C \dist (z, \partial P))$. Thus $\gamma$ contains
a desired subcurve which moreover contains $z$.

If $z_1$ exists then we repeat the process with $z_1$ instead of $z$.
Eventually we get a sequence of points $z_0, z_1, ... , z_n$ on $\gamma$
with $z_0 = z$ such that for all indices $i $ :
$$\dist (z_i , \partial P) = 2^i \dist (z,\partial P)~\textrm{and}~
d(z_i,z_{i+1})  < 2 C \dist (z_i, \partial P).$$
Since $ \gamma \subset N_r (\partial P)$ and $\dist (z,\partial P) >0$ the process has to stop,
let $z_n$ be the last point.
>From the first case argument we get a subcurve $\gamma'$ of $\gamma$ 
such that letting $d=\diam (\gamma')$ we have :  
$$ z_n \in \gamma' \subset B(z_n, 2  C \dist (z_n , \partial P) )
~, ~d \ge 2  C \dist (z_n , \partial P),$$
and 
$$\gamma' \subset N_{d/C}(\partial P)
~,~ \gamma' \nsubseteq N_{d/(8C)}(\partial P).$$
Finally we compute
$$ d(z, z_n) \le \sum _{i=0} ^{n-1} d(z_i, z_{i+1}) \le 
2^{n+1} C \dist (z,\partial P) =  2 C \dist (z_n , \partial P) \le d.$$
Therefore $\gamma' \subset B(z, 4d)$. 
\end{proof}

Pick for every $\ell \in \hn$
a collection $\mathcal B _{\ell} $ of balls in $Z$ centered
on $\partial P$ and of radius $2 ^{- \ell }$ 
such that the set $\{\frac{1}{2} B ~;~ B \in \mathcal B _{\ell} \}$
is a minimal covering of 
$\partial P$.  Let $ \mathcal B = \bigcup _{\ell \in \hn} \mathcal B _{\ell} $.
Recall that the radius of a ball 
$B$ is denoted by $r(B)$.
 From the previous lemma we get :

\begin{lemma}\label{subcurve}
Let $r>0, C \ge 1$ and let $\gamma \subset N_r (\partial P)$
be a curve with $\diam (\gamma) > C r$. For every
$z \in (\gamma \setminus \partial P)$ there exists a ball $B \in \mathcal B$
and a subcurve $\gamma'$ of $\gamma$ such that letting $d = \diam (\gamma')$
we have :
$$ \gamma' \cup \{z\} \subset B ~,~  r(B) \le 36 d ~,~ 
\gamma' \subset  N_{d/C}(\partial P)
~\textrm{and}~ \gamma' \nsubseteq N_{d/(8C)}(\partial P).$$
\end{lemma}

\begin{proof}
Consider the subcurve $\gamma'$ obtained in the previous lemma.
We compute :
$$ \dist (z, \partial P) \le 4 d + \dist (\gamma', \partial P) 
\le 5d.$$
Let $w$ be a point in $\partial P$ which realises $\dist (z, \partial P)$,
we have $\gamma' \cup \{z\} \subset B(w, 9d)$.
Pick  a ball $B \in \mathcal B$ such that $w \in \frac{1}{2} B$ and
$r(B)/4 \le 9d \le r(B)/2$. Then $B$ contains $B(w, 9d)$ and hence 
it contains $\gamma' \cup \{z\}$. In addition we have $r(B) \le 36 d$.
\end{proof}

\begin{proof}[Proof of Theorem \ref{theocircular}]
Let $\partial P \in \mathcal L$ be the connected proper limit set of the statement,
and let $p > \Confdim (\partial P)$.
For appropriate $r>0$ we wish to construct 
on every scale a $\mathcal {N} _{r} (\partial P)$-admissible
function $\rho_k : G_k ^0 \to \hr _+$ with controlled $p$-mass. 
For this purpose we equip $\partial P$
with an Ahlfors regular metric $\delta _{\partial P}$, which is quasi-Moebius 
equivalent\footnote{Two metrics $d_1, d_2$ on a space $X$ are said to be
\emph{quasi-Moebius 
equivalent} if the identity map $(X,d_1) \to (X,d_2)$ is a quasi-Moebius
homeomorphism. } to 
the restriction of the $Z$-metric to $\partial P$,
and whose Hausdorff dimension $Q$ satisfies $Q <p$. Let $\mu$
be the $Q$-Hausdorff measure of $(\partial P, \delta _{\partial P})$.
Given any ball  $B \subset Z$ centered on $\partial P$,  
its trace on $\partial P$ is \emph{comparable} to a ball
of $(\partial P, \delta _{\partial P})$.  
In other words there are balls $B'_1, B'_2 \subset (\partial P, \delta _{\partial P})$
of radii $\tau _1, \tau _2$, such that :
$B'_1 \subset (B \cap \partial P) \subset B'_2$, with $\frac{\tau _2}{\tau _1}$
smaller than an absolute constant.  
We denote by  $\tau (B)$ the radius of a minimal ball
of $(\partial P, \delta _{\partial P})$ containing $B \cap \partial P$. One has
$$ \mu (B \cap \partial P) \asymp \tau (B)  ^Q . $$

Let $\mathcal B$ be a ball family as considered in the statement of the previous lemma.
We inflate
every element $B
\in \bigcup _{\ell = 0} ^{k} \mathcal B _{\ell} $ 
to essentially unit diameter as in Definition \ref{homothety}.
Let $g_B$ be the corresponding group element. 
Define $G_k \cap B$ to be the incidence graph of the covering
of $B$ by the subsets $v \in G_k ^0$ with
$v \cap B \neq \emptyset$. Note that for $B \in \mathcal{B} _{\ell}$,
the graph $G_k  \cap B$ 
may be considered as a $(2 \kappa L_0)$-approximation of $g_B (B)$ on scale $k-\ell $
(via the group element $g_B$).

\smallskip

For $k \in \hn$ and $r>0$ write 
$m _k := \Mod _p (\mathcal {F} _{r } (\mathcal L) , G _k)$ for simplicity. 
Suppose that $r$ is small compared to $L_0, C_0, d_0$
($L _0$ and $C _0$ are defined respectively in statements \ref{homothety}
and \ref{C_0}) .
Then using  the group
element $g_B$ and the $\Gamma$-invariance of $\mathcal L$, we may  pull back to $G_k \cap B$ 
a normalized  
minimal $\mathcal {F} _{r } (\mathcal L)$-admissible function on scale $k-\ell $, 
in order to get for every ball
$B \in \bigcup _{\ell = 0} ^{k} \mathcal B _{\ell} $ a function  
$\rho _{B} : (G_k \cap B)^0 \to \hr _+$, with the following properties 
\begin{itemize}
\item [(i)] its $p$-mass satisfies, with $D \ge 1$ independent of $k$ and of $B$ :
$$M_{p} (\rho _{B} ) \le D \cdot m_{k- \ell} \cdot \tau (B) ^p ,$$ 
\item [(ii)] every curve $\gamma \subset B$ 
whose diameter is larger than $r(B)/36$
picks up a $(\rho _{B})$-length larger than  $\tau (B)$, unless it lies 
in the $(10^{-3}C_0 ^{-1} r(B))$-neigh\-borhood of a limit set of $\mathcal L$. 
\end{itemize}
Define $\rho_k$ by
$$\forall v \in G_k ^0, ~ \rho _k (v) = \max \rho _{B} (v),$$
where the maximum is over all 
$B \in \bigcup _{\ell = 0} ^{k} \mathcal B _{\ell} $ 
with $v \cap B \neq \emptyset$.
Its $p$-mass is linearly bounded by above by  
$$\sum _{\ell = 0}^{k} \sum _{B \in \mathcal B _{\ell}}
M_p (\rho _{B}).$$ 
Since the metric $\delta _{\partial P}$ and the rectriction of the $Z$-metric
are quasi-Moebius
equivalent, they are H\"older
equivalent (see \cite{He}). In particular there is $\alpha >0$ such that  
for every $B \in \mathcal B_\ell$ one has $\tau (B) \lesssim 2^ {- \alpha \ell}$.
With item (i) one obtains that the $p$-mass of $\rho _k$ is 
linearly bounded by above by :
$$\sum _{\ell = 0}^{k} m_{k-\ell} \sum _{B \in \mathcal{B} _{\ell} } 
\tau (B) ^{p}
\le \sum _{\ell = 0}^{k} m_{k-\ell} \sum _{B \in \mathcal{B} _{\ell} } 
\tau (B) ^{Q} a^ {\ell}
\asymp \mu (\partial P)  \sum _{\ell = 0}^{k} m_{k-\ell}~ a^{\ell},$$
with $a = 2^{- \alpha (p-Q)}$.
Because $p>Q$ one has $0< a <1$ .  

\smallskip

It remains to prove that $\rho _k$ is a 
$\mathcal {N} _{r} (\partial P)$-admissible function 
up to a multiplicative constant
independent of $k$. 
%We divide the proof in some steps.
Let $\gamma \in \mathcal {N} _r (\partial P)$. 

\medskip

\noindent \emph{Claim} : For every $z \in \gamma$ there exists
a ball $B_z \in \bigcup _{\ell = 0} ^{k} \mathcal B _{\ell} $ 
with $z \in B_z$ and such that the $\rho_k$-length 
of $\gamma \cap B_z$ is larger than $\tau (B_z)$.

\medskip

Indeed, assume at first that $z \in \gamma \cap \partial P$. Pick $B_z \in \mathcal B_k$ with $z \in B_z$. 
Since for $B\in \mathcal B_k$ we may assume that the function $\rho_B$ is larger than $\tau (B)$,
the ball $B_z$ admits the desired property. 

Now assume that $z \in \gamma \setminus \partial P$. 
Because  $r$ is small compared with $d_0$ the hypotheses 
of Lemma \ref{subcurve} are satisfied with $C = C_0$. Let $B \in \mathcal B$
as in this lemma. If $B$ belongs to 
$\bigcup _{\ell \ge k}  \mathcal B _{\ell} $, then choosing  $B_z \in \mathcal B_k$
with $B \subset B_z$ we are back to the previous case.

If $B \in \bigcup _{\ell = 0} ^{k-1} \mathcal B _{\ell} $,  
let $B_z = B$ and consider the subcurve $\gamma'$ given by 
Lemma \ref{subcurve} applied with $C=C_0$. It lies in the $(d/C_0)$-neighborhood
of $\partial P$, hence our hypotheses and Lemma \ref{C_0} assert that there
is no other limit set in $\mathcal L$ whose $(d/C_0)$-neighborhood
contains $\gamma'$. In addition we have
$$10^{-3}C_0 ^{-1} r(B_z) \le 10^{-3}C_0 ^{-1} 36 d \le d/(8C_0),$$
thus $\gamma'$ does not lie in any $(10^{-3}C_0 ^{-1}r(B_z))$-neighborhood
of a limit set of $\mathcal L$. Therefore the claim follows from
item (ii).

\medskip

According to $5r$-covering theorem
(see \cite{M}) we may extract from the collection $\{10 B_z ~;~ z \in \gamma\}$ a finite cover $\{10 B_1, ... , 10 B_n\}$
of $\gamma$ such that the balls $2B_1, ... , 2B_n$ are pairwise disjoint. This 
last property ensures that for every $v \in G_k ^0$
the number of elements of $\{B_1, ... , B_m\}$ which meet $v$
is bounded by above by a constant $C_1$ which depends only on
$\kappa$ and on the geometry of $Z$.  With the claim we obtain 
$$C_1 L_{\rho_k} (\gamma) \ge \sum _{i=1} ^n L_{\rho_k} (\gamma \cap B_i)
\ge \sum _{i=1} ^n \tau (B_i) .$$ 
Reordering if necessary 
we may assume that $\gamma$ enters successively the $10 B_i$'s. Let $w_i \in 
\partial P$ be the center of $B_i$. For every $i \in \{1, ... ,n-1\}$
the balls $10 B_i$ and $10 B_{i+1}$ intersect, so one of two balls
$20 B_i$ and $20 B_{i+1}$ contains the pair $\{w_i, w_{i+1}\}$. Thus 
$20 B_i \cap \partial P$ and $20 B_{i+1} \cap \partial P$
intersect.
Since the metric $\delta _{\partial P}$ and the rectriction of the $Z$-metric  are quasi-Moebius,
it follows from the triangle inequality that the last sum is linearly bounded from below 
by the $\delta _{\partial P}$-diameter of $\bigcup _{i=1} ^n 20 B_i \cap \partial P$.
Therefore $L_{\rho_k} (\gamma)$ is bounded from below in terms of $d_0$. 
\end{proof}

\section{Cannon Coxeter groups}\label{cannonsection}

Using the techniques developed in the previous section we prove a special case of the Cannon's conjecture
(Theorem \ref{cannon}). This result 
was already known, we review a proof due to M. Davis at the end of the section.
However our methods are new and of different flavour. 
An application to the hyperbolic Coxeter groups whose boundaries are 
homeomorphic to the 
Sierpinski carpet is presented in Corollary \ref{sierpinski} (the definition
of the Sierpinski carpet is recalled in Section \ref{menger}).

\begin{theorem} \label{cannon} Let $\Gamma$ be a hyperbolic Coxeter group
whose boundary is homeomorphic to the $2$-sphere. Then there is
a properly discontinuous, cocompact, and isometric action of
$\Gamma$   on $\hh ^3$,
the real hyperbolic space.
\end{theorem}

Assume $\Gamma$ is a hyperbolic Coxeter group
whose boundary is homeomorphic to the $2$-sphere. 
Let $\{ G_k \} _{k \in \hn}$ be a  $\kappa$-approximation of $ \partial \Gamma$,
and let $d_0$ be a fixed small positive constant. As before we denote by $\mathcal F _0$ the family of curves
$\gamma \subset \partial \Gamma$ with $\diam \gamma \ge d_0$.
The following proposition is the main step of the proof.

\begin{proposition} \label{3.6} 
For $p=2$ the modulus $\Mod _p (\mathcal F _0 , G_k)$ is bounded independently of the scale $k$.
\end{proposition}

Its proof will use the following general dichotomy :

\begin{lemma} Let $\Phi$ be a hyperbolic group whose boundary is connected.
Then either $\partial \Phi$ is a topological circle, or it contains a subset
homeomorphic to the capital letter  $\vdash \!\!\! \dashv$.
\end{lemma} 

\begin{proof}
Since $\partial \Phi$ is locally connected and nontrivial, it contains an
arc $X$ (see \cite{BK2}).  If $ x \in X$ is not an endpoint and if some 
neighborhood of $x$ in $\partial \Phi$  is contained in $X$, then 
--  since the action of $\Phi$ on $\partial \Phi$ is minimal
 -- $\partial \Phi$ 
is a 1-manifold
and so is homeomorphic to the circle.  Otherwise, if $x_1, x_2 \in X$
are distinct points which are not endpoints of $X$, then we may find
$y_1, y_2 \in \partial \Phi \setminus X$ and disjoint arcs $J_i \subset
\partial \Phi $, ($i = {1,2}$), 
joining $y_i$ to $x_i$.  Trimming the $J_i$'s and taking the union
$X \cup J_1 \cup J_2$, we get a subset homeomorphic to the capital letter 
$\vdash \!\!\! \dashv$, provided
$y_i$ is sufficiently close to $x_i$.
\end{proof}  

\begin{proof}[Proof of Proposition \ref{3.6}.]
With the notation of Section \ref{cox}, let $\mathcal L$ be the collection of circular 
limit sets.  We shall prove 
that for every $r>0$ the modulus $\Mod _2 (\mathcal {F} _{r} (\mathcal L), G_k)$
is bounded independently of the scale $k$.
Proposition \ref{3.6} will then follow from Corollary \ref{corcircular}. 
To establish the desired bound, we notice that  
$\mathcal {F} _{r} (\mathcal L)$ decomposes as a finite union of curve families 
$\mathcal F _{\delta _i, r _i} (\partial P _i)$, with $\partial P_i$
non-circular connected parabolic limit set and $\delta _i \le \delta _0 (r _i)$
--- see the discussion that follows Corollary \ref{comparable}. Therefore it is enough
to show
that for every $r> 0$, $\delta \le \delta _0 (r)$, and 
every non-circular connected parabolic limit set $\partial P
\subset \partial \Gamma$,
the modulus 
$\Mod _2 (\mathcal F _{\delta, r} (\partial P), G_k)$
is bounded independently of the scale $k$.

\smallskip

Let $\partial P$ be a non-circular connected parabolic limit set.
According to the above lemma, $\partial P$ contains a graph homeomorphic
to the capital letter $\vdash \!\!\! \dashv$ . 
Express the capital letter $\vdash \!\!\! \dashv$ as the union of two arcs $\alpha, \alpha '$
where $\alpha$ joins  the upper left
endpoint to the lower right endpoint, and 
$\alpha '$ joins the upper right
endpoint to the lower left endpoint.
Bending the two vertical segments of $\alpha '$ 
we obtain a horizontal segment denoted by $\beta$.
Consider a thin horizontal rectangle $\mathcal R$ 
which is a small planar neighborhood of the horizontal 
segment $\alpha \cap \beta$.  Say that a planar curve $\gamma$ 
\emph{cross-connects in
$\mathcal R$ two opposite sides of $\mathcal R$}, if 
$\gamma$ admits a subcurve contained in $\mathcal R$ and joining these sides.
Then every curve which lies 
within sufficiently small $C^0$ distance from $\alpha$ cross-connects
in $\mathcal R$ its horizontal sides, and every curve 
which lies 
within sufficiently small ${C} ^0$ distance from $\beta$ cross-connects
in $\mathcal R$ its vertical sides.

\smallskip

By assumption a similar picture appears in $\partial \Gamma$ and $\partial P$ too :
there exist two arcs $\eta_1, \eta_2$ in $\partial P$
and there exists a topological rectangle $\mathcal R \subset \partial \Gamma$ such that    
every curve in $\partial \Gamma$ which lies 
within sufficiently small ${C} ^0$ distance from $\eta_1$ cross-connects
in $\mathcal R$ its horizontal sides, and every curve 
which lies 
within sufficiently small ${C} ^0$ distance from $\eta_2$ cross-connects
in $\mathcal R$ its vertical sides. 
Let $\mathcal F _h (\mathcal R)$ (resp. $\mathcal F _v (\mathcal R)$)
be the family of curves contained in $\mathcal R$ and joining 
its horizontal (resp. vertical) sides.
With Proposition \ref{3.1} we get that for $\epsilon > 0$ small enough:
$$ \Mod _2 (\mathcal U _{\epsilon} (\eta _1 ), G_k ) \le 
\Mod _2 (\mathcal F _h (\mathcal R), G_k ),$$
$$ \textrm{and} ~ \Mod _2 (\mathcal U _{\epsilon} (\eta _2 ), G_k) \le 
\Mod _2 (\mathcal F _v (\mathcal R), G_k ).$$
The following lemma shows that 
$\min _{i = 1, 2} \Mod _2 (\mathcal U _{\epsilon} (\eta _i ), G_k )$
is bounded independently of $k$. Therefore Theorem \ref{3.4} 
applied to $\eta_1$ or $\eta_2$ shows that 
$\Mod _2 (\mathcal F _{\delta, r} (\partial P), G_k)$
is bounded independently of the scale $k$.
\end{proof}

\begin{lemma} \label{rect} 
There exists a constant $C \ge 1$ such that for any topological rectangle
$\mathcal R \subset \partial \Gamma$ one has for every $k \in \hn$ large enough :
$$\Mod _2 (\mathcal F _h (\mathcal R), G_k ) \cdot  
\Mod _2 (\mathcal F _v (\mathcal R), G_k ) \le C .$$
\end{lemma}

\begin{proof}
 A short and straightforword proof  
for annuli instead of rectangles can be found in  
\cite{H} preuve du lemme 2.14. The same proof applies \emph{verbatim} to
rectangles too. 

We now sketch a less direct but enlightening proof.  
Let $\mathrm{mod} (\cdot)$ denotes 
the classical analytic 
modulus on the Euclidean sphere $S^2$.  A well-known result asserts 
that  for $\mathcal R \subset S^2$ : 
$$\mathrm{mod} (\mathcal F _v (\mathcal R))\cdot  
\mathrm{mod} (\mathcal F _h (\mathcal R)) = 1.$$
The lemma follows from this fact and from some rather elaborate process to 
relate the combinatorial $2$-modulus on $\partial \Gamma$ 
with the analytic modulus on $S^2$ (see \cite{CFP} Th. 1.5 or \cite{BK1}
Cor. 8.8 for more details). 
\end{proof}   

\begin{proof}[Proof of Theorem \ref{cannon}. ]
According to Proposition \ref{3.6} and Corollary \ref{corunif}
one obtains that $\partial \Gamma$ is quasi-Moebius homeomorphic
to the Euclidean $2$-sphere. As explained in the 
sketch of proof in  Subsection \ref{presentation}, 
a theorem of Sullivan completes the proof
(\cite{S} p. 468, see also \cite{T}).
\end{proof}

\begin{corollary} \label{sierpinski} Let $\Gamma$ be a hyperbolic Coxeter group whose
boundary is homeomorphic to the Sierpinski carpet, then $\Gamma$
acts properly discontinuously by isometries on $\hh ^3$, and 
cocompactly on the convex-hull of its limit set.
\end{corollary}

\begin{proof} A nonseparating
topological circle in  $\partial \Gamma$ is called a \emph{peripheral circle}. 
A \emph{peripheral subgroup} of $\Gamma$ is the stabilizer of
a peripheral circle. In \cite{KaK} Th. 5, the following results are proved :
\begin{itemize}
 \item[(i)] there is only a finite number of conjugacy classes of peripheral subgroups ;
 \item[(ii)]  let $H_1, ... , H_k$ be a set of representatives of these classes, then the
double group $\Gamma \star _{H_i} \Gamma $ (\footnote{ $\Gamma
  \star _{H_i} \Gamma $ 
denotes the fundamental group of a graph of groups with two 
vertices and $k$ edges between them. The vertices groups are copies
of $\Gamma$ and the edges groups are the subgroups $H_i$.}) is a hyperbolic group 
with 2-sphere boundary. 
\end{itemize}
Observe that Theorem \ref{1.1} shows that every peripheral subgroup is a parabolic
subgroup. Choose subsets $I_1, ... , I_k$ of $S$ such that the subgroups $P_i$ 
generated by $I_i$ form a set of representatives of conjugacy classes of peripheral
subgroups. Then $\hat \Gamma := \Gamma \star _{P_i} \Gamma $ is obviously an index $2$ 
subgroup of a Coxeter group. Therefore the corollary
follows from item (ii) and Theorem \ref{cannon}.
\end{proof}

\subsection*{Another proof of Theorem \ref{cannon} } The following proof of Theorem \ref{cannon}
has been communicated to us by M. Davis. A theorem of Bestvina-Mess \cite{BeM} and the boundary
hypothesis show that $\Gamma$ is a virtual $3$-dimensional Poincar\'e 
duality group.
Then Theorem 10.9.2 of \cite{D} implies that  $\Gamma$ decomposes as 
$\Gamma = \Gamma _0 \times \Gamma _1$,
where $\Gamma _0$ is a finite Coxeter group and where $\Gamma _1$ is 
a Coxeter group
whose nerve is a $2$-sphere. Applying Andreev's theorem to the dual 
polyhedron to the nerve, 
one obtains that $\Gamma _1$ acts on $\hh ^3$ as a cocompact 
reflection group.

\section{The Combinatorial Loewner Property for Coxeter groups} \label{ex}

This section gives a sufficient condition for the boundary of a hyperbolic Coxeter group 
to satisfy the Combinatorial Loewner Property (Theorem \ref{CLP}).
Some examples of groups for which the condition applies are presented in the next section. 

\subsection{Generic curves} \label{gen}

Let $\Gamma$ be a hyperbolic Coxeter group whose boundary is connected.
Let  $Z$ be the  metric space $\partial \Gamma$ equipped with a self-similar metric (see 
Definition \ref{self}).
We fix in the sequel some $\kappa$-approximation 
$\{ G_k \} _{k \in \hn}$ of $Z$. 

\medskip

As before $d_0$ is a fixed small positive constant. For $r>0$  consider the family of curves 
$\gamma \subset Z$ with $\diam (\gamma) \ge d_0$,  and such that 
$\gamma \nsubseteq N_r (\partial Q)$ for any connected proper parabolic limit set $\partial Q \subset Z$. 
It is a non-empty curve family provided $r$ is small enough.
According to  Corollary \ref{comparable} its combinatorial modulus is comparable 
to the one of any of its subsets
of the form $\mathcal {U} _{\epsilon} (\eta)$.
In particular it does not depend on $r$ up to a multiplicative constant
independent of the scale.
We shall denote  such a family of curves  by $\mathcal {F} ^g $ and we shall call it
\emph{a family of generic curves}. Similarly we will 
call its combinatorial modulus the \emph{combinatorial modulus of generic curves}.

\medskip

When $p \ge 1$ is fixed we denote 
$\Mod _p (\mathcal {F} ^g , G_k)$ by $m_k$ for simplicity. 
It will be also convenient to set  $L_k := m_k ^{-1/p}$.  It satisfies :
$$L_k = \sup _{\rho} \frac{L_{\rho} (\mathcal {F} ^g)}{M_p (\rho) ^{1/p}}
~~ \textrm{with}~~ L_{\rho} (\mathcal {F} ^g) := \inf _{\gamma \in \mathcal {F} ^g} L_{\rho} (\gamma), $$
where the supremum is over all positive functions of $G_k ^0$.

\medskip

We now study the asymptotic behaviour of the sequence 
$\{L_k \}_{k \in \hn}$ depending on  $p \ge 1$.
Our main result establishes a weak type submultiplicative inequality 
for the sequence  $\{L_k \}_{k \in \hn}$ (Proposition \ref{super}).
The results of this paragraph must be compared with
those of paragraph \ref{hyp} concerning the family $\mathcal{F} _0$.

\medskip

For the rest of the paragraph $p$ is an arbitrary number in $[1, +\infty)$.
In the statement of the following two lemmata, $A_0$ denotes a fixed (large) positive number.

\begin{lemma} \label{gen1} There exist a constant $\Lambda \ge 1$ and a positive function $\Phi$
on $(0, +\infty)$ with the following properties. Let $B$ be a ball in $Z$ and let $k \in \hn$
such that the radius  $r(B)$ satisfies 
$$A_0 ^{-1} \le \frac{r(B)}{2 ^{-k}} \le A_0.$$ 
Consider two balls of same radius
$B_1, B_2 \subset B $,  
and let $t:= \frac{r(B_i)}{r(B)}$. Then for every $\ell \in \hn$
the $G_{k+\ell} $-combinatorial $p$-modulus of the family 
$$\{\gamma \in \mathcal F (B_1 , B_2) ~;~ \gamma \subset \Lambda B \}$$
is greater than 
$ m_{\ell} \cdot \Phi (t)$.
\end{lemma}

\begin{proof}
Using self-similarity property \ref{homothety} we can restrict ourself to the case 
$k=0$. 
For every $0<r \le \diam Z$ pick a maximal $\frac{r}{2}$-separated subset $E_r
\subset Z$. 
Since $Z$ is linearly connected there exists
a constant $\lambda >0$ such that every pair of points $x,y \in Z$ can be joined by a path of diameter
less than $\lambda d(x, y)$. For every pair of points $x,y \in E_r$ we choose
such a path and we call it $\eta_{x,y}$.

\smallskip

Let $\Lambda$ be a large number compared with $\lambda$. Consider any two $r$-balls $B_1, B_2$
contained in $B$.
Pick two points $z_1, z_2 \in E_r$ such that $B(z_i, r/2) \subset B_i$ for $i=1,2$. Then with the notation
of Section \ref{cox} we have
$$ \mathcal{U} _{r/2} (\eta_{z_1,z_2}) \subset 
\{\gamma \in \mathcal F (B_1 , B_2) ~;~ \gamma \subset \Lambda B \}.$$
Given $r>0$ there is only a finite number of left handside terms. Therefore 
the desired inequality follows from Theorem \ref{3.4} and from the fact that 
$r \asymp t$ by the rescaling assumption.
\end{proof} 

\begin{lemma} \label{gen2} 
There exist constants $\Lambda, D \ge 1$ and $b \in (0,1)$ with the following properties.
Let $B$ be a ball in $Z$ and let $k \in \hn$
such that the radius $r(B)$ satisfies 
$$A_0 ^{-1} \le \frac{r(B)}{2 ^{-k}} \le A_0.$$ 
Consider two continua $F_1, F_2 \subset Z$
with $F_i \cap \frac{1}{4} B \neq \emptyset$ and  $F_i \setminus B \neq \emptyset$.
Then for every $\ell \in \hn$ and every positive function $\rho$ on $G_{k+\ell} ^0$
there exists a path in $\Lambda B$ joining $F_1$ to $F_2$ whose $\rho$-length is smaller than
$$D \cdot M_p(\rho)^{1/p} \cdot \sum _{n=0} ^{\ell}  L_{\ell - n} ~ b^n.$$
\end{lemma}

\begin{proof}
The arguments are basically the same as those used in the proof of Proposition \ref{ball}.
Indeed pick $q \in \hn$ such that $\lambda := 2^{-q}$ and $a := 2 \cdot (8 \lambda)^{1/p}$ satisfy
$\lambda < 1/8$ and $a<1$. Let $\Lambda$ be as in Lemma \ref{gen1}.
Using Lemma \ref{gen1} we can construct -- like in the proof of
Proposition \ref{ball} -- a path joining $F_1$ to $F_2$ in $\Lambda B$ whose $\rho$-length is smaller than
$$D \cdot M_p(\rho)^{1/p} \cdot \sum _{n=0} ^{[\ell /q]}  L_{\ell - nq} ~ a^n ,$$
where $D >0$ depends only on $\Lambda, q, p, \kappa  $ and the geometry of $Z$, and where $[ \cdot ]$
denotes the integer part. Letting $b :=  a^{1/q}$ the lemma follows.
\end{proof}

With the above two lemmata we obtain :

\begin{proposition}\label{super}
There exist constants $b\in (0,1)$ and $C \ge 1$ such that for every pair of integers $k, \ell$
one has :
$$L_{k+\ell} \le C \cdot L_k \cdot \sum _{n=0} ^{\ell}  L_{\ell - n} ~ b^n.$$
Morever when $p$ belongs to a compact subset of $[1, +\infty)$ the constants $C$ and $b$
may be choosen independent of $p$.
\end{proposition}

\begin{proof}
Given a positive function $\rho _{k+\ell}$ of $G_{k+\ell} ^0$ we wish to construct
a curve $\gamma \in \mathcal{F} ^g $ with controlled $\rho _{k+\ell}$-length. 
For this purpose we pick two disjoint balls $E_1, E_2 \subset Z$ such that 
$\mathcal F (E_1, E_2) \subset \mathcal{F} ^g$, and we will look for $\gamma$ in
$\mathcal F (E_1, E_2)$. 

For any $v \in G_k ^0$, let $B_v $ be a ball containing $v$ and whose radius
is approximately $2^{-k}$.
Let $\rho _k : G_k ^0 \to \hr _+ $ be defined by 
$$ 
\rho _k (v) ^p = \sum _{w \cap 4\Lambda B_v \neq \emptyset} \rho _{k+\ell} (w) ^p ,
$$
where $\Lambda$ is the constant appearing in Lemma \ref{gen2}.
Since $Z$ is a doubling metric space one has 
\begin{equation}
\label{comparison2} 
M_p (\rho _{k+\ell}) \asymp M_p (\rho _{k}) \,.
\end{equation} 
>From Theorem \ref{3.4} there exists a curve $\delta \in \mathcal F (E_1, E_2)$
whose $\rho _k $-length is linearly bounded by above by
\begin{equation}
\label{length2} 
M_p(\rho_k)^{1/p} \cdot L_k  \,.
\end{equation}
Let $v_i \in G_k ^0$ so that $\delta$ enters successively $v_1, ..., v_N$,
and set $B_i := B_{v_i}$ for simplicity.
Then, providing $2^{-k}$ is small enough compared to $\diam E_i$, Lemma
\ref{gen2} allows one to construct by induction on $s \in \{1, ..., N\}$ a curve 
$\gamma _s \subset \cup _{i=1} ^s 4 \Lambda B_i$, joining $E_1$ to 
$ \cup _{i=s+1} ^N B_i \cup E_2$, whose
$\rho_{k+\ell}$-length is bounded linearly by above by 
$$
\sum _{i=1} ^s \rho _k (v_i)
\cdot \sum _{n=0} ^{\ell} L_{\ell - n} ~ b^n . 
$$ 
Therefore $\gamma := \gamma_N$
belongs to $\mathcal F (E_1, E_2)$ and its $\rho_{k+\ell}$-length is bounded linearly
by above by $L_{\rho _k} (\delta)
\cdot \sum _{n=0} ^{\ell} L_{\ell - n} ~ b^n $. 
Thanks to estimates (\ref{comparison2}) and (\ref{length2}), the proposition
follows. 
\end{proof}

Associated to generic curves, we introduce the following critical exponent 
$$ Q_m := \sup \{p \in [1, +\infty ) ~;~ \lim _{k \to +\infty} \Mod _p (\mathcal{F} ^g , G_k) = +\infty \}.$$
The previous theorem combined with the monotonicity and the continuity 
of the functions $p \mapsto \Mod _p (\mathcal{F} ^g , G_k)$ gives 

\begin{corollary} \label{Q_m}
The set of $p \ge 1$ such that $\lim _{k \to +\infty} \Mod _p (\mathcal{F} ^g , G_k) = +\infty$
is equal to the interval $[1, Q_m)$, in particular it is void if $Q_m =1$.
\end{corollary}

\begin{proof}  As a consequence of the weak submultiplicative inequality
stated in Proposition \ref{super} one has :
$$\lim _{k \to +\infty} L_k = 0 ~~\Longleftrightarrow~~ \exists \ell \in \hn ~~\textrm{with}~~
\sum _{n=0} ^{\ell}  L_{\ell - n} ~ b^n < C ^{-1} ,$$
where $C$ is the constant appearing in the weak submultiplicative inequality.
Hence the set of $p \ge 1$ such that $\lim _{k \to +\infty} \Mod _p (\mathcal{F} ^g , G_k) = +\infty$
is open in $[1, +\infty)$.
\end{proof}

\begin{proposition} \label{cut}
One has $Q_m >1$ if and only if $Z$ has no local cut point.
\end{proposition}

\begin{proof} The above corollary shows that 
$Q_m >1$ is equivalent to the fact that for $p=1$ one has 
$\lim _{k \to +\infty} \Mod _p (\mathcal{F} ^g , G_k) = +\infty $.

\smallskip

Assume $z_0$ is a local cut point of $Z$, and let $\eta \subset Z$ be a curve 
containing $z_0$ in its interior. Enlarging $\eta$ if necessary, we may assume that  for
$\epsilon > 0$ small enough, the curve family 
$\mathcal {U} _{\epsilon} (\eta)$ is contained in a family
$\mathcal{F} ^g$ of generic curves.
Therefore for any $p \in [1, +\infty)$ the moduli
$\Mod _p (\mathcal{F} ^g , G_k)$ and $\Mod _p (\mathcal {U} _{\epsilon} (\eta), G_k)$
are comparable. 
Since $z_0$ is a local cut point, every curve which belongs 
to $\mathcal {U} _{\epsilon} (\eta)$ passes through $z_0$, 
provided $\epsilon$ is small enough.
Choose for every $k \in \hn$ an element $v_k \in G_k ^0$ containing $z_0$.
The function of $G_k ^0$ whose value is  $1$ on $v_k$ and $0$ otherwise,
is a $\mathcal {U} _{\epsilon} (\eta)$-admissible function of $p$-mass $1$.
Hence the sequence $\{ \Mod _p (\mathcal{U} _{\epsilon} (\eta), G_k) \}_{k \in \hn}$
is bounded.

\smallskip

Conversely assume that $Z$ has no local cut point.
Let $\mathcal {F} ^g$ be a family of generic curves and let $\eta \subset Z$, $\epsilon > 0$ 
such that $\mathcal {U} _{\epsilon} (\eta) \subset  \mathcal{F} ^g $. 
Since $Z$ has no local cut point, a construction of 
J. Mackay \cite{Mac} shows that the family
$\mathcal {U} _{\epsilon} (\eta) $ contains an infinite collection of pairwise disjoint curves.
It implies obviously that 
$\lim _{k \to +\infty} \Mod _p (\mathcal {U} _{\epsilon} (\eta), G_k) = +\infty$
for $p=1$.
\end{proof}

\subsection{A sufficient condition for the CLP}

As in the previous subsection $\Gamma$ is a hyperbolic Coxeter group 
with connected boundary, and $\partial \Gamma$
is equipped with a self-similar metric.

Among the curve families already considered, the combinatorial modulus of generic curves is 
the lowest one, while the combinatorial modulus of $\mathcal{F} _0 $
is the largest one. Intermediate curve families
are the $\mathcal {F} _r (\mathcal L)$'s introduced in Section \ref{cox}.

In the following statement we allow $\mathcal L$ to be the empty collection
$\emptyset$. In this case we set $\mathcal {F} _r (\emptyset) := \mathcal{F} _0$
and $\Confdim (\emptyset) := 1$. The critical exponents $Q_M$ and $Q_m$ are respectively
defined in paragraphs \ref{hyp} and \ref{gen}.

\begin{theorem} \label{CLP}
Let $\mathcal{F} ^g$ be a family of generic curves and let $\mathcal L$
be a $\Gamma$-invariant collection of 
connected proper parabolic limit sets. 
Assume that for every pair of distinct elements $\partial P, \partial Q \in \mathcal L$, 
the set $\partial P \cap \partial Q$ is totally disconnected or empty.
Suppose that there exists $p \in [1, +\infty)$ satisfying 
$$p > \Confdim (\mathcal L) ~~ \textrm{and}~~ Q_m \le p \le Q_M,$$
such that for $r >0 $ small enough  and for every $k \in \hn$ :
$$ \Mod _p (\mathcal {F} _r (\mathcal L) , G_k) \le D  \Mod _p (\mathcal{F} ^g , G_k),$$
with $D = D(r) \ge 1$ independent of $k$. Then $p = Q_M$ and $\partial \Gamma$ 
satisfies the Combinatorial Loewner Property.
\end{theorem}

This  leads to

\begin{corollary} \label{simple}
Assume that $\Confdim (\partial \Gamma) > 1$, and suppose that for every 
proper, connected, parabolic limit set $\partial P \subset \partial
\Gamma$, one has 
$$\Confdim (\partial P) < \Confdim (\partial \Gamma)\, .$$
Suppose furthermore that for every pair $\partial P, \partial Q$ of
distinct, proper, connected, parabolic limit sets, the subset
$\partial P \cap \partial Q$ is totally disconnected or empty.
Then $\partial \Gamma$ satisfies the CLP. 
\end{corollary}

\begin{proof} One knows that $Q_M = \Confdim (\partial \Gamma)$
from Remark 2 at the end of Section
\ref{hyp}.
Therefore the statement follows from Theorem \ref{CLP} 
by considering the collection $\mathcal L$ of all proper
connected parabolic limit sets. 
\end{proof}

The rest of the paragraph is devoted to the proof of the theorem.
We will slightly abuse notation 
writing $L_k$ instead of  $C \cdot L_k $ where $C$ be the constant appearing in the statement 
of Proposition \ref{super}.  With this convention the weak submultiplicative inequality
in  Proposition \ref{super} simply writes 
$$L_{k+\ell} \le L_k \cdot \sum _{n=0} ^{\ell}  L_{\ell - n} ~ b^n.$$
Letting $u_{\ell} := \sum _{n=0} ^{\ell}  L_{\ell - n} ~ b^n$, it becomes
$L_{k+\ell} \le L_k \cdot u_{\ell}$.

\begin{lemma} \label{u}
Let $p \in [1, +\infty)$ and $a \in (0, 1)$ 
such that $\sum _{\ell = 0} ^{+\infty} u_{\ell} ^p ~ a^{\ell} = + \infty$.
Then the sequence $\{ L_k \} _{k \in \hn}$ tends to $+\infty$ exponentially fast.
\end{lemma}

\begin{proof}
Observe that  $\{ u_k \} _{k \in \hn}$ is submultiplicative. Indeed
from its definition and with the weak submultiplicative inequality one has 
$$u_{k + \ell} = \Big( \sum _{n = 0} ^{\ell -1} L_{k+\ell -n} ~ b^n \Big) +  b ^{\ell}u_k
~\le~  \Big( \sum _{n = 0} ^{\ell -1} L_{\ell -n} ~ b^n \Big)u_k +  b ^{\ell}u_k.$$
The weak submultiplicative inequality applied with $k = \ell = 0$ shows that $L_0 \ge 1$.
With the above inequalities we get that $u_{k + \ell} \le u_k \cdot u_{\ell}$.

\smallskip

Therefore there exists $\alpha \in \hr _+$ such that $\lim _{k \to +\infty} \frac{\log u_k}{k} = \log \alpha$.
With our hypotheses we obtain that $\alpha >1$. It follows that there exists a constant $C \ge 1$
such that for every $k \in \hn $ 
$$ C^{-1} \alpha ^{\frac{3}{4} k} \le u_k \le C \alpha ^{\frac{5}{4} k} .$$
In another hand one has for every $k, \ell \in \hn$ :
$$u_{k + \ell} = \Big( \sum _{n = 0} ^{\ell -1} L_{k+\ell -n} ~ b^n \Big) +  b ^{\ell}u_k
~\le~  L_k \Big( \sum _{n = 0} ^{\ell -1} u_{\ell -n} ~ b^n \Big) +  b ^{\ell}u_k.$$
Because $b \in (0,1)$ the above inequalities and a simple computation show the existence of constants
$C_1, C_2 \ge 1$ with
$$\alpha ^{\frac{3}{4} (k + \ell)} \le C_1 L_k \alpha ^{\frac{5}{4} \ell} + C_2 \alpha ^{\frac{5}{4} k} .$$
Letting $k = \ell$ yields the desired conclusion.
\end{proof}

\begin{proof}[Proof of Theorem \ref{CLP}.]
Let $p$ satisfy the hypotheses of the statement. 
At first we prove that the sequence  $\{L_k\} _{k \in \hn}$ is bounded.
Since $p \le Q_M$ part (ii) of Corollary \ref{Q_M}
shows that the sequence $\{\Mod _p (\mathcal{F} _0 , G_k)\} _{k \in \hn}$ admits a positive
lower bound. With our assumptions and with Corollary \ref{corcircular} we obtain
the existence of constants $C>0$ and $a \in (0, 1)$ such that for every $k \in \hn$
$$ C \le \sum _{\ell = 0} ^k m_{k - \ell} ~ a^{\ell} = \sum _{\ell = 0} ^k L_{k - \ell}^{-p} ~ a^{\ell}.$$
In particular $\{ L_k \} _{k \in \hn}$ does not tend to $+ \infty$, and so according to the previous lemma
the sum $\sum _{\ell = 0} ^{+\infty} u_{\ell} ^p ~ a^{\ell}$ is finite.
The weak submultiplicativity inequality writes $L_{k - \ell} \ge L_k \cdot u_{\ell} ^{-1}$,  and hence 
$$ C \le L_k ^{-p} \cdot \sum _{\ell = 0} ^k  u_{\ell} ^p ~ a^{\ell}.$$
Therefore $\{L_k\} _{k \in \hn}$ is bounded.

\smallskip

We now claim that $\{L_k\} _{k \in \hn}$ admits a positive lower bound. For this purpose observe that 
from its definition and with the weak submultiplicative inequality one has 
$$u_{k + \ell} = \Big( \sum _{n = 0} ^{\ell -1} L_{k+\ell -n} ~ b^n \Big) +  b ^{\ell}u_k
~\le~  L_k \Big( \sum _{n = 0} ^{\ell -1} u_{\ell -n} ~ b^n \Big) +  b ^{\ell}u_k.$$
Since $\{L_k\} _{k \in \hn}$ is bounded so is $\{u_k\} _{k \in \hn}$. Let $M$ be an upper
bound for $\{u_k\} _{k \in \hn}$, we obtain
$$u_{k + \ell} \le  L_k  S_{\ell}  +  b ^{\ell} M ~~\textrm{with}~~
S_{\ell} := \sum _{n = 0} ^{\ell -1} u_{\ell -n} ~ b^n.$$
Assume by contradiction that $\liminf _{k \to +\infty} L_k = 0$. At first choose
$\ell$ such that $ b ^{\ell}  M < 1/2$, then pick $k$ such that $ L_k  S_{\ell} < 1/2$ ;
we get that $u_{k + \ell} <1$. As seen already in the proof of Corollary \ref{Q_m}, this implies
that $\{L_k\} _{k \in \hn}$ tends to 0. Since $p \ge Q_m$ it contradicts 
Corollary
\ref{Q_m}, the claim follows.

\smallskip

The claim combined with our assumption imply that 
$\Mod _p (\mathcal {F} _r (\mathcal L) , G_k)$ is bounded independently of $k$. 
With Corollary \ref{corcircular} we get that the same holds for
$\Mod _p (\mathcal{F} _0 , G_k)$.
Summarizing there exits a constant $C>0$ such that for every 
$k \in \hn$ :
$$ C^{-1} \le \Mod _p (\mathcal F ^g , G _k) \le \Mod _p (\mathcal F _0 , G _k) \le C.$$
The $p$-combinatorial Loewner property 
comes now from Lemma \ref{gen1} and from  Propositions \ref{ball} and \ref{rd}. Note that 
Proposition \ref{rd} requires that $p >1$. This inequality is satisfied because by assumption
$p > \Confdim (\mathcal L) \ge 1$. Finally one has $p=Q_M$ according to 
Corollary \ref{Q_M}.
\end{proof}

\subsection*{Remarks and questions }   1) The equality $Q_M = Q_m$  is a
necessary condition for the CLP.  This follows from Remark 1 at the end of
Section \ref{hyp}, by choosing the continua $A, B$ in order that 
$\mathcal F (A, B)$ is contained in a family of generic curves.

\medskip

2) Recall that $Q_M$ is a quasi-Moebius  invariant of $Z$ 
(see Remark 2 at the end of Section \ref{hyp}). Is it true for $Q_m$ too ?

\section{Examples}\label{examples}

We give various examples of boundaries of hyperbolic Coxeter groups satisfying the Combinatorial
Loewner Property.

\subsection{Simplex groups } Let $(\Gamma, S)$ be an infinite hyperbolic
Coxeter group such that every proper parabolic subgroup is finite. This last condition
is equivalent to the fact that the Davis chamber of $(\Gamma, S)$ is a simplex 
of dimension $\vert S \vert - 1$ (see \cite{D}).
Using  Proposition \ref{cut} one sees easily that $Q_m >1$ if and only if $\vert S \vert \ge 4$.
Since the only infinite parabolic subgroup is $\Gamma$ itself, the assumptions
of Corollary \ref{simple} are clearly satisfied for $\vert S \vert \ge 4$. Therefore when $\vert S \vert \ge 4$,
the boundary of $\Gamma$ admits the CLP. 

This result was already known by other methods. Indeed it is a well-known theorem due to Lann\'er
that such a group acts by isometries, properly discontinuously and cocompactly
on the real hyperbolic space $\hh ^n$ of dimension $n = \vert S \vert -1$. 
Moreover the boundary of $\hh ^n$ is the Euclidean $(n-1)$-sphere
which is a Loewner space for $n-1 \ge 2$ (see \cite{HeK}).

\subsection{Prism groups }  Let $n \ge 3$,
consider an $n$-simplex and truncate an open neighborhood of one of its vertices. The
resulting
polytope is an $n$-prism (\emph{i.e.} isomorphic
to the product of a segment with a $(n-1)$-simplex). Let $\Sigma \subset \hh ^n$ 
be a geodesic $n$-prism whose dihedral angles are submultiples 
of $\pi$, and those of the truncated face are equal to $\frac{\pi}{2}$. 
Such polytopes have been classified by Kaplinskii
(\cite{Kap}, see also \cite{V} table 4), they exist only
when $n \le 5$. Let $\Gamma$ be the discrete subgroup of  $\Isom (\hh ^n)$
generated by reflections along the codimension $1$ faces of $\Sigma$ except
the truncated one.
The subset $\Gamma \cdot \Sigma$ is equal
to $\hh ^n$ minus a countable
union of disjoint totally geodesic half-spaces. Therefore $\Gamma$ is a (word)
hyperbolic Coxeter group and
its boundary identifies with $S^{n-1}$ minus a countable
union of disjoint Euclidean $(n-1)$-balls. 

Up to conjugacy the only infinite proper parabolic subgroup is the simplex group
generated by the faces adjacent to the truncated one. Its limit set
is an Euclidean $(n-2)$-sphere. 
Let $\mathcal L$ be the collection of all proper parabolic limit sets. 
It follows that $\Confdim (\mathcal L) = n-2$ and that  $\partial P
\cap \partial Q = \emptyset$ for every distinct elements $\partial P, 
\partial Q \in \mathcal L$.
Therefore the following lemma, in combination with Corollary \ref{simple} and the fact
that $\Confdim (\partial \Gamma) = Q_M$ (see Remark 2 at the end of Section
\ref{hyp}),
shows that $\partial \Gamma$ admits the CLP. 

\begin{lemma}
One has  $Q_M > n-2$.
\end{lemma}

\begin{proof} When  $n = 3$ this follows from Proposition \ref{cut} since in this case
$\partial \Gamma$ is homeomorphic to the Sierpinski carpet which doesn't possess
any local cut point.

In general we will prove that $\partial \Gamma$ equipped with the induced Euclidean metric 
satisfies $Q_m  > n-2$. 
Since $Q_m \le Q_M$ and since $Q_M$ is a quasi-Moebius 
invariant (see Remark 2 at the end of Section \ref{hyp}), the lemma will follow.
For this purpose we will exhibit a collection $\{F_i \} _{i \in \hn}$ of pairwise 
disjoint closed subsets of $\partial \Gamma$ with the following properties :
\begin{itemize}
\item[(1)] The $F_i$'s are uniformly bi-Lipschitz equivalent to the unit Euclidean sphere 
$S^{n-2}$,
 \emph{i.e.} for every $i \in \hn$
there is a $C$-bi-Lipschitz homeomorphism from
$F_i$ to $S^{n-2}$ where $C > 0$ is independent of $i$,
\item[(2)] There exists a constant $D$ with  $d_0 \le D \le \frac{1}{2} \diam(F_i)$, 
such that for every $i \in \hn$  one has 
$~\mathcal F _i := \{\gamma \subset F_i ~;~ \diam (\gamma) \ge D \} \subset \mathcal{F} ^g$. 
\end{itemize}
As we shall see, existence of such a collection implies that 
\begin{equation} \label{equaprism} 
\lim _{k \to +\infty} \Mod _{n-2} (\mathcal F ^g, G_k) = + \infty, 
\end{equation}
which in turn implies that  $Q_m > n-2$ thanks to Corollary \ref{Q_m}.

\smallskip

To establish the relation (\ref{equaprism}), notice first 
that $S^{n-2}$ admits the CLP with exponent
$n-2$ (see Theorem \ref{theoloewner}.1). Therefore item (1) implies that $\Mod _{n-2} (\mathcal F _i , G _k)$
is bounded away from $0$ independently of $(i, k) \in \hn ^2$. For every $k \in \hn$ pick a maximal
subcollection $\{F_{k_i}\}_{i = 1} ^{n_k} \subset \{F_{i}\}_{i \in \hn}$ such that every 
$v \in G_k ^0$ intersects at most one of the $F_{k_i}$'s. Then we get with item (2) :
$$\Mod _{n-2} (\mathcal F ^g , G _k) \gtrsim \sum _{i = 1} ^{n_k} \Mod _{n-2} (\mathcal F _{k_i} , G _k)
\gtrsim n_k .$$
Since $n_k$ tends to $+\infty$ the relation (\ref{equaprism}) follows.

\smallskip

We now construct the collection $\{F_{i}\}_{i \in \hn}$.  Pick a wall $M$ of $\Gamma$, its limit set in $\partial \hh ^n$ is contained in an $(n-2)$-sphere.  
Identifying $\partial \hh ^n$
with the unit $(n-1)$-Euclidean sphere, we may assume that $\partial M$ is contained
in the equator. Let $\mathcal B$ be the set of connected components of $\partial \hh ^n \setminus 
\partial \Gamma$. Because $\partial \Gamma$ is invariant by the reflection along $M$, 
every ball $B \in \mathcal B$ with $\overline {B} \cap \partial M \neq \emptyset$ is centered on the equator.
Thus by adding to $\partial M$ a countable union of  disjoint $(n-2)$-half-spheres contained in the north hemisphere,
we obtain a subset $F \subset \partial \Gamma$ homeomorphic to $S^{n-2}$.  
Since the parabolic subgroups are quasi-convex in $\Gamma$, there exists a constant $C >0$ such that for any distinct 
$B, B' \in \mathcal B$ one has 
$\Delta (B, B') \ge C$ (see \cite{B2} p. 89). One deduces easily from this fact that $F$ is bi-Lipschitz-equivalent 
to $S^{n-2}$.

\smallskip

We construct $F_i$ by induction, as a small deformation of $F$. Let $S_0$ be a smooth manifold diffeomorphic to 
$S^{n-2}$,
contained in the north hemisphere, disjoint from $F$ and lying within small Hausdorff distance from $F$.
For every $B \in \mathcal B$ with $B \cap S_0 \neq \emptyset$, replace $B \cap S_0$ by the spherical boundary
of the connected component of $B \setminus S_0$ which does not contain the center of $B$. Let $F_0$ be the resulting
subset of $\partial \Gamma$.  Since the balls $B \in \mathcal B$ are pairwise disjoint, $F_0$ is still a 
topological $(n-2)$-sphere disjoint from $F$.   Only finitely many  $B \in \mathcal B$ have
diameter larger than a given positive constant. Therefore by taking $S_0$ close enough to $F$,
one can ensure that $F_0$ lies 
within arbitrary small Hausdorff distance from $F$. Note that for $B \in \mathcal B$
the frontier of $B \cap S_0$ in $S_0$ belongs to $F_0$. Their union (for $B \in \mathcal B$)
is called the \emph{singular locus of $F_0$}. We claim that one can choose $S_0$ properly in order that 
$F_0$ is bi-Lipschitz equivalent to $S^{n-2}$ with
controlled Lipschitz constants. Indeed the bi-Lipschitz regularity holds  
provided the diedral angles of the singular locus in $F_0$
admit an uniform positive lower bound. Using the previous lower bound estimates 
for the relative distances between elements of $\mathcal B$, one can 
construct  $S_0$ in such a way that at every location its curvatures are smaller than those of the balls 
$B \in \mathcal B$ with $B \cap S_0
\neq \emptyset$. The control on the diedral angles follows, which in turn implies the claim.
We repeat this procedure starting with a manifold $S_1$ whose Hausdorff distance to $F$ is much smaller 
than the one between $F_0$ and $F$, in order that the resulting subset $F_1
\subset \partial \Gamma$ 
is disjoint from $F_0$. 
The construction of the $F_i$'s is now clear.
\end{proof}

\subsection{Highly symmetric Coxeter groups } 
Let $L$ be a finite graph whose girth is greater or equal to $4$
and such that the valency of each vertex is at least equal to $3$. 
Let $r$ be an even integer greater or equal to $6$. Consider 
the Coxeter system $(\Gamma, S)$ with one generating reflection
for every vertex in $L$, and such that the order of the product of
two distinct vertices $v, w$ is $r/2$ if $(v,w) \in L^1$, and $+\infty$
otherwise.
Its Davis complex is a $2$-cell contractible complex $X$ where
every $2$-cell is isomorphic to the regular right-angled $r$-gon
in $\hh ^2$, and where the link of every vertex is isomorphic to $L$. Equipped with
the geodesic distance induced by its  $2$-cells, $X$ is a CAT($-1$)-space
on which $\Gamma$ acts properly discontinuously, by isometries, and
cocompactly. In particular $\Gamma$ is a hyperbolic
group. Moreover the Cayley graph of $(\Gamma, S)$ identifies
with the  $1$-skeleton of $X$. The walls of $(\Gamma, S)$ identify with 
totally geodesic subtrees of $X$, namely those generated by a geodesic segment which joins
the middles of two opposite edges in a $2$-cell. We shall call them the \emph{walls of X}.
Consider $X$ minus the union of its walls, the closures of the connected
components 
are called the \emph{Davis chambers
of $X$}. Each of them is homeomorphic to the cone over $L$ (see \cite{D}, \cite{Ha} for more details).

Assume now that $L$ is a highly symmetric graph in the following sense: for every
pair of adjacent vertices $v, w \in L^0$ the pointwise stabilizer of the star of $v$
acts transitively on the set of remaining edges in the star of $w$. 
Recall that the star of a vertex is the union of its incident edges.
Examples of such graphs include the full bipartite graph with $s+t$ vertices
($s, t \ge 3$), the Mouffang generalized polygons (see \cite{R}), every $3$-transitive
trivalent graph (such as the Petersen graph), the odd graphs (see \cite{Bi}), \emph{etc}.

We have
 
\begin{proposition}\label{symmetric}
$\partial \Gamma$ satisfies the CLP.
\end{proposition}

When $L$ is a generalized polygon the statement was already  known. Indeed in this case $\partial \Gamma$
possesses a  self-similar metric $\delta$ such that $(\partial \Gamma, \delta)$ is an 
Loewner space \cite{BP1}.

\begin{proof}
Let $K \leqq \Isom (X)$ be the pointwise stabilizer of the Davis chamber containing the identity of $\Gamma$. 
Endow $\partial \Gamma$ with the self-similar metric $d$ induced by  a visual
metric on $\partial X$ ; then $K$ acts on $(\partial \Gamma, d)$ by
bi-Lipschitz homeomorphisms.

We will show that the assumptions of Theorem \ref{CLP} are satisfied with
$\mathcal L = \emptyset$. Since $\partial \Gamma$ does not admit any local cut point
-- it is homeomorphic to the Menger curve, see \cite{KaK} -- Proposition \ref{cut}
implies that $Q_M  > 1 = \Confdim (\mathcal L)$. Let $p > 1$
with $Q_m \le p \le Q_M$.

\smallskip

Let $I \subset S$ be such that $\partial \Gamma _I$ is a connected parabolic limit set.
We will prove by induction on $\vert S \vert - \vert I \vert$, that for $\delta, r >0$ small enough,
there exists a constant
$C \ge 1$ such that for every $k \in \hn$ one has
\begin{equation}\label{equasymmetric}
\Mod _p (\mathcal {F}_{\delta, r} (\partial \Gamma _I), G_k) \le C \Mod _p
(\mathcal {F}^g , G_k).
\end{equation}
For $\vert S \vert - \vert I \vert = 0$ there is nothing to prove.  Assume that
inequality (\ref{equasymmetric}) holds
for $\vert S \vert - \vert I \vert \le N$. We wish to prove it for $\vert S \vert - \vert I \vert \le N+1$.
Pick a curve $\gamma \subset 
\partial \Gamma _I$ which crosses every $\partial M_v$ for $v \in I$, and such that for $\delta, r, \epsilon$ small enough
$\mathcal {F}_{\delta, r} (\partial \Gamma _I)$  and
$\mathcal {U}_{\epsilon} (\gamma)$ have comparable modulus
(see Corollary \ref{comparable}).  Using the group
$K$ we will deform $\gamma$ into a curve contained in a larger parabolic limit
set.
 To do so we first need some information on the dynamics of the action of $K$ on
$\partial \Gamma$.

\smallskip

Let $v \in I$ and consider a special subgroup $\Gamma _J$ with $v \in J$. Its Davis complex,
denoted by $X_J$, embeds equivariantly and totally geodesically in $X$. The wall of $v$ in $X_J$
is equal to $ M_v \cap X_J$, where $ M_v$ denotes the wall of $v$ in $X$. Hence we see that 
either $\partial \Gamma _J \cap \partial M_v$ 
is of empty interior in $\partial M_v$, or $J$ contains all the neighbors of $v$ in $L$. 
For $v \in S$, let $U_v$ be the open dense subset of $\partial M_v$ whose complement is 
the union of all empty interior subsets of the form $\partial \Gamma _J
\cap \partial M_v$. 
One has :

\begin{lemma} \label{open} Every orbit of the $K$-action on $\Pi _{v \in S} \partial M_v$
is open.
\end{lemma}

Assuming the lemma we finish the proof of the proposition. Thanks to the lemma there exists   
$g \in K$ such that for every $v \in I$ one has $g \gamma \cap U_v \neq \emptyset$.  
Let $Q \leqq \Gamma$ be a parabolic subgroup such that $g \gamma \subset \partial Q$. Enlarging $\gamma$ if necessary
we may assume that its convex hull contains the identity of $\Gamma$. This property remains for $g \gamma$,
thus $Q$ is a special subgroup, $\Gamma _J$ say. Since $g \gamma$ crosses every $\partial M_v$ 
for $v \in I$, one has $I \subset J$ (see Remark 1 at the end of Section 
\ref{preliminary}). 
Therefore the  above discussion implies that $J$ contains the neighbors in $L$ of every
vertex $v \in I$. Hence, assuming  $\vert I \vert < \vert S \vert$, one obtains that 
$\vert I \vert < \vert J \vert$. Moreover, since $g$ is bi-Lipschitz,  
$\mathcal {U}_{\epsilon} (\gamma)$ and  $\mathcal {U}_{\epsilon} (g \gamma)$ have comparable modulus.
Therefore Corollary \ref{comparable} combined with the induction assumption
give the expected inequality (\ref{equasymmetric}).
\end{proof}  

It remains to give the 

\begin{proof} [Proof of Lemma \ref{open}. ]
Let $\Sigma \subset X$ be the Davis chamber containing the identity of $\Gamma$.
For a pair of adjacent vertices  $v, v' \in L^0$ consider the walls  $M, M' \subset X$ of the reflections $v, v'$.
They intersect exactly at the center $o \in X$ of the $2$-cell corresponding to the edge $(v, v') \in L^1$.
Denote by $H_{-}, H_{+}$ the closed half-spaces bounded by $M$ with $\Sigma \subset H_{-}$.
We will show that the pointwise stabilizer of $ H_{+}$ in $\Isom (X)$ acts transitively on the set of the edges of
$M' \cap  H_{-}$ adjacent to  the one containing $o$. Since $v, v'$ is an arbitrary pair of adjacent vertices, this result combined
with the $\Gamma$-action on $X$ gives the lemma.
To do so we shall use an argument of F. Haglund (\cite{Ha}, d\'emonstration de la Proposition A.3.1).
For simplicity assume first that $r$ is a multiple of $4$. Every automorphism of $L$ gives rise to a
group automorphism of $\Gamma$ which in turn induces an isometry of $X$. Let $a \in \Aut(L)$ 
which pointwise stabilizes
the star of $v$. Then the associated isometry $f$ acts trivially on the wall $M$ 
(the hypothesis $r$ is a multiple of $4$ is used at this point).
Therefore one can define $g \in \Isom (X)$ by letting $g = \textrm{id}$ on $H_{+}$ and $g = f$ on $H_{-}$.
Since $L$ is highly symmetric the result follows.

\smallskip

When $r$ is only a multiple of $2$, $f$ is not anymore the identity on $M$ but still it acts trivially 
on a tree $T$ containing $o$ and 
which zig-zag in the $2$-cells (joining the center to the middles of two edges at even distance). 
Moreover one can choose $T$ such that
$H_{+}$ and $M' \cap  H_{-}$ are contained in the closure of two distinct connected components of $X \setminus T$.
Hence the previous argument generalizes.
\end{proof}
   
\subsection{More planar examples }
Let $\Gamma$ be a hyperbolic Coxeter group with planar connected boundary.
Planarity is exploited to define trans\-versal intersection of parabolic limit sets.
Some modulus estimates are derived in Proposition \ref{transversal}. They are useful to establish the CLP
in some examples.

\begin{definition}
Two curves $\gamma _1, \gamma _2 \subset \partial \Gamma$ \emph{intersect tranversely} if there exists 
$\epsilon >0$ such that every pair of curves $\eta _1 \in \mathcal{U} _{\epsilon} (\gamma _1),
~\eta _2 \in \mathcal{U} _{\epsilon} (\gamma _2)$ intersect.  
We say that two connected parabolic limit sets $\partial P_1, \partial P_2 \subset \partial \Gamma$ \emph{intersects tranversely}
if there exists two curves  $\gamma _1 \subset \partial P_1, \gamma _2 \subset \partial P_2$ which
intersect transversely, and such that the smallest parabolic limit set containing $\gamma _i$
is $\partial P_i$ for $i=1,2$.
\end{definition}

\begin{proposition}\label{transversal} 
Let  $\partial P_1, \partial P_2$ be two connected parabolic limit sets which intersect transversely and 
let $\partial Q$ be the smallest parabolic limit set containing $\partial P_1 \cup \partial P_2$.
Then for every $p \ge 1$ and for $\delta, r$ small enough  one has 
$$ \inf _{i=1,2} \Mod _p (\mathcal {F}_{\delta, r} (\partial P_i), G_k) \asymp
\Mod _p (\mathcal {F}_{\delta, r} (\partial Q), G_k).$$
\end{proposition}

\begin{proof} Pick a pair of curves $\gamma _1 \subset \partial P_1, \gamma _2 \subset \partial P_2$ as in the 
definition.
Let $\epsilon >0$ be small enough and 
assume by contradiction that for $i=1,2$ the modulus  $\Mod _p (\mathcal {U}_{\epsilon} (\gamma _i), G_k)$
is large compared to $\Mod _p (\mathcal {F}_{\delta, r} (\partial Q), G_k)$.
Let $\rho : G_k \to \hr _+$ be a $\mathcal {F}_{\delta, r} (\partial Q)$-admissible 
function with minimal $p$-mass. 
Then for $i=1,2$ the family $\mathcal {U}_{\epsilon} (\gamma _i)$ contains a curve $\eta _i$ of small $\rho$-length
(see Lemma \ref{dual}).
By transversality the curves $\eta _1$ and $\eta _2$ intersect, thus $\mathcal {F}_{\delta, r} (\partial Q)$
contains a curve of small $\rho$-length, contradicting the admissibility of $\rho$.
Therefore  the moduli $\inf _{i=1,2} \Mod _p (\mathcal {U}_{\epsilon} (\gamma_i), G_k)$ and 
$\Mod _p (\mathcal {F}_{\delta, r} (\partial Q), G_k)$
are comparable, so the proposition follows from Corollary \ref{comparable}.
\end{proof}

As an illustation consider a 
$3$-dimensional cube and truncate an open neighborhood of every vertex. The
resulting polyhedron possesses $8$ triangular faces and $6$ hexagonal ones. 
Let $\Sigma \subset \hh ^3$ be a regular geodesic truncated cube whose diedral
angles are submultiples of $\pi$,
and those of the triangular faces are equal to $\frac{\pi}{2}$. By regular we
mean that $\Sigma$ admits all the cube symmetries.
Let $\Gamma$ be the discrete subgroup of $\Isom (\hh ^3)$ generated by 
reflections along the hexagonal faces of $\Sigma$.
The subset $\Gamma \cdot \Sigma$ is equal to $\hh ^3$ minus a countable disjoint
union of totally geodesic half-spaces. Therefore $\Gamma$ is a (word) 
hyperbolic Coxeter group whose boundary
is homeomorphic to the Sierpinski carpet. We will check that $\partial \Gamma$
satisfies the CLP.

For this purpose, equip $\partial \Gamma$ with a self-similar metric $d$ 
such that the symmetries of the cube act on $(\partial \Gamma, d)$ by bi-Lipschitz homeomorphisms.

Let $L$ be the graph whose vertices are the generators of $\Gamma$ and whose edges 
are the pairs $(s, s')$ with $s \neq s'$ and $ss'$ of finite order. Then $L$
is the $1$-skeleton of the octahedron. One sees easily
that there is only one type of proper parabolic subgroups with connected non-circular limit set: those whose 
graph is equal to $L$ minus a vertex and its adjacent edges. By applying an
appropriate cube symmetry
to such a parabolic $P$, we get another parabolic $P'$ such that $\partial P$ and $\partial P'$ intersect transversely.
Moreover the smallest parabolic limit set containing $\partial P \cup \partial P'$ is $\partial \Gamma$ itself. 
Hence Proposition \ref{transversal} and the
Lipschitz invariance of $d$, show that $\mathcal {F}_{\delta, r} (\partial P)$ and $\mathcal {F} ^g$ have  comparable 
modulus. The CLP comes now from Theorem \ref{CLP} with $\mathcal L$ equal
to the collection of all circular parabolic limit sets.

\subsection*{Remarks and questions }
1) In \cite{Ben} Y. Benoist exhibits examples of hyperbolic Coxeter groups with the following properties :
\begin{itemize}
\item Their Davis chambers are isomorphic to the product of two simplices of dimension 2,
\item They are virtually the fundamental group of a compact locally CAT($-1$) $4$-dimensional manifold,
\item They are not quasi-isometric to $\hh ^4$,
\item They admit a properly discontinuous cocompact projective action 
on a strictly convex open subset of the real projective space
$\hp ^4$.
\end{itemize}
An interesting question is to determine whether their boundary  satisfies 
the CLP or even the analytical Loewner
property. Previous examples of hyperbolic Coxeter groups enjoying 
the first three properties above appear
in Moussong's thesis (unpublished, see \cite{D} example 12.6.8). In 
connection with these problems, we remark that F. Esselmann has classified the
Coxeter polytopes in $\hh ^4$ isomorphic to the 
product of two simplices of dimension 2 (see \cite{Es}).

\medskip

2) Suppose $Z$ is a 
 Sierpinski carpet sitting in the standard $2$-sphere.
M. Bonk (\cite{Bonk} Prop. 7.6) has shown that if 
 the peripheral circles of $Z$ are uniformly quasi-Moebius
homeomorphic to the standard circle, and have pairwise relative distance 
bounded away from zero, then $Z$ satisfies a version
of the Loewner property for  transboundary $2$-modulus.
Since a Sierpinski carpet Coxeter group boundary $\partial\Gamma$
is quasi-Moebius
homeomorphic to such a $Z$, 
this should imply that $\partial\Gamma$
 satisfies a transboundary variant of the
CLP for  $2$-modulus.  It would be interesting
to have examples which do not satisfy the (usual) CLP.

\section{$\ell_p$-Equivalence relations} \label{equivalence}

This section covers some applications of the previous results to $\ell_p$-equivalence relations.
These equivalence relations are of great interest because of their invariance by quasi-Moebius
homeomorphisms. Moreover they can be used to provide examples of spaces which do not
admit the CLP (see the remark at the end of the section).

\medskip

We start  by defining the $\ell_p$-equivalence relation, 
it requires some preliminary materials. 
Let $Z$ be a compact doubling metric space. Assume in addition that it is \emph{uniformly perfect },
\emph{i.e.} there exists a constant
$0< \lambda <1$ such that for every ball $B(z, r)$ of $Z$ with $0 < r \le \diam Z$ one has
$ B(z,r) \setminus B(z, \lambda r) \neq \emptyset$. 

We will associate to $Z$ a Gromov hyperbolic graph $G$ such that 
$\partial G$ -- namely its boundary at infinity -- identifies canonically with $Z$.
For this purpose fix a constant $\kappa \geq 1$, and pick 
for each $k \in \hn ^*$ a finite covering
$\mathcal U _k$ of $Z$ with the following properties :
\begin{itemize} 
\item For every $v \in \mathcal U _k$ there exists $z_v \in Z$ such that :
$$B(z_v, \kappa ^{-1} 2^{-k}) \subset v \subset B(z_v, \kappa 2^{-k}),$$
\item For all distinct $v, w \in \mathcal U _k$ one has :
$$  B(z_v, \kappa ^{-1} 2^{-k}) \cap B(z_w, \kappa ^{-1} 2^{-k}) = \emptyset,$$
\item For every $z \in Z$ there exists $v \in \mathcal U _k$ with 
$B(z, \kappa ^{-1} 2^{-k}) \subset v$. 
\end{itemize}
Let $\mathcal U _0 = \{ Z \}$ be the trivial cover, and let $G$ be the 
graph whose vertex set is $\cup _{k \in \hn} \mathcal U _k$ and whose edges are
defined as follows : two distinct vertices $v$ and $w$ are joined by an edge if
\begin{itemize}
\item [--] $v$ and $w$ both belong to  $\mathcal U _k $, ($k \in \hn)$, and 
 $v \cap w \neq \emptyset $, or if
\item [--] one belongs to $\mathcal U _k $ the other one belongs to
$\mathcal U _{k+1} $ and  $v \cap w \neq \emptyset $.
\end{itemize}
Then $G$ is a Gromov hyperbolic graph with
bounded valency. The metric space $Z$ is bi-Lipschitz equivalent
to $\partial G$ equipped with a visual metric \cite{BP2}.

\medskip

For a countable set $E$ and for $p\in [1,\infty)$, we denote by $\ell_p(E)$
the Banach space of $p$-summable real functions on $E$.
The \emph{first $\ell_p $-cohomology group} of $G$ is
$$\ell_p H^1(G ) = 
\{f: G^0 \to \hr \; ; \; df\in \ell_p(G^1)\}/ \ell_p(G^0 ) + \hr ,$$
where $df$ is the 
function on $G^1$ defined by 
$$\forall a \in G^1, \; df(a) = f(a_+) - f(a_-),$$
and where $\hr$ denotes the set of constant fonctions on $G^0 $.
Equipped with the semi-norm induced by the $\ell_p$-norm of $df$
the topological vector space $\ell_pH^1(G )$ is a Banach space.
It is a quasi-isometric
invariant of $G$ and a quasi-Moebius invariant of $Z$. 
In addition $\ell_p H^1(G)$ injects 
in  $\ell_q H^1(G )$ for $1 \le p \le q < +\infty $.
See \cite{G}, \cite{BP2} for a proof of these results.
 
\medskip

Recall that $Z$ being a compact, doubling, uniformly perfect metric space, it admits 
a \emph{doubling measure}, that is a finite measure $\mu$ such that for every ball 
$B \subset Z$ of positive
radius one has : $0 < \mu (2B) \le C \mu (B)$, with $C \ge 0$ independent of $B$ (see \cite{He}).
By a result of R. Strichartz \cite{St}  
(see also \cite{Pa1}, \cite{BP2}) there is a continuous monomorphism
$$
\begin{array}{lll}
\ell_p H^1(G ) 
        & \hookrightarrow   & L^p(Z, \mu)/\hr 
\\
\ \ \ \ \ \ \    {[f]} & \longmapsto   &  f_\infty \mod \hr ,
\end{array}
$$
where  $f_\infty $ is defined  $\mu $-almost everywhere as follows :
for $z \in Z $ and for any geodesic ray
$r_z $ of $G $ with endpoint  $z $,
$$f_\infty (z ) = \lim_{t\to +\infty} f(r_z (t)).$$
Following  M. Gromov (\cite{G} p. 259, see also \cite{E} and \cite{B2})
we set
$$B^0_p (Z) :=
\{u : Z \to \hr \; \textrm{continuous} \; ; \;
 u = f_\infty \; \textrm{with}  \;
[f] \in \ell_p H^1(G)\},$$
and we define the \emph{$\ell_p $-equivalence relation} on $Z$ by ~:
$$x \sim_p y \Longleftrightarrow \forall u \in B^0_p(Z),
\ u(x) = u(y).$$
This is a closed equivalence relation which is invariant by the group of quasi-Moebius
homeomorphisms of $Z$. 

\begin{proposition}\label{2.1}
The cosets of the $\ell_p $-equivalence relation are continua.
\end{proposition}

\begin{proof}
We will use the following obvious properties of the space $B^0_p(Z)$:
\begin{itemize}
\item [(i)] if $u_1$ and $u_2$ belong to $B^0_p(Z) $ then $\max \{u_1, u_2 \}$ 
does too ;
\item [(ii)] let $u_1, u_2 \in B^0_p(Z) $, and let $U_1, U_2$ be open
subsets of $Z$ with $U_1 \cup U_2 = Z$.
Assume that $u_1 = u_2$ on $U_1 \cap U_2$.
Then the function $u$ defined by $u = u_1$ on $U_1$ and $u = u_2$
on $U_2$, belongs to $B^0_p(Z) $.
\end{itemize}
First we claim that for any coset $F$ and any compact subset $K \subset Z$ 
disjoint from $F$, there exists $u \in B^0_p(Z) $
such that $u(F) = 0$ and $u(K) = 1$. Indeed for any $z \in K$
there exists a function $u_z \in  B^0_p(Z) $ with $u_z (F) = 0$ and
$u_z (z) > 1$. Extract a finite cover of $K$ from the 
open subsets $\{u_z > 1 \}$.  Let $U_1,...,U_n$ be such a cover and let
$u_1,...,u_n$ be the corresponding functions. Then 
the function 
$$ v = \sum_{i=1}^n \max \{0, u_i \}$$
belongs to $B^0_p(Z)$ and satisfies $v(F) = 0$ and $v \geq 1$ on $K$.
Letting $u = \min \{1, v\}$ the claim follows.

\smallskip

Assume now by contradiction that a coset $F$ 
is a disjoint union of two non-empty compact subsets
$F_1$ and $F_2$. Let $r > 0$. By the previous claim
there exists a function $u \in B^0_p(Z) $ such that $u(F)= 0$ and 
$u (Z \setminus N_{r} (F)) = 1$.
For $r$ small enough $N_{r} (F)$
is the disjoint union of $N_{r} (F_1)$ and $N_{r} (F_2)$.
Define a function  $v$ on $Z$ by letting
$v(z) = 1$ for $z \in N_{r} (F_2)$, and $v(z) = u(z)$
otherwise. Then the above property (ii) applied to
the open subsets $Z \setminus \overline {N} _{r} (F_i)$ 
shows that the function
$v$ belongs to $B^0_p(Z) $. Moreover  it satisfies  $v(F_1) = 0$
and $v(F_2) = 1$. Contradiction.
\end{proof} 

The following result relates  
the $\ell_p$-equivalence relation with the combinatorial $p$-modulus. 

\begin{proposition}\label{3.3} 
Let $\{ G_k \} _{k \in \hn}$ be a  $\kappa$-approximation of $Z$.
Assume $x, y \in Z$ satisfy $x \nsim _p y$, then 
there exist open subsets  $U, V \subset Z$ containing respectively $x$ and $y$ such 
that 
$$ \lim_ {k \rightarrow + \infty} \Mod _p (U, V, G_k) = 0.$$
\end{proposition}

\begin{proof} 
Recall that $Z$ being a
doubling metric space, up to a multiplicative constant the 
$G_k$-combinatorial $p$-modulus does not depend
on the $\kappa$-approximation (see Proposition \ref{3.2}).   
Consider the graph $G$ associated to a family
of covers $\{ \mathcal U _k \} _{k \in \hn} $ as described at the 
beginning of the section.
For $k \in \hn$, let $G _k$ be the subgraph of $G$ which is the 
incidence graph of the 
covering $\mathcal U _k$. 
The family $\{ G_k \} _{k \in \hn}$
is a $\kappa$-approximation of $Z$. 
In addition $G_k ^0$ identifies with the sphere in $G$
of radius $k$ centered at the unique vertex of $G_0$. 

\smallskip

Let $u \in B_p ^0 (Z) $ with $u(x) \neq u(y)$. Changing $u$
to $au +b$ ($a, b \in \hr$) if necessary,  we can assume that there exist open
subsets $U, V$ of $Z$ with $x \in U, y \in V$ and $u \le 0$ on $U$,
$u \ge 1$ on $V$.
Pick a function $f : G ^0 \to \hr$ such that $df \in \ell _p (G ^1)$
and $f _{\infty} = u$. Choosing  $f$ properly one can ensure that 
for $k$ large enough and for every $v \in G_k ^0$ the following holds :
$$f(v) \le 1/3 ~~\textrm{if}~~ v \cap U \neq \emptyset , ~~
 f(v) \ge 2/3 ~~\textrm{if}~~ v \cap V \neq \emptyset ,$$
(see \cite{BP2} preuve du Th. 3.4).
Let $\rho _k : G _k ^0 \to \hr _+ $ be defined by 
$$ \rho _k (v) = 3 \max _{(v, w) \in G_k ^1} \vert f(v) - f(w) \vert. $$
Obviously it is an $\mathcal F (U, V)$-admissible function.
In addition its $p$-mass satisfies
$$ M_p (\rho _k ) \le 3 ^p  \sum _{a \in G _k ^1} \vert df (a) \vert ^p 
.$$
For $df \in \ell _p (G ^1)$ the right handside term tends to $0$
when $k$ tends to $+ \infty$.
\end{proof}

Combining several previous results we now collect some applications 
to hyperbolic Coxeter groups.

\begin{corollary} \label{coset}
Let $\Gamma$ be a hyperbolic Coxeter group. 
Then each coset of the $\ell_p $-equivalence relation on $\partial \Gamma$
is either a point or a connected parabolic limit set.
\end{corollary}

\begin{proof}
It is a straitforward consequence of Proposition \ref{2.1} and 
Corollary \ref{1.3}.
\end{proof}
 
\begin{corollary} \label{Confdim}
Let $\Gamma$ be a hyperbolic Coxeter group with connected boundary. Let $p \ge 1$
and suppose that $ \sim _p$ admits a coset different from a point and the whole $\partial \Gamma$.
Then $ \sim _p$ admits a coset $F$ with $\Confdim (F) = \Confdim (\partial \Gamma)$.
\end{corollary}

\begin{proof} Let $p$ be as in the statement and let $\mathcal L$ be the collection of the cosets 
of $ \sim _p$ which are different from a point.
 From the above corollary its elements are connected proper parabolic limit sets. 
Obviously the $\Gamma$-invariance and separation hypotheses
of Corollary \ref{corcircular} are satisfied. 
In addition, for every $r>0$, Proposition \ref{3.3} shows that 
$\Mod _p (\mathcal F _{r} (\mathcal L), G_k)$ tends to $0$ when $k$ tends to $+ \infty$.
The same holds for every $q \ge p$ by monotonicity of the function
$q \mapsto \Mod _q $. 
Assume by contradiction that $\Confdim (\mathcal L) < \Confdim (\partial \Gamma)$.
Since $\Confdim (\partial \Gamma) = Q_M$ (see the remark 2 at the end of Section
\ref{hyp}), one can apply Corollary  \ref{corcircular} with exponent $Q_M$.
One obtains that  $\Mod _{Q_M} (\mathcal F _{0} , G_k)$ tends to $0$ when $k$ tends to $+ \infty$,
contradicting Corollary \ref{Q_M}(ii).
\end{proof}
 
Finally we return to general approximately self-similar metric spaces and 
to the combinatorial Loewner property. 
The second part of the following corollary is the combinatorial analog of Theorem 0.3  in \cite{BP2}.

\begin{corollary} \label{confdim} Let $Z$ be an approximately self-similar
  metric space. Assume $Z$ is connected and let $p \ge 1$.
Then $p > Q_M$  if and only if  $(Z/ \sim _p ) = Z$.
If in addition  $Z$ satisfies the CLP, then for $1\le p \le Q_M$ the quotient  $Z/ \sim _p$ is a singleton.
\end{corollary}

\begin{proof} The second part of the statement 
follows from Proposition \ref{3.3} and from the monotonicity of $p \mapsto \Mod _p$. 
To establish the ``only if'' part of the first, one invokes that $Q_M$ is equal to the Ahlfors regular
conformal dimension of $Z$ (see the remark 2 at the end of Section
\ref{hyp}), and the fact that  
$B_p ^0 (Z)$ separates the points of $Z$  for 
$p$ strictly larger than the Ahlfors-regular conformal dimension \cite{BP2}. 
Conversely if $(Z/ \sim _p ) = Z$, then Proposition \ref{3.3} implies that 
$\Mod _p (\mathcal {F}_0 , G_k)$ tends to $0$ when $k$ tends to $+\infty$.
With Proposition \ref{Q_M}(ii) we get that $p > Q_M$. 
\end{proof}

\subsection*{Remarks and questions }
Corollary \ref{confdim} may be used to produce examples of spaces for which the CLP fails. 
Indeed suppose that for a given hyperbolic group $\Gamma$, the family
of quotient spaces $\partial \Gamma /\!\sim _p $, ($p \in [1, +\infty)$),
contains intermediate states between the singleton and the whole $\partial \Gamma$.
Then, according to Corollary \ref{confdim},  $\partial \Gamma$ does not
admit the CLP.

Currently all known examples of hyperbolic groups for which the CLP fails are of two
types.  Either their boundaries 
admit local cut points -- in which case the CLP fails ``trivially'' (see Proposition \ref{linear}.2) -- 
or  they decompose as $\Gamma = A \star _C B $ 
and there exists $p \in [1, Q_M]$ with : 
$$ x \sim_p y \Longleftrightarrow x = y ~~\textrm{or}~~\exists g \in \Gamma ~~\textrm{such that}~~
x, y \in g \partial B.$$
Examples of the second type, including some Coxeter groups, are described in \cite{B2}.
In the Coxeter group case, Corollary \ref{Confdim} shows that 
$\Confdim (\partial B) = \Confdim (\partial \Gamma)$. 

It would be desirable to have a better understanding of the relations between the combinatorial
modulus and the $\ell_p$-cohomology.

%%%%%%%%%%%%%%%%%%%%%%%%%%%%%%%%%%%%%%%%%%%%%%%%%%%%%
%\selectlanguage{frenchb} 

\bibliographystyle{alpha}
\bibliography{coxeter}

\medskip

\medskip

\noindent Marc Bourdon, Universit\'e Lille 1, D\'epartement de math\'ematiques, Bat. M2, 59655 Villeneuve d'Ascq,
France, bourdon@math.univ-lille1.fr 

\medskip

\noindent Bruce Kleiner, Courant Institute, New York University,
251 Mercer Street,
New York, N.Y. 10012-1185, USA, 
bkleiner@courant.nyu.edu 

\end{document}